\documentclass[preprint,10pt]{elsarticle_nojournal}
\usepackage{amsmath,amsfonts,amssymb}

\newcommand{\coefmat}{\mathbf{A}}           
\newcommand{\coefev}{\lambda}               
\newcommand{\evmin}{\coefev_\mathrm{min}}   
\newcommand{\evmax}{\coefev_\mathrm{max}}   

\newcommand{\dom}{\Omega}
\newcommand{\domint}{\dom_{\mathrm{int}}}
\newcommand{\dompml}{\dom_{\mathrm{pml}}}
\newcommand{\domext}{\dompml^\infty}

\newcommand{\bdpml}{\Sigma}

\newcommand{\bx}{\mathbf{x}}
\newcommand{\bq}{\mathbf{q}}

\newcommand{\by}{\mathbf{y}}
\newcommand{\bs}{\mathbf{s}}
\newcommand{\bp}{\mathbf{p}}

\newcommand{\bB}{\mathbf{B}}
\newcommand{\bA}{\mathbf{A}}
\newcommand{\bR}{\mathbf{R}}

\newcommand{\radpml}{R_{\mathrm{pml}}}  
\newcommand{\lpml}{L}                   
\newcommand{\spml}{\sigma}              
\newcommand{\shpml}{\tilde\sigma}         
\newcommand{\rpml}{r_{\sigma}}            
\newcommand{\dtpml}{\tilde{d}_\spml}           %
\newcommand{\dpml}{d_\spml}
\newcommand{\hpml}{h_{\spml}}                   %
\newcommand{\Apml}{\coefmat_{\sigma}}     
\newcommand{\Apmld}{\coefmat_{\sigma,\delta}}
\newcommand{\sfpml}{a_{\sigma}}                  
\newcommand{\sfpmld}{a_{\sigma,\delta}}
\newcommand{\oppml}{\mathcal{A}_{\sigma}}    
\newcommand{\oppmld}{\mathcal{A}_{\sigma,\delta}}                     
\newcommand{\vecxs}{\vecx_\sigma} 
\newcommand{\bc}{\boldsymbol{c}}
\newcommand{\bz}{\boldsymbol{z}}

\newcommand{\bbeta}{\boldsymbol{\beta}}
\newcommand{\ppml}{\mathbf{p}_{\mathrm{pml}}}  

\newcommand{\rad}[1]{{#1}_{\operatorname{rad}}}
\DeclareMathOperator{\diag}{diag}
\DeclareMathOperator{\sign}{sign}
\DeclareMathOperator{\supp}{supp}
\DeclareMathOperator{\spa}{span}
\renewcommand{\div}{\operatorname{div}}

\newcommand{\testf}[1]{{#1}^\dagger}
\newcommand{\vecs}[1]{\mathbf{#1}} 
\newcommand{\innerprod}[2]{\left\langle#1,#2\right\rangle}

\newcommand{\vecx}{\vecs x}

\newcommand{\vecxi}{\boldsymbol{\xi}}

\newcommand{\setR}{\mathbb{R}}
\newcommand{\setRp}{\setR^+}
\newcommand{\setRpz}{\setR^+_0}

\newcommand{\setRmz}{\setR^-_0}
\newcommand{\setN}{\mathbb{N}}
\newcommand{\setC}{\mathbb{C}}
\newcommand{\setCp}{\mathbb{C}^+}
\renewcommand{\Re}{\operatorname{Re}}
\renewcommand{\Im}{\operatorname{Im}}
\newcommand{\Px}{\Pi_{\parallel}}
\newcommand{\Pxp}{\Pi_{\perp}}

\newcommand{\radsub}{\parallel}
\newcommand{\tangsub}{\perp}

\newcommand{\Js}{\mathbf{J}_{\sigma}}
\newcommand{\Jsg}{\mathbf{J}_{\sigma,\gamma}}

\newcommand{\iepar}{\eta} 

\newcommand{\scalsp}{\mathcal V}
\newcommand{\vecsp}{\boldsymbol{\mathcal Q}}
\newcommand{\scalspint}{\scalsp_{\mathrm{int}}}
\newcommand{\vecspint}{\vecsp_{\mathrm{int}}}
\newcommand{\vecspext}{\vecsp_{\mathrm{ext}}}
\newcommand{\scalspext}{\scalsp_{\mathrm{ext}}}
\newcommand{\scalspextrad}{\scalsp_{\mathrm{rad}}}
\newcommand{\scalspexttang}{\scalsp_{\Sigma}}
\newcommand{\vecspextrad}{\mathcal Q_{\mathrm{rad}}}
\newcommand{\vecspexttang}{\vecsp_{\Sigma}}

\usepackage{natbib}
\usepackage{placeins}
\usepackage{amsmath,amsthm,amssymb,amsfonts}
\usepackage{marginnote}
\usepackage{bbm}
\usepackage{fullpage}
\usepackage[font=footnotesize,labelfont=bf]{caption}
\usepackage[font=footnotesize,labelfont=bf]{subcaption}
\usepackage[]{xcolor}
\usepackage{hyperref}
\usepackage{enumitem}
\usepackage{tabularx}
\usepackage{graphicx}
\graphicspath{{./images/}}
\theoremstyle{plain}
\newtheorem{thm}{Theorem}[section]
\newtheorem{lem}[thm]{Lemma}
\newtheorem{cor}[thm]{Corollary}
\newtheorem{prop}[thm]{Proposition}

\theoremstyle{definition}

\newtheorem{assump}[thm]{Assumption}

\theoremstyle{remark}
\newtheorem{rmk}[thm]{Remark}

\journal{CMAME}
\def\appendixname{}

\begin{document}

	\begin{frontmatter}
		
		
		
		\title{Radial perfectly matched layers and infinite elements for the anisotropic wave equation}
		
		\author[adr1]{Martin Halla}
		\ead{m.halla@math.uni-goettingen.de}
		\affiliation[adr1]{organization={Institut f\unexpanded{\"u}r Numerische und Angewandte Mathematik, Georg-August-Universit\unexpanded{\"a}t G\unexpanded{\"o}ttingen},
			addressline={Lotzestr.\ 16-18}, 
                        city={G\unexpanded{\"o}ttingen},
			postcode={37083}, 
			state={Niedersachsen},
			country={Germany}}
		\author[adr2]{Maryna Kachanovska}
		\ead{maryna.kachanovska@inria.fr}
		\affiliation[adr2]{organization={POEMS, ENSTA Paris, CNRS, Inria, Institut Polytechnique de Paris},
			addressline={828 Boulevard des Mar\unexpanded{\'e}chaux}, 
			city={Palaiseau},
			postcode={91120}, 
			country={France}}
		\author[adr3]{Markus Wess}
		\ead{markus.wess@tuwien.ac.at}
		\affiliation[adr3]{organization={Institute of Analysis and Scientific Computing, TU Wien},
			addressline={Wiedner Hauptstra\ss e 8-10}, 
			postcode={1040}, 
			state={Vienna},
			country={Austria}}
		
		\begin{abstract}
				We consider the scalar anisotropic wave equation.
				Recently a convergence analysis for radial perfectly matched layers (PML) in the frequency domain was reported and in the present article we continue this approach into the time domain.
				First we explain why there is a good hope that radial complex scalings can overcome the instabilities of PML methods caused by anisotropic materials.
				Next we discuss some sensitive details, which seem like a paradox at the first glance:
				if the absorbing layer and the inhomogeneities are sufficiently separated, then the solution is indeed stable.
				However, for more general data the problem becomes unstable. 
				In numerical computations we observe instabilities regardless of the position of the inhomogeneities, although the instabilities arise only for fine enough discretizations.
				As a remedy we propose a complex frequency shifted scaling and discretizations by Hardy space infinite elements or truncation-free PMLs. We show numerical experiments which confirm the stability and convergence of these methods. 
		\end{abstract}

		
	\end{frontmatter}
	


\section{Introduction}
In the 1990s B\'erenger \cite{MR1294924} introduced his perfectly matched layer method as an approximate transparent boundary condition for transient electromagnetic wave equations.
Soon it was recognized \cite{ChewWeedon:94} that the PML equations can be derived by means of a complex scaling technique, which was already extensively used under the names complex scaling/analytic dilation/spectral deformation since the 1970s in mathematical physics for analysis and resonance computations, see \cite{aguilar_combes,balslev_combes,simon_complex_scaling}, and a detailed review in \cite{HislopSigal:96}.
We refer to the introduction of \cite{Halla:19Diss} for an extensive literature review on PML and to the articles  \cite{joly:12,pled_desceliers_22} for an introduction and review of the PMLs. 
No doubt, the reason for the popularity of the PML method stems from the fact that it is very easy to implement: for transient equations one only needs to introduce some additional auxiliary unknowns (without any knowledge of a fundamental solution or Dirichlet-to-Neumann operator).
For this reason it is easily possible to apply the PML method to all kinds of equations.
However, as soon as one deviates from classical applications it is a delicate question if the PML method yields physically correct and stable solutions.
In particular, backward waves which can occur in dispersive materials \cite{BecacheJolyVinoles,BecacheKachanovska:17} and waveguide geometries \cite{BonnetBDChambeyronLegendre:14,DuruKreiss:14b} lead to challenges for the PML.
Another important challenge, which is the focus of this article, are anisotropic materials.
Indeed also the application/construction of absorbing boundary conditions for anisotropic equations requires a careful analysis and has received quite some attention, see e.g.,\ \cite{BecacheGivoliHagstrom:10,SavadattiGuddati:10b,SavadattiGuddati:12b,RabinovichEtal:19,Lee:21}.

Up to our knowledge, the instability of B\'erenger's Cartesian perfectly matched layers for anisotropic media was noticed as early as in 1996 by Fang Q. Hu \cite{FangQHu96} when applying PMLs to Euler equations, and, as a first remedy, the author suggested numerical filtering. The reason for the  appearance of such instabilities was investigated in particular in \cite{hesthaven98,AbarbanelGottliebHesthaven99}, and at first was attributed to the possible loss of the strong well-posedness of the PML problem. In \cite{FangQHu01} it was shown that a strongly well-posed PML system can be constructed, but it is still unstable: the PML instability is induced by so-called backward propagating waves in the direction of the PML absorption. While the explanation in \cite{FangQHu01} was given using somewhat semi-heuristic arguments, it found its mathematical justification in the seminal work by B\'ecache, Fauqueux and Joly \cite{BFJ}. There, it was proven that the classical Cartesian B\'erenger's PMLs in one direction with constant absorption can be unstable when applied to anisotropic media. The analysis in \cite{BFJ} is based on a plane-wave analysis of the Cauchy problem for the constant-coefficient PML media filling the free space. It is proven that some (but not all) of the instabilities of the PMLs are high-frequency instabilities, caused by so-called backward propagating waves in the direction of the PMLs, occurring in particular in many anisotropic materials. 

At the numerical level, the fact that instabilities appear for high frequencies implies that coarse discretizations are not likely to exhibit them, see \cite{KreissDuru:13}. Since the usual goal is to construct numerical methods for which the error can be made arbitrarily small, this is not a very satisfactory solution (though it indeed can be used in practice). 
While there exist several works discussing and predicting the behaviour of  Cartesian B\'erenger's PMLs
 in anisotropic media (whereas the treatment against the instabilities is typically problem-dependent),
there exists only little knowledge about the behaviour of the other types of PMLs, in particular radial PMLs, introduced in \cite{MR1638033}. The analysis of \cite{BFJ} does not apply to this case, and the numerical experiments are not always conclusive.

The principal goal of this work is to address this question in detail. More precisely, we would like to answer the following questions:
\begin{enumerate}[label={Q\arabic*. }]
	\item Are radial PMLs stable when applied to anisotropic problems?
	\item If not, what is the reason for the instability?
	\item If the answer to Question 1 is negative, can we find a workaround, \textbf{which would not be specific to the anisotropic \textit{scalar} problem in question}?
\end{enumerate}
To answer these questions we concentrate on the simplest model problem, namely the anisotropic wave equation, for which many computations can be done explicitly. This model already contains some of the difficulties which occur when applying PMLs to more challenging problems (e.g., anisotropic elasticity). 

\textbf{Answers to Questions 1 and 2. }Our answers to Questions 1 and 2 are given in Section \ref{sec:well_posedness_and_stability}.
For theoretical purposes, we work with radial perfectly matched layers of \textit{infinite length} (which we will refer to as \textit{truncation-free} PMLs), which already exhibit many interesting phenomena observed in truncated PMLs. 
If the spatial support of the source term (compactly supported in space and time) is located far enough from the damping layer, then the solutions to the PML problem are stable, i.e. exhibit at most time-polynomial growth. 
On the other hand, the PML system itself allows for unstable solutions. We believe that this is related to the presence of an essential spectrum of the underlying operator in the right-half of the complex plane. We prove the existence of such a spectrum. Our numerical experiments (Section \ref{sec:two_different_types_behaviour}) indicate that the instabilities manifest themselves at the discrete level, when the discretizations are chosen fine enough, independently of the support of the source term.
Our findings for the radial PMLs complement those by K. Duru and G. Kreiss \cite{DuruKreiss:12, KreissDuru:13} who observed similar phenomena for Cartesian PMLs applied to anisotropic elasticity and wave equations.

\textbf{One answer to Question 3. } Indeed, stable Cartesian PMLs for the anisotropic wave equation were constructed by a clever change of variables in \cite{MR3138002}.
Nonetheless, because our goal is to propose a method that \textbf{is possibly} suitable for other anisotropic problems apart from the scalar anisotropic wave equation, we choose a different path. 
Our starting point is the PML time-harmonic work by one of the authors of the present article \cite{Halla:22PMLani}, where it is suggested to use complex frequency shifted PMLs. However, we show that, like in the time-harmonic case, in the time domain the complex frequency shift has to be chosen proportionally to the damping parameter of the layer, with the proportionality coefficient limited 
by the anisotropy. In this setting increasing the damping parameter does not allow to decrease the PML error, and the convergence can be ensured by increasing the width of the layer only. See Section \ref{sec:remedy} for details. 
We propose two solutions to this problem in Section \ref{sec:numerics}, and test their performance numerically. They are outlined below. 
\begin{enumerate}[label=S\arabic*.]
  \item We do not truncate the perfectly matched layer, but rather use the complex-shifted change of variables combined with the discretization by infinite element methods. To do so we use the Hardy space infinite elements (HSIE), see \cite{hsm,Halla:16} for their introduction, analysis and numerical experiments in the time-harmonic case, and  \cite{RuprechtSchaedleSchmidt:13}, \cite[Chapter~10.2.4]{WessDiss}, \cite{td_ie} for the time-domain formulation and experiments. In particular, we use a time-domain equivalent of the two-scale method presented in \cite{HallaNannen:18}.
  \item We exploit so-called exact PML methods, where the exactness is ensured either by choosing a non-integrable damping parameter \cite{BermudezHervellaNPrietoRodriguez:08}, or a coordinate transformation mapping an infinite domain to a finite one \cite{HugoninLalanne:05,yang_wang_gao} and \cite[Chapter~4.5.1]{Halla:19Diss}.
		We concentrate on this latter choice. 
\end{enumerate}

Let us finally remark that the complex frequency shift appeared already in the literature of PMLs to prevent long time instabilities of Cartesian PMLs \cite{BecachePetropoulosGedney:04}, or as a practical approach to stabilize perfectly matched layers in anisotropic elasticity, cf. \cite{DuruKreiss:12}. Studies combining classical, Cartesian PMLs, and infinite finite elements, can be found in \cite{inf_finite_pettigrew}.

The remainder of this article is structured as follows.
In Section \ref{sec:setting} we specify the wave propagation problem under consideration and the associated PML equations derived by Cartesian and radial complex scalings.
We recap the instability results of \cite{BFJ} for Cartesian PMLs and give some heuristic arguments why radial PMLs might allow to overcome these instabilities.
We present different computational radial PML examples, some with stable and some with unstable behaviour.
In Section \ref{sec:well_posedness_and_stability} we investigate the cause of the instabilities: the continuous problem is well-posed and for convenient configurations the solutions are stable.
Nevertheless for small spectral parameters $s$ with $\Re s>0,$ there exists an essential spectrum, which proves the unstable character of the PML equations.
As a remedy we introduce in Section \ref{sec:remedy} a complex frequency shifted Hardy space infinite element method, as well as a corresponding truncation-free PML. We report computational results, which confirm stability and convergence of these methods in Section \ref{sec:numerics}.

\section{Problem setting and motivation}\label{sec:setting}
\subsection{Problem setting}
\subsubsection{The model and the goals of the article}
We consider wave propagation in anisotropic media, described by a scalar anisotropic wave equation. Namely, we look for $u\colon\setRpz\times \setR^2\rightarrow \setR$ and $\bp\colon\setRpz\times \setR^2\rightarrow \setR^2$, where $\setRpz:=\{t\in\setR:t\geq 0\}$, s.t.
\begin{align}
\partial_t u &= \div\bp + f \quad\text{in }\setRp\times \setR^2,\qquad
\bA^{-1}\partial_t \bp = \nabla u  \quad\text{in }\setRp\times \setR^2,\\
u(0,\bx)&=0, \quad \bp(0, \bx)=0 \text{ in }\setR^2,
\end{align}
where the source $f$ is sufficiently regular, $f(0,\cdot)=0$, and $f(t,\cdot)$ is compactly supported inside a bounded domain $\Omega_{f}$ for all $t>0$. The matrix $\bA\in \setR^{2\times 2}$ is symmetric strictly positive definite.
By eliminating $\bp$ we write this system as a second order equation, and equip it with homogeneous initial conditions:
\begin{align}
\label{eq:main_problem}
&\partial_t^2 u-\operatorname{div}(\bA\nabla u)=\partial_tf \quad\text{in }\setRp\times \setR^2,\qquad
u(0,\bx)=0, \quad \partial_t u(0, \bx)=0 \quad\text{in }\setR^2.
\end{align}
We are interested in the solution to the above problem inside a bounded convex domain $\domint$, s.t.\ $\Omega_{f}\subset \domint$. To bound the computational domain we surround $\domint$ by an absorbing perfectly matched media $\domext:=\setR^2\setminus \overline{\domint}$, which is usually then truncated to a bounded layer $\dompml$. Inside this layer, the original PDE is modified in such a way that, on one hand, its solution decays exponentially fast in space, and, on the other hand, the waves propagating from $\domint$ into $\dompml$ do not get reflected at the interface $\bdpml$ between $\domint$ and $\dompml$.  
The absorbing system is constructed via a certain frequency-dependent change of variables applied in the direction normal to $\bdpml$ (this will be discussed in detail later). In the literature there are several common choices of the geometry of $\domint$ and the construction of the corresponding PMLs:
\begin{enumerate}
\item Cartesian PMLs \cite{MR1294924, MR1412240} where $\domint$ is chosen as a box, e.g., $(-a,a)^2$, and $\dompml$ as $(-a-L, a+L)^2\setminus [-a, a]^2$;
\item radial PMLs (see \cite{MR1638033}), where $\domint$ is a circle, e.g., $B_{\radpml}:=\{\bx\in\setR^2:\|\bx\|<\radpml\}$, and $\dompml=B_{\radpml+L}\setminus \overline{B_{\radpml}}$;
\item general convex PMLs (see e.g., \cite{MR1862449} or  \cite{MR1638033}), where $\domint$ is a general convex domain. 
\end{enumerate}
For anisotropic wave propagation problems, and in particular for the problem \eqref{eq:main_problem}, in \cite{BFJ} it has been shown that the Cartesian PMLs may exhibit time-domain instabilities, i.e., exponentially growing in time solutions (see Section \ref{sec:instability_cartesian_pml} for a more detailed discussion). The respective analysis is done in a simplified setting (by using a plane-wave analysis), and  cannot be applied to radial PMLs or general non-polygonal convex PMLs. 

As discussed before, the goal of this article is to examine the question of stability of radial PMLs for anisotropic problems. Because in its full generality this problem is quite technically challenging, we concentrate on the simplest model \eqref{eq:main_problem}. We first introduce the Cartesian and radial PMLs and next discuss the failure of Cartesian PMLs for anisotropic problems. We finally explain in Section \ref{sec:radial_pml} why the use of the  radial PMLs may seem, at least at a first glance, to lead to a stable
time-domain system, contrary to the Cartesian PMLs. 

\subsubsection{Perfectly Matched Layers}
\label{sec:pml_intro}
In this section we provide a very brief introduction to PMLs, first concentrating on the classical, Cartesian PMLs, and subsequently discussing the radial PMLs.
\paragraph{Cartesian PMLs} To construct a PML system with the PML applied in the $x$-direction, for $|x|>a$, we start by transforming the original problem \eqref{eq:main_problem} into the Laplace domain. For $v\in L^1(\setR)$ s.t.\ $v(t)=0$ for $t<0$ (i.e., $v$ is a causal function) the Fourier-Laplace ( Laplace) transform is defined by 
\begin{align*}
\hat{v}(s):=(\mathcal{L}v)(s):=\int_{-\infty}^{+\infty}\mathrm{e}^{-st}v(t)dt=\int_{0}^{+\infty}\mathrm{e}^{-st}v(t)dt, 
\; s \in \setCp:=\{s\in\mathbb{C}\colon \Re s>0\}.
\end{align*}
For functions $v$ defined on $\mathbb{R}^+$ only, we define the Fourier-Laplace transform by extending $v$ by $0$ to $\mathbb{R}^{-}$. 

This yields the following reformulation of \eqref{eq:main_problem}: 
\begin{align*}
	s^2 \hat{u}-\operatorname{div}(\bA \nabla \hat{u})=s\hat{f} \quad\text{in }\setR^2.
\end{align*}
Next one assumes that $\hat{u}(s, x, y)$ is analytic in $x\in \setCp$ and thus satisfies the same relation as above but with a complexified spatial variable:
\begin{align}
	x_{\sigma}:=
	\left\{
	\begin{array}{ll}
	x+s^{-1}\int_{a}^x\sigma(x')dx', & x>a, \\
          x, & x\in [-a, a],\\
	x+s^{-1}\int_{-a}^x\sigma(x')dx', & x<-a. 
	\end{array}
	\right.
\end{align}
Here the absorption parameter $\sigma$ is a non-negative $L^{\infty}$-function, s.t.\ $\sigma=0$ on $(-a, a)$. 
We then conclude that $\hat u^{\sigma}(s,x,y):=u(s,x_{\sigma}, y)$ satisfies the following PDE: 
\begin{align*}
s(s+\sigma)\hat{u}^{\sigma}-\operatorname{div}\left( \left(1+\frac{\sigma}{s}\right)\Js^{-1}\bA \Js^{-\top}\nabla \hat{u}^\sigma\right)=s\hat{f} \quad\text{in }\setR^2, 
\qquad\Js:=D_{\bx} \vecxs=\operatorname{diag}\left(1+\frac{\sigma}{s},1\right).
\end{align*}
Next the computational domain is truncated to $(-a-L, a+L)\times \setR$ and homogeneous boundary (e.g., Dirichlet or Neumann) conditions are imposed on the artificial boundary $|x|=a+L$. 
To simulate the free-space problem, the analogous change of variables and domain truncation is applied in the $y$-direction.
It remains to reformulate the problem in the time domain, by performing the inverse Laplace transform and introducing auxiliary variables where necessary to obtain a first or second order system in time. We will not present this step here, and the interested reader can consult, e.g., \cite{TeixeiraChew:97,duru_kreiss_efficient,grote2010efficient}, for various reformulations of the PMLs for the wave-type problems. 
\paragraph{Radial PMLs} They are based on the same idea as the Cartesian PMLs, with the difference being that the first step is rewriting the respective system in the radial coordinates $(r, \phi)\in \setRpz\times [0,2\pi)$. 
When appropriate, we consider these coordinates as functions of $\bx$, i.e., $r=r(\bx)=\|\bx\|$ and $\phi=\phi(\bx)$.
Denoting by $\nabla_{r,\phi}v=(\partial_r, r^{-1}\partial_{\phi})^\top v$, $\operatorname{div}_{r, \phi}\boldsymbol{v}=r^{-1}(\partial_r (rv_{r})+\partial_{\phi}v_{\phi})$, we obtain the following Laplace-domain radial coordinates counterpart of \eqref{eq:main_problem}:
\begin{align*}
	s^2\hat{u}-\operatorname{div}_{r, \phi}(\bA^{\phi}\nabla_{r, \phi}\hat{u})=s\hat{f} \quad\text{in } \setR_0^+\times [0, 2\pi),
\end{align*}
with $\bA^{\phi}:=\bR_{\phi}^\top\bA \bR_{\phi}$ and the rotation matrix $\bR_{\phi}:=\begin{pmatrix} \cos\phi&-\sin\phi\\\sin\phi&\cos\phi\end{pmatrix}$. Remark that here we used the same notation  $\hat{u}$ for the unknown in the Cartesian and polar coordinates.
The radial PMLs are based on a change of variables  
\begin{align}
  \rpml(s,r):=
	\left\{
	\begin{array}{ll}
		r+s^{-1}\int_{\radpml}^r\sigma(r')dr', & r>\radpml, \\
		r, & r\leq \radpml.
	\end{array}
	\right.
        \label{eq:radpml}
\end{align}
Here $\sigma$ satisfies the following assumption. 
\begin{assump}
	\label{assum:sigma}
	Let $\sigma\in L^{\infty}(\setR^+)$ satisfy $\sigma(r)=0$ for $r\leq\radpml$ and $\sigma(r)>0$ for $r>\radpml$.
\end{assump}
Let us additionally introduce 
\begin{align}
  \begin{aligned}
  \dtpml(s,r)&:=\frac{\rpml(s,r)}{r}=1+s^{-1}\tilde{\sigma}(r),&\tilde{\sigma}(r)&:=r^{-1} \int_{\radpml}^r\sigma(r')\,dr',\\
  \dpml(s,r)&:=\frac{d\rpml(s,r)}{dr}=1+s^{-1}\sigma(r), & \vecxs(s,\bx)&:=\bx\dtpml(s,\|\bx\|).
  \end{aligned}
  \label{eq:pmldefns}
\end{align}
Where convenient, we will abuse the notation as follows: we will write $\dpml(s)$ for $\dpml(s,r)$, when $\dpml$ is considered as a function of the Laplace variable $s$ and $r$ is fixed, or, resp., $\dpml(r)$. The same applies to $\dtpml$ and $\vecxs$ respectively.
In addition, we use the overloaded notation $\dtpml(s,\bx):=\dtpml(s,\|\bx\|)$, $\sigma(\bx):=\sigma(\|\bx\|)$, etc.
With this new change of variables and the above notation,  we obtain the following problem in the Laplace domain, satisfied by
$\hat{u}^{\sigma}(s,\bx):=\hat{u}(s,\vecxs)$:
\begin{align}
s^2\dtpml \dpml \hat{u}^{\sigma}-\operatorname{div}_{r, \phi}(\bA_{\sigma}^{\phi}\nabla_{r, \phi}\hat{u}^{\sigma})= s\hat{f} \quad\text{in } \setR_0^+\times [0, 2\pi),
\qquad\bA_{\sigma}^{\phi}:=\operatorname{diag}(\dtpml, \dpml)\bA^{\phi}\operatorname{diag}(\dpml^{-1}, \dtpml^{-1}). 
\label{eq:rad_lpl_v0}
\end{align}
Let us translate the above in the time domain, and put it in a form more suitable for the numerical implementation, based on a reformulation of \eqref{eq:rad_lpl_v0} in the Cartesian coordinates. Given $\bx\in \setR^2\setminus\{0\}$, we denote by $\Pi_{\parallel}(\bx)$ and $\Pi_{\perp}(\bx)$ the orthogonal projections on the spaces $\operatorname{span}\{\bx\}$ and $\left(\operatorname{span}\{\bx\}\right)_{\perp}$ respectively.
Then the Jacobian can be expressed as
\begin{align}
  	\label{eq:js}
	\begin{split}
\Js&:=D_{\bx}\bx_{\sigma}={\dtpml}\operatorname{Id}+\frac{\bx\, \bx^{\top}}{\|\bx\|}\partial_r{\dtpml}=\dtpml(\Px+\Pxp)+\|\bx\|\Px \partial_r\dtpml\\
&=(\dtpml+\|\bx\|\partial_r\dtpml)\Px+\dtpml\Pxp=\dpml\Px+\dtpml\Pxp.
\end{split}
\end{align} 
Its determinant is $\det \Js=\dpml \dtpml$.
Then \eqref{eq:rad_lpl_v0} rewrites, based on the usual Piola transform \cite[pp. 59-60]{boffi},
\begin{align}
	\label{eq:basigma_main}
  s^2 \operatorname{det}\Js\hat{u}^{\sigma}-\operatorname{div}\left(\bA_{\sigma}\nabla \hat{u}^{\sigma}\right)=s\hat f,\quad \bA_{\sigma}(s,\bx):=\Js (s,\bx)^{-1}\bA\Js(s,\bx)^{-\top}\det\Js(s,\bx).
\end{align}
The above can be written as a first-order system:
\begin{align}
	\label{eq:ford}
  s \operatorname{det}\Js\hat{u}^{\sigma}&=\operatorname{div}\hat\bp^{\sigma}+\hat f,\qquad
  s\bA_{\sigma}^{-1}\hat{\bp}^{\sigma}=\nabla\hat u^{\sigma}.
\end{align}
This rewriting preserves the skew symmetry of the first-order counterpart of our original problem \eqref{eq:main_problem}.
Using
\begin{align}
  \begin{split}
  s\dpml(s)\tilde \dpml(s) = s+ \sigma+\tilde\sigma + \frac{\sigma\tilde\sigma}{s},\qquad
  s\frac{\dpml(s)}{\tilde \dpml(s)} &= s + \sigma-\tilde\sigma -\frac{(\sigma-\tilde\sigma)\tilde\sigma}{s+\tilde \sigma},\\
  s\frac{\tilde \dpml(s)}{\dpml(s)} &= s - (\sigma-\tilde\sigma) +\frac{(\sigma-\tilde\sigma)\sigma}{s+ \sigma},
  \end{split}
  \label{eq:sdd}
\end{align}
we obtain that the left-hand-sides of \eqref{eq:ford} can be rewritten as
\begin{align}
  \label{eq:detJ_s}
  s \operatorname{det}\Js\hat{u}^{\sigma}={}&s\hat u^\sigma +\left(\sigma+\tilde\sigma+\frac{\sigma\tilde\sigma}{s}\right)\hat u^\sigma,\\
  \label{eq:detJJJ_s}
  \begin{split}
  s\det\Js^{-1}\Js^{\top}\bA^{-1}\Js\hat\bp^{\sigma} ={}&s\bA^{-1}\hat \bp^{\sigma}+\left(\sigma-\tilde\sigma -\frac{(\sigma-\tilde\sigma)\tilde\sigma}{s+\tilde \sigma}\right)\Px \bA^{-1}\Px\hat\bp^{\sigma},\\
  &+\left(- (\sigma-\tilde\sigma) +\frac{(\sigma-\tilde\sigma)\sigma}{s+ \sigma}\right)\Pxp \bA^{-1}\Pxp\hat\bp^{\sigma}.
  \end{split}
\end{align}
Thus, introducing the auxiliary unknowns
\begin{align*}
  s\hat v &:= \hat u^\sigma,\\
  \hat\bq &:= \frac{\tilde\sigma}{s+\tilde \sigma}\Px\hat\bp^{\sigma} +\frac{\sigma}{s+ \sigma}\Pxp \hat\bp^{\sigma}\quad \iff\quad s\hat{\bq}=\tilde{\sigma}\left(1-\frac{\tilde{\sigma}}{s+\tilde{\sigma}}\right)\Px\hat{\bp}^{\sigma}+\sigma\left(1-\frac{\sigma}{s+\sigma}\right)\Pxp\hat{\bp}^{\sigma},
\end{align*} 
and applying the inverse Laplace transform (where the respective time-domain functions are indicated by omitting the hats) leads to the first order system 
\begin{align}
  \begin{split}
  \partial_t u^\sigma &=-(\sigma+\tilde\sigma) u^\sigma - \sigma\tilde\sigma  v+\div\bp^{\sigma}+f,\\
  \partial_t v &= u^\sigma, \\
  \bA^{-1} \partial_t\bp^{\sigma}&=(\sigma-\tilde\sigma)\left(\Pxp \bA^{-1}\Pxp\bp^{\sigma}-\Px \bA^{-1}\Px\bp^{\sigma}-\Pxp \bA^{-1}\Pxp\bq+\Px \bA^{-1}\Px\bq\right) +\nabla u^\sigma,\\
  \partial_t \bq&=\tilde\sigma\left(\Px\bp^{\sigma}-\Px\bq\right)+\sigma\left(\Pxp\bp^{\sigma}-\Pxp \bq\right).
  \end{split}
  \label{eq:rad_time}
\end{align}
In practice, the radial PMLs are truncated just like the Cartesian PMLs, by imposing homogeneous boundary conditions at $r=\radpml+L$.
\subsection{Instability of Cartesian perfectly matched layers}
\label{sec:instability_cartesian_pml}	
	A convenient way to start the stability analysis of the initial-value problem for the PML system is to consider the case of frozen coefficients, i.e.,\ when $\sigma=\operatorname{const}$, in the free space setting. Under this simplified assumption, the analysis can be performed with a help of a plane wave (Fourier) approach \cite{kreiss}.
The plane-wave analysis consists of studying plane-wave solutions $e^{i\omega t-i\mathbf{k}\cdot\mathbf{x}}$ of the underlying problem; in particular, the problem admits such non-trivial solutions if $(\omega, \mathbf{k})$ satisfy the so-called dispersion relation $F(\omega, \mathbf{k})=0$. In our case \eqref{eq:main_problem}, the dispersion relation reads 
	\begin{align*}
          F(\omega,\mathbf{k})=\omega^2-\mathbf{k}^{\top} \mathbf{A}\mathbf{k}=0.
	\end{align*}	
	 This dispersion relation defines multiple branches of the solution $\omega_j(\mathbf{k})$, $j=1, \ldots, K$, $\mathbf{k}\in \setR^2$ (modes). The plane wave analysis relies on studies of $\Im\omega_j(\mathbf{k})$: in particular, if we assume that non-vanishing $\omega_j(\mathbf{k})$ satisfy $\omega_j(\mathbf{k})\neq \omega_{\ell}(\mathbf{k})$ for all $\mathbf{k}\neq 0$, the underlying initial-value problem is stable if and only if  $\Im\omega_j(\mathbf{k})\geq 0$ for all $\mathbf{k}\in \setR^2$.
This approach was developed in \cite{BFJ}, where stability of Cartesian constant coefficient PMLs for anisotropic systems was analyzed, based on a perturbation argument applied to the PML dispersion relation. It was shown that if backward propagating waves (defined below) are present, the resulting PML system is unstable. 
More precisely, the dispersion relation for the original system allows to define for each mode $\omega_j(\mathbf{k})$ two quantities: 
\begin{align*}
\text{the group velocity }
\mathbf{v}_{g, j}( \mathbf{k}):=\nabla_\mathbf{k}\omega_j(\mathbf{k})
\qquad\text{ and the phase velocity }
\mathbf{v}_{p, j}(\mathbf{k}):=\frac{\omega_j(\mathbf{k})}{|\mathbf{k}|^2}\mathbf{k}.
\end{align*}
Let $\mathbf{e}\in\setR^2$ be a unit vector. Then one of the results of \cite{BFJ} reads:  if there exists $j$ and a direction $\mathbf{k}$ such that $(\mathbf{v}_{g, j}( \mathbf{k})\cdot\mathbf{e})(\mathbf{v}_{p, j}(\mathbf{k})\cdot\mathbf{e})<0$, then the PML in direction $\mathbf{e}$ is unstable.
The above phenomenon is referred to as the existence of a backward wave in the direction $\mathbf{e}$ and strongly depends on the chosen direction $\mathbf{e}$.

\begin{figure}
	\centering
	\includegraphics[width=0.35\textwidth]{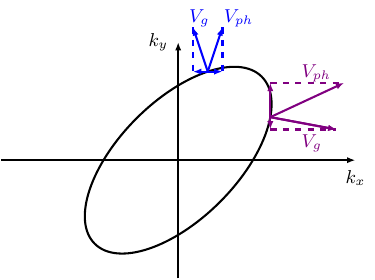}
	\caption{An example slowness curve for an anisotropic wave equation. With blue we mark the vectors of the group and phase velocity that have projections of different signs on the direction $\mathbf{e}_x$, and with violet those that have projections of different signs on the direction $\mathbf{e}_y$. 
Remark that the phase velocity is directed in the radial direction, and the group velocity is directed along the exterior normal to the slowness curve. Clearly, the projection of the group velocity on the phase velocity is always positive. }
	\label{fig:asf}
\end{figure}

A convenient way to study the existence of backward propagating waves in a specified direction is to consider the slowness curves  $S:=\{\mathbf{k}\in\setR^2\colon F(\omega,\mathbf{k})=\omega F(1,\mathbf{k}/\omega)=0\}$. Figure \ref{fig:asf} depicts such a curve for the equation \eqref{eq:main_problem}, where we show the directions of the phase and the group velocities. We see that we can expect the Cartesian PMLs to be unstable both in the directions $\mathbf{e}_x$ and $\mathbf{e}_y$.
However, the notion of the backward wave in the direction $\mathbf{e}$ is something very different from the physical definition of a backward wave: In this sense a wave is called backward if the angle between its phase and group velocity is obtuse, i.e.,  $\mathbf{v}_g(\mathbf{k})\cdot\mathbf{v}_p(\mathbf{k})<0$.
For the configurations which we want to approach in this article, i.e.,\ (classical) non-dispersive wave equations in homogeneous exterior domains, there exist no backward waves in the sense of the physical definition, because their slowness curves are always star-shaped. Therefore, this observation may give us a hope to construct stable PMLs for general anisotropic media. 

\subsection{Radial perfectly matched layers: arguments for stability}\label{sec:radial_pml}

For the radial PML system the frozen coefficient analysis is much more subtle, and, for the moment, we do not know whether it is possible to do the stability analysis based on plane-wave arguments. However, one could naively apply the same ideas about the relation of the instability of the PMLs to the presence of the backward propagating waves in the direction of the PML. 
For radial complex scalings there exists no distinguished direction $\mathbf{e}$ and it is more natural to look in the radial direction, i.e.,\ $\mathbf{e}=\mathbf{k}$, for which, as discussed in the previous section, no backward propagating wave is present, see also Figure \ref{fig:asf} for illustration. 

These ideas can be made more rigorous by computing the fundamental solution associated to the PML system and by showing its decay.
The study of radial PMLs for time-harmonic anisotropic wave equations was initiated in \cite{Halla:22PMLani}, wherein the convergence of a radial PML for anisotropic scalar resonance problems was proven rigorously.
However, the Fredholm and convergence results \cite{Halla:22PMLani} for the time-harmonic equation do not allow any direct deduction of stability for the time-dependent equation.

Another argument which may indicate a potential stability of radial PMLs for anisotropic problems is the following. If we apply the radial PMLs in the free space with $\sigma=\operatorname{const}>0$ (which is the setting in which the Cartesian PMLs are often analyzed), we have that $d=\tilde{d}$, and thus we obtain the following PML problem: 
\begin{align*}
    \partial_t^2 u^{\sigma}+2\sigma \partial_t u^{\sigma}+\sigma^2 u^{\sigma}-\operatorname{div}\bA\nabla u^{\sigma}=0, \quad (t, \bx)\in\setR^+\times \setR^2,  
\end{align*}
and we easily see that the associated energy is non-increasing: 
\begin{align*}
\frac{d}{dt}E=-2\|\sigma^{1/2}\partial_t u^{\sigma}\|^2_{L^2(\setR^2)}, \quad E:=\frac{1}{2}\left(\|\partial_t u^{\sigma}\|^2_{L^2(\setR^2)}+\|\sigma u^{\sigma}\|^2_{L^2(\setR^2)}+\|\bA^{1/2}\nabla u^{\sigma}\|^2_{L^2(\setR^2)}\right). 
\end{align*}

\FloatBarrier 
\subsection{Numerical experiments and their interpretation}
\label{sec:first_experiments}
To motivate our following analysis we present the results of some preliminary numerical experiments which already showcase the quite surprising behavior of radial PMLs for anisotropic materials.
\subsubsection{Two different types of behaviour}
\label{sec:two_different_types_behaviour}
We implement a discretization of the radial PML system \eqref{eq:rad_time} using fourth order finite elements (for details see Section \ref{sec:numerics_ie}, discretization of the interior domain). 
 We choose $\domint=B_1:=\{\bx\in\setR^2:\|\bx\|<1\}$, $\dompml= B_2\setminus \overline {B_1}$ and an anisotropy with $\bA=\diag(1,9)$. We use a piecewise constant damping with $\sigma_c=20$ (cf. Assumption \ref{assump:piecewise_const}) and a time-harmonic source 
$f(t,\bx)=1200\sin(10t)\exp\left(-200\|\bx\|^2\right)$.
The resulting problem is discretized in time with the help of the implicit Crank-Nicholson time-stepping scheme (with the time step $\tau=0.02$) to rule out possible instabilities caused by a CFL condition.
To reduce computational time we exploit the symmetry of the problem and simulate merely one quarter of the domain.

\begin{figure}
  \centering
  \begin{subfigure}[t]{0.45\textwidth}
  \includegraphics[width=\textwidth]{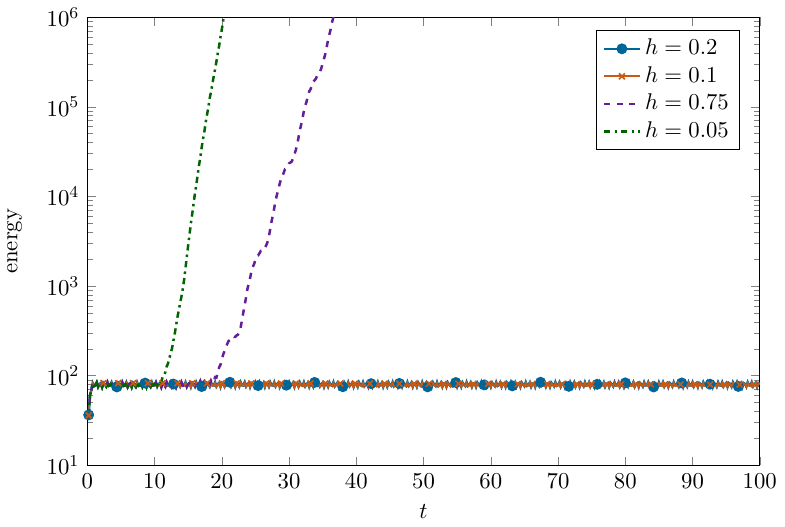}
  \caption{Long-time (in-)stability of PML discretizations for different mesh sizes. The energy curve corresponding to the coarsest mesh size $h=0.2$ is hardly visible, because it is overlapped by the energy curve for $h=0.1$.}
  \label{fig:coarse_stability}
\end{subfigure}\hfill
  \begin{subfigure}[t]{0.45\textwidth}
  \centering
  \includegraphics[width=\textwidth]{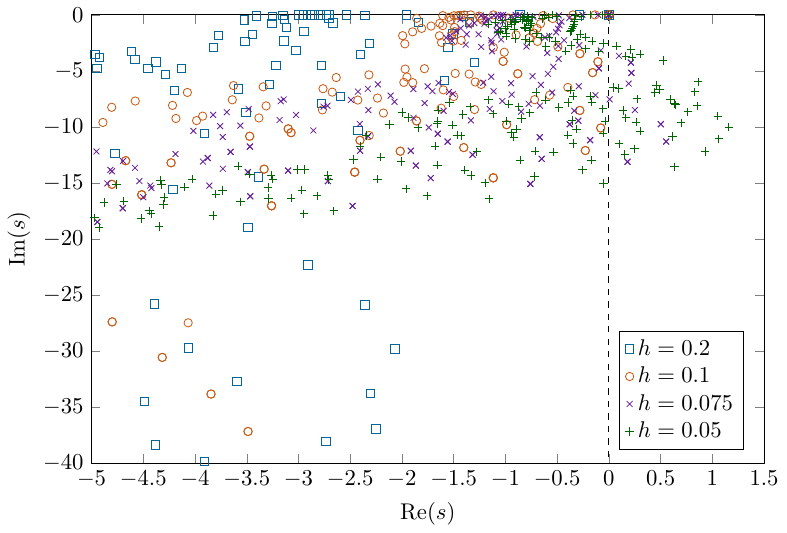}
  \caption{Spectra of the time harmonic problems corresponding to Figure \ref{fig:coarse_stability}.}
  \label{fig:coarse_stability_res}
\end{subfigure}
  \caption{(In-)stability of standard PMLs.}
  \label{fig:coarse_stability_all}
\end{figure}

\begin{figure}
  \centering
  \begin{subfigure}{0.3\textwidth}
    \includegraphics[width=\textwidth,trim={0 150 0 150 },clip]{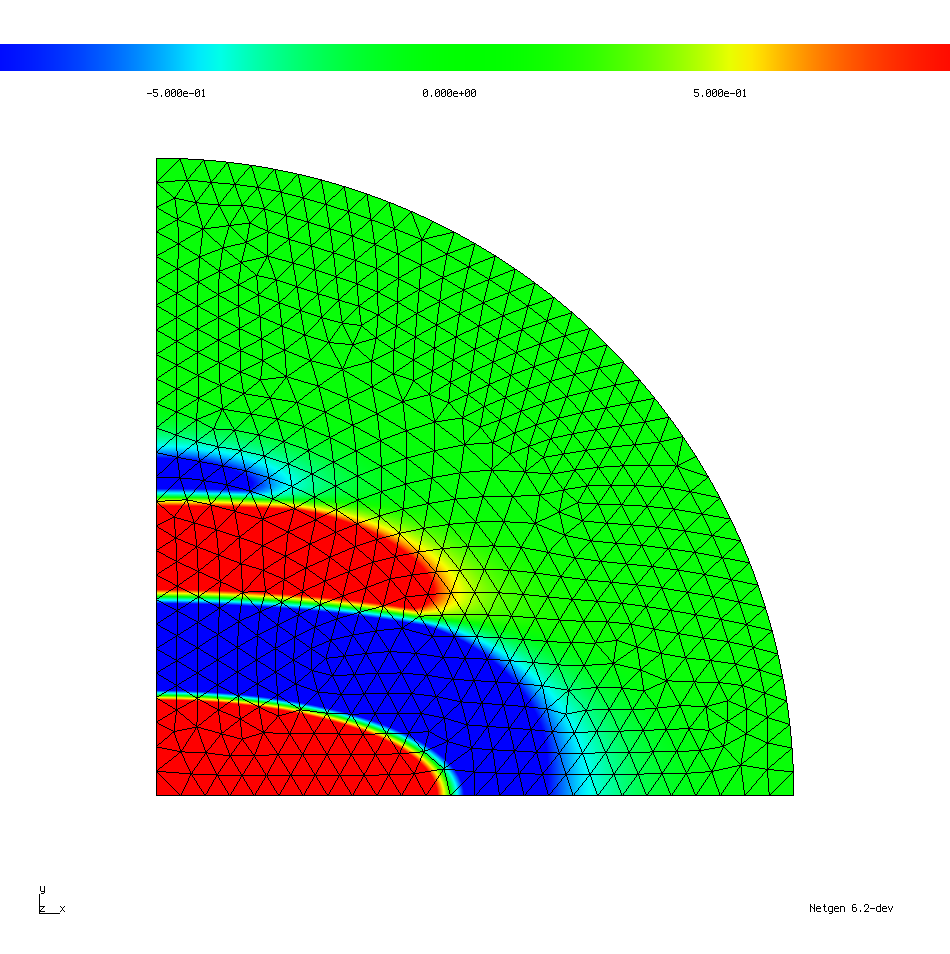}
    \label{fig:pml_unstable_4}
    \caption{$t=4$}
  \end{subfigure}
  \begin{subfigure}{0.3\textwidth}
  \includegraphics[width=\textwidth,trim={0 150 0 150 },clip]{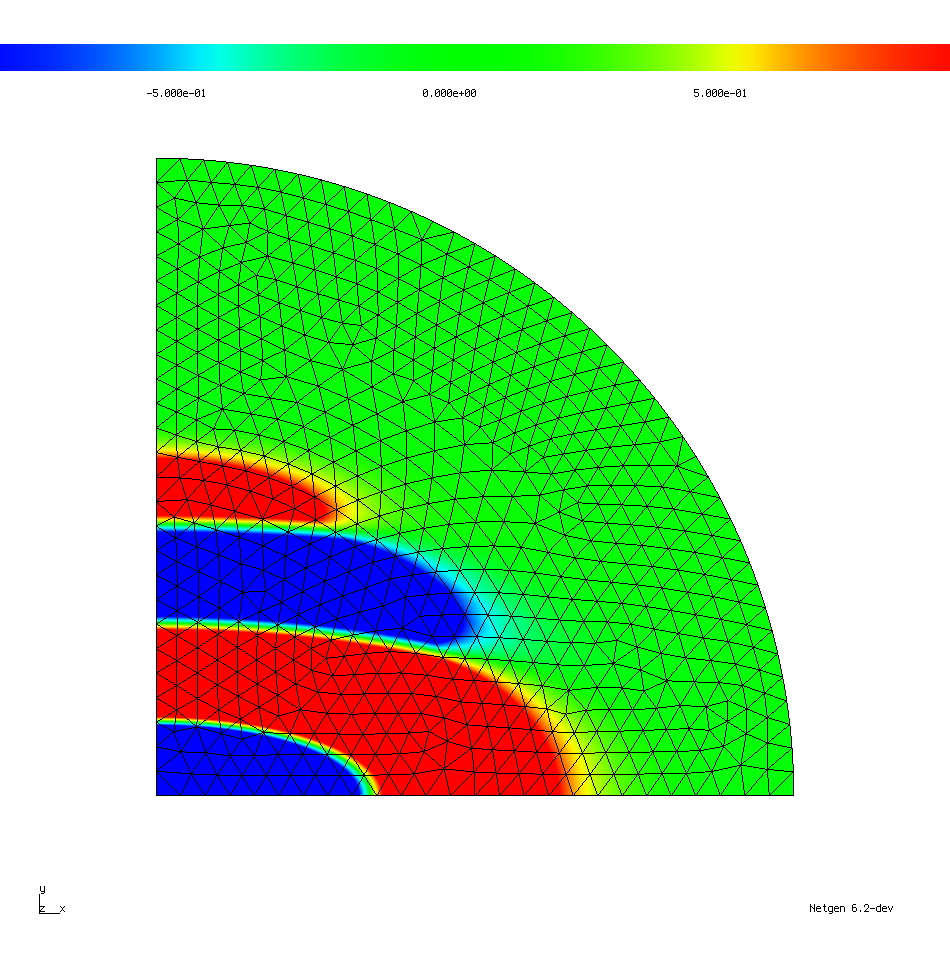}
    \label{fig:pml_unstable_8}
    \caption{$t=8$}
  \end{subfigure}
  \begin{subfigure}{0.3\textwidth}
  \includegraphics[width=\textwidth,trim={0 150 0 150 },clip]{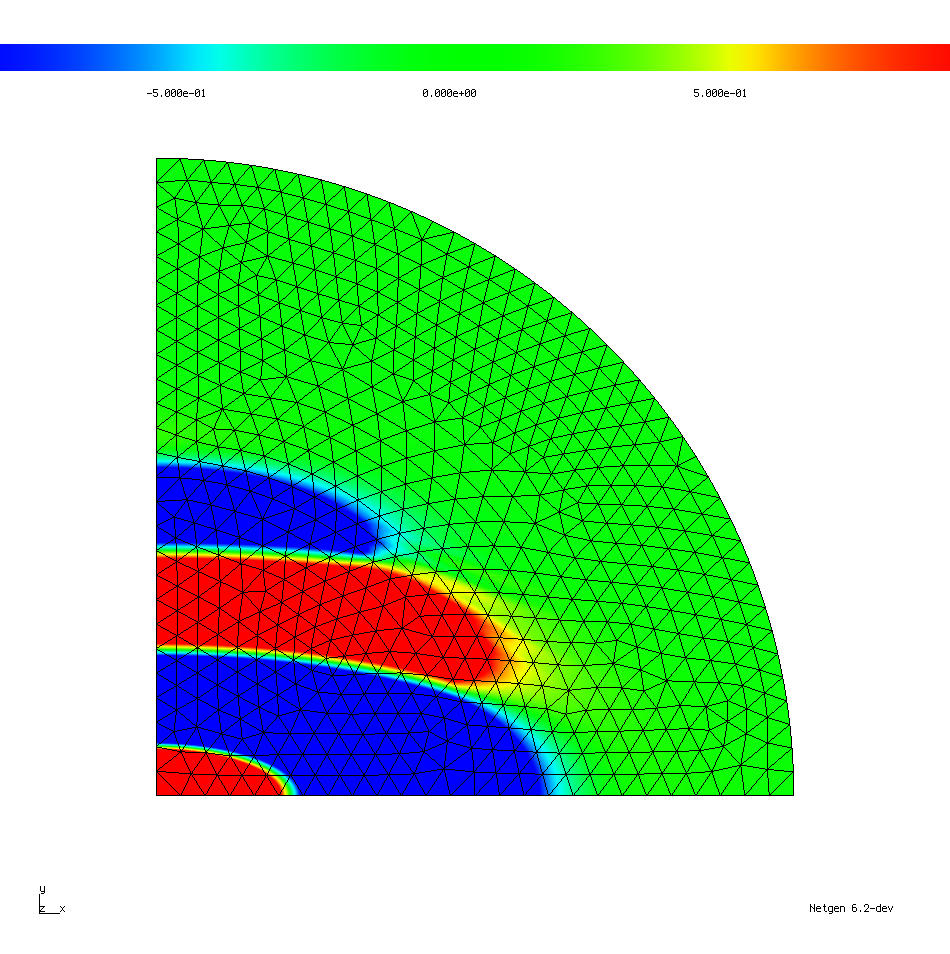}
    \label{fig:pml_unstable_12}
    \caption{$t=12$}
  \end{subfigure}\\
  \begin{subfigure}{0.3\textwidth}
  \includegraphics[width=\textwidth,trim={0 150 0 150 },clip]{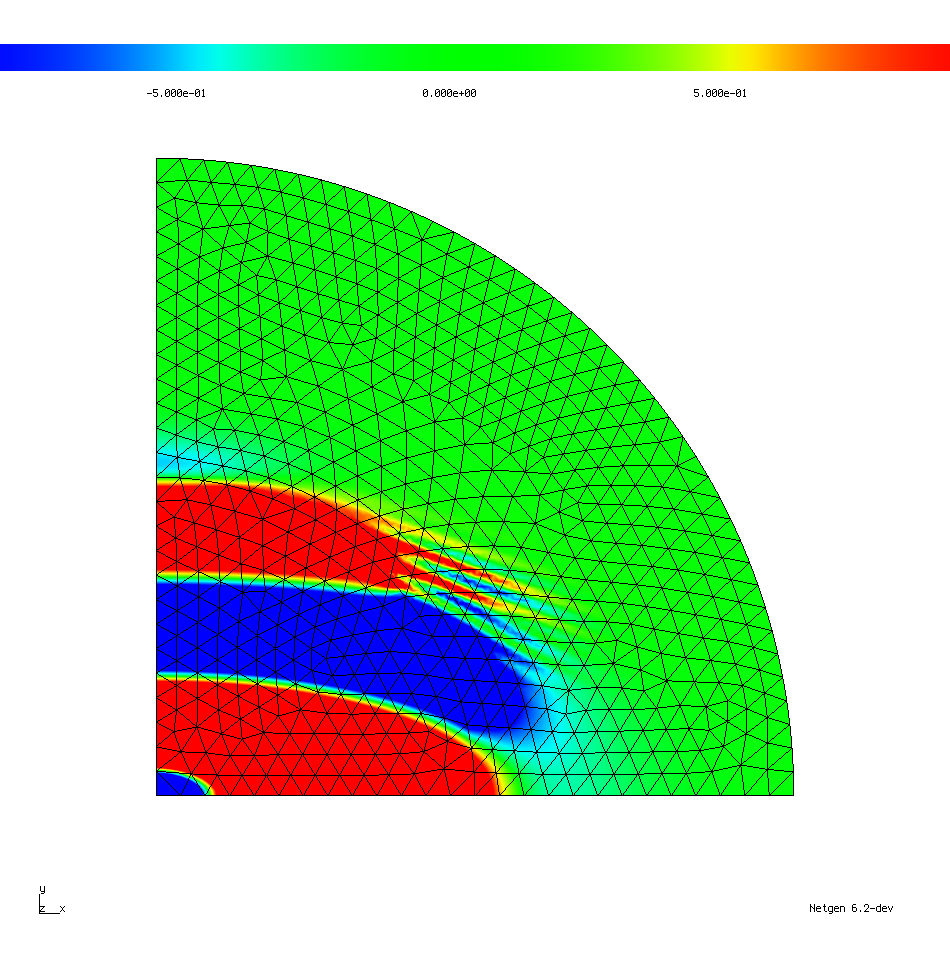}
    \label{fig:pml_unstable_16}
    \caption{$t=16$}
  \end{subfigure}
  \begin{subfigure}{0.3\textwidth}
  \includegraphics[width=\textwidth,trim={0 150 0 150 },clip]{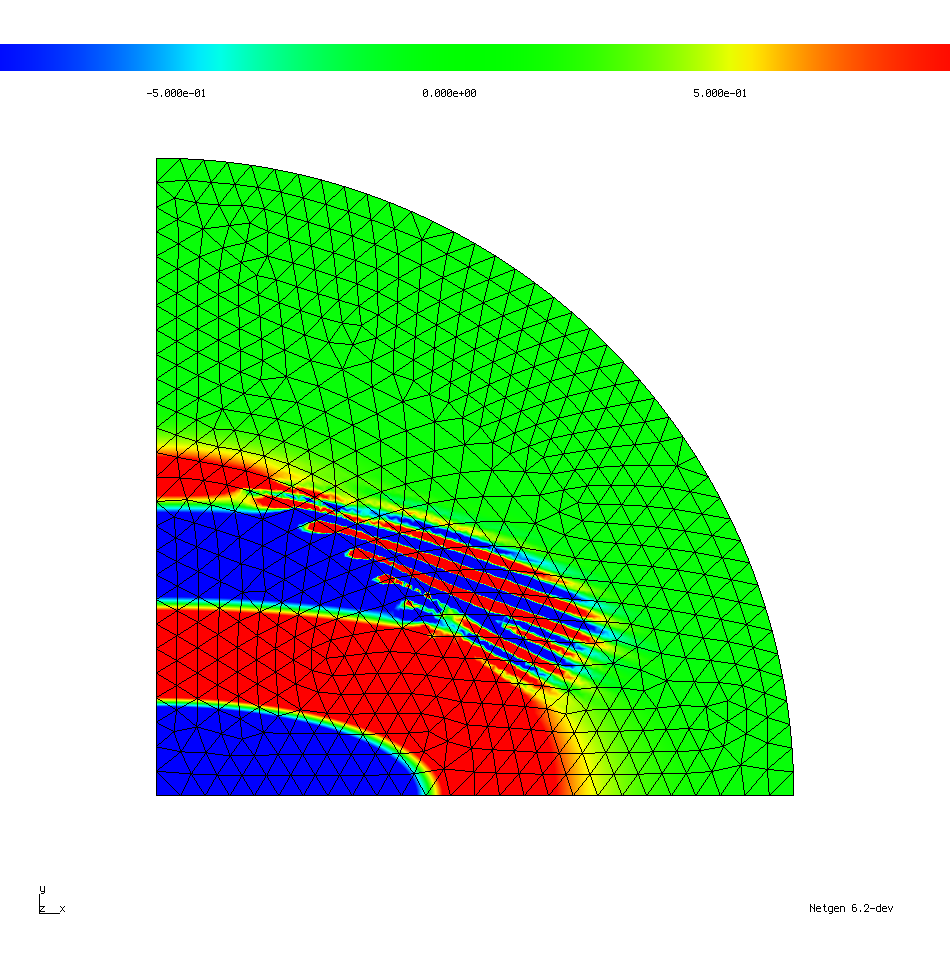}
    \label{fig:pml_unstable_20}
    \caption{$t=20$}
  \end{subfigure}
  \begin{subfigure}{0.3\textwidth}
  \includegraphics[width=\textwidth,trim={0 150 0 150 },clip]{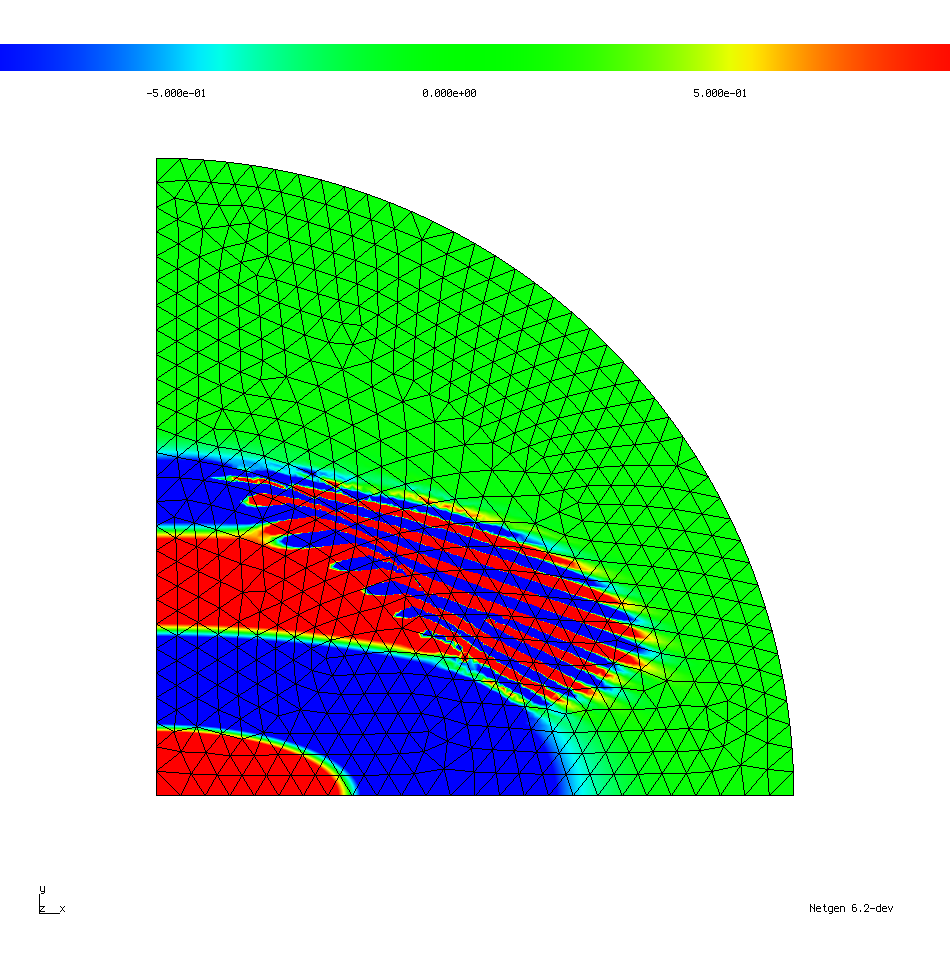}
    \label{fig:pml_unstable_24}
    \caption{$t=24$}
  \end{subfigure}
  \caption{Unstable PML solution with $h=0.075$ at different times $t$ (cf. also Figure \ref{fig:coarse_stability}).}
  \label{fig:pml_unstable}
\end{figure}

Figure \ref{fig:coarse_stability} shows the energy $E$ of the time-domain solutions in $\domint$ for different mesh sizes $h$. We observe that for the coarser meshes ($h=0.2,0.1$) the solution appears to be stable even for large computation times. For finer meshes ($h=0.075,0.05$) we observe an exponential growth of the energy (i.e., unstable solutions, see also Figure \ref{fig:pml_unstable} for snapshots of the unstable time domain solution for $h=0.075$). However, these instabilities occur at different times and with different exponential rates.  
These observations are confirmed by studying the eigenvalues of the eigenvalue problem corresponding to the Laplace transform of the system \eqref{eq:rad_time} where the Laplace parameter $s$ is the unknown spectral parameter.
Figure \ref{fig:coarse_stability_res} shows parts of the spectra of the eigenvalue problems corresponding to the time domain problems from Figure \ref{fig:coarse_stability}. We observe that for the coarser meshes ($h=0.2,0.1$) all of the  approximated eigenvalues have negative real part, which indicates stable time-domain solutions.
The spectra computed using finer meshes contain eigenvalues with positive real parts which correspond to exponentially growing time-domain solutions. Note that the finer the discretization, the larger the maximal real part of the eigenvalues. This explains the different exponential growth rates in the time-domain experiments.
\subsubsection{Discussion of the results}
The numerical results above clearly show that PMLs are not unconditionally stable for anisotropic media. 
Ad hoc there seem to be three possible causes for the instability:
\begin{description}
  \item[a)]{{the untruncated, continuous system \eqref{eq:rad_time} is already unstable,}}
  \item[b)]{{the instabilities are caused by the truncation,}}
  \item[c)] {{the instabilities are caused by the discretization (in space and time).}}
\end{description}
Due to the fact that the instability worsens for finer discretizations there seems to be none to little hope for stability for sufficiently fine discretizations. In the following we will argue that there {exist some configurations where the solution to the untruncated, continuous system is stable, and some where instabilities already occur on the continuous level, i.e., in general the system \eqref{eq:rad_time} is unstable. }Even in the configurations with stable continuous solutions, we will show the existence of an essential spectrum
 in $\setCp$ that causes discrete instabilities for stable continuous solutions.

A similar behavior was already observed in \cite{KreissDuru:13} for Cartesian PMLs and anisotropies which are not aligned with the coordinate axes. Therein the authors exploit that for coarse meshes the discrete solutions are stable. Moreover they provide stability conditions on the mesh size. In the article at hand our goal is to construct a stable numerical method for arbitrarily fine discretizations.

\FloatBarrier 
\section{Well-posedness and stability analysis of radial PMLs}
\label{sec:well_posedness_and_stability}
The instabilities observed in the numerical results of the previous section can be explained either by the ill-posedness or instability of the original problem, or by a wrong choice of the numerical method. In this section we argue that these instabilities are not entirely numerical (the meaning of 'entirely' will become clear later). 
\begin{itemize}
	\item  well-posed, if there exist $m_t, \, m_s, C, \, a, \, p\geq 0$, s.t.\ for each $\boldsymbol{f}\in H^{m_t}(\mathbb{R}^+; (H^{m_s}(\setR^2))^6)$, with $\boldsymbol{f}(0)=\ldots= \partial_t^{m_t-1}\boldsymbol{f}(0)=0$, for each $T>0$,  there exists a unique solution $u^{\sigma}\in L^{\infty}(0, T; L^2(\setR^2))$, s.t.\ 
	\begin{align}
		\label{eq:upml_bound_well_posedness}
		\|u^{\sigma}\|_{L^{\infty}(0, T; L^2(\setR^2))}\leq C\mathrm{e}^{aT}(1+T)^p\|\boldsymbol{f}\|_{H^{m_t}(0, T; H^{m_s}(\setR^2))}.
	\end{align}
	\item stable, if the above holds with $a=0$. 
\end{itemize}
We formulate the well-posedness and stability definitions only for $u^{\sigma}$, since the remaining unknowns can be expressed via $u^{\sigma}$ in a unique way. 
\begin{rmk}
	The initial conditions for the data $\boldsymbol{f}$ ensure that $\boldsymbol{f}$ can be continued to a causal function from the space $H^{m_t}(\mathbb{R};(H^{m_s}(\setR^2))^6)$. 
\end{rmk}
\begin{rmk}
In the above definition, we allow the support of the sources to have a non-empty intersection with the PML, even though in practical applications the starting radius of the PML can always be chosen to avoid this setting.
\end{rmk}
Since we are in the setting of radial PMLs, we can assume that the matrix $\bB=\bA^{-1}$ is diagonal, i.e., $\bB=\diag(\lambda_1,\lambda_2)$, with $\lambda_1, \lambda_2>0$.  Indeed, an arbitrary anisotropy $\bB=\mathbf{U}^{\top}\diag(\lambda_1,\lambda_2)\mathbf{U}$, with $\mathbf{U}$ being an orthogonal transformation, can be handled as follows. Recall that $\mathbf{U}$ is necessarily a product of a rotation $\bR_{\phi_0}$, for some $\phi_0\in [0, 2\pi)$, and a (possible) {coordinate permutation}  $\left(\begin{matrix}
	0 & 1\\
	1 & 0	
\end{matrix}\right)$. Then all the arguments below can be extended by performing a {permutation} of coordinates and a shift $\phi_0$ in the angular variable $\phi$ of $\bA^\phi_\sigma$, for which we have to deal anyway with all possible $\phi\in[0,2\pi)$.  In view of the above discussion, we assume from now on without loss of generality that 
\begin{align}
  \bB:=\bA^{-1}=\diag(\evmax,\evmin),
  \label{eq:a_diag}
\end{align}
with $\evmax\geq\evmin>0$.
For simplicity, we will work with a piecewise constant absorption parameter $\sigma$, i.e., satisfying the following assumption. 
\begin{assump}
	\label{assump:piecewise_const}
	The absorption parameter $\sigma$ satisfies $\sigma(r)=0$ for $r<\radpml$ and $\sigma(r)=\sigma_c>0$ otherwise.
\end{assump}
Some of the results will be valid for more general classes of $\sigma$, and we will make this precise in the statements of those results. 
\begin{thm}
\label{theorem:summary}
Let $\sigma$ satisfy Assumption \ref{assump:piecewise_const}.
Then:
\begin{enumerate}
 \item The problem \eqref{eq:rad_time} is well-posed.
 \item {If for all $t>0$, $\operatorname{supp}\boldsymbol{f}(t, .)\subset B_{r_*}$, where $r_*<\radpml/\mu_*$, with 
\begin{align}
	\label{eq:cstar_def}
\mu_*:=\frac{\evmax+\evmin}{2\sqrt{\evmax\evmin}}\geq 1,
\end{align} 
    i.e., the source is supported sufficiently far away from the absorbing layer, then the bound \eqref{eq:upml_bound_well_posedness} holds with $a=0$ and some $m_t, m_s, C, p\geq 0$.}
\label{enum:thm_summary}
\end{enumerate}
\end{thm}
\begin{proof}
Statement 1 is proven in Section \ref{sec:rad_pml_well_posed} (see Corollary \ref{cor:wp}). Statement 2 is proven in Section \ref{sec:rad_pml_stable_source_proof}, see Corollary \ref{cor:main_stability}.
\end{proof}
The above result shows that the observed numerical instability of radial PMLs is definitely not due to a lack of well-posedness. On the other hand, at the continuous level, the solution to the PML system \eqref{eq:rad_time} is stable for a particular class of sources which are located sufficiently far away from the PML. Since such sources are very special (while, perhaps, being the only ``interesting'' sources for the application of the PMLs), this does not imply that the corresponding problem is stable. Indeed, instabilities manifest themselves at the discrete level, independently of the support of the source term, as illustrated in Section \ref{sec:two_different_types_behaviour} (it is easily verified that in the setting of Section \ref{sec:two_different_types_behaviour} $\mu_*=5/3$ while $|f|$ is far smaller than machine precision outside of $B_{3/5}$, and thus the conditions of Theorem \ref{theorem:summary} are met). 

A possible origin of time-domain instabilities is the presence of singularities in $\setCp$ of the Laplace transform $s\mapsto \hat u^\sigma(s)\in L^2(\setR^2)$ of the time-domain solution $u^\sigma$. These singularities are related to the points $s$ where the operator corresponding to \eqref{eq:bilinear_form} is not invertible. 
The existence of such a spectrum is stated in the following theorem (with the notation $\langle u, u^\dagger\rangle=\int_{\setR^2}u(\bx)\overline{u^\dagger(\bx)}d\bx$).
\begin{thm}
\label{theorem:fredhomlness}
Let $\sigma$ be non-decreasing and continuous for $r\geq \radpml$ (not necessarily piecewise-constant).
  For $s\in \mathbb{C}$ let the sesquilinear form $a_{\sigma}^s(\cdot,\cdot)$ be defined by
\begin{align}
\label{eq:bilinear_form}
a^s_{\sigma}(u, u^\dagger):=\langle \bA_{\sigma}\nabla u, \nabla u^\dagger\rangle
+ s^2 \langle \dtpml \dpml u, u^\dagger\rangle, \quad u, u^{\dagger}\in H^1(\setR^2). 
\end{align}
  Then there exists $s_0\in\setCp$, s.t.\ the Riesz representation $\oppml(s_0)\in \mathcal{L}(H^1(\setR^2))$ of $\sfpml^{s_0}(\cdot,\cdot)$ is not Fredholm. As a corollary, the essential spectrum of $\oppml(\cdot)$ in $\setCp$ (i.e., the set of $s\in\setCp$ where $\oppml(s)$ is not Fredholm) is non-empty.  		
\end{thm}
This theorem is proven in Section \ref{subsubsec:ess-spec}. 
The above result implies that, in general, for suitably chosen sources, $s\mapsto \hat{u}^{\sigma}(s)$ is not $H^1(\setR^2)$-holomorphic in $\setCp$ (see \cite[p.365]{kato}). We believe that a stronger result holds: it is not $L^2(\setR^2)$-holomorphic.
A precise theoretical justification to this is out of scope of the present article.
This lack of analyticity then indicates that we cannot bound the time-domain solution $t\mapsto {u}^{\sigma}(t)\in L^2(\setR^2)$ by a polynomial in time bound, uniformly for any admissible source $\boldsymbol{f}$. 



\subsection{A remark on the analysis}
The sections that follow are dedicated to the proofs of Theorems \ref{theorem:summary} (Sections \ref{sec:rad_pml_well_posed} and \ref{sec:rad_pml_stable_source_proof}) and \ref{theorem:fredhomlness} (Section \ref{subsubsec:ess-spec}). A direct time-domain analysis of well-posedness and stability is quite complicated, and thus we will work using Laplace domain arguments.  For this we consider the system \eqref{eq:rad_time}, with non-vanishing sources $(f, f_v, \boldsymbol{f}_p,\, \boldsymbol{f}_q)$ in the right-hand side. We further rewrite it in the Laplace domain w.r.t. the unknown $u^{\sigma}$, more precisely, we replace \eqref{eq:rad_lpl_v0} by:
\begin{align}
	\label{eq:lpl_u}
	&s^2\dtpml \dpml \hat{u}^{\sigma}-\operatorname{div}_{r, \phi}(\bA_{\sigma}^{\phi}\nabla_{r, \phi}\hat{u}^{\sigma})=\hat{F},
\end{align}
where $\hat{F}=D(s, \dpml(s), \dtpml(s),\partial_x,\partial_y)(\hat{f}, \hat{f}_v, \hat{\boldsymbol{f}}_p,\, \hat{\boldsymbol{f}}_q)^\top$, with the operator $D$ being defined as 
$D(a_1,a_2,a_3,a_4,a_5)=\sum\limits_{j=1}^5 \boldsymbol{b}_j^\top a_j$, for some  
$\boldsymbol{b}_j\in \setR^{6\times 1}$.
Assuming that $\hat{F}\in L^2(\setR^2)$, we will look for a solution of \eqref{eq:lpl_u} belonging to $H^1(\setR^2)$. In terms of the Laplace domain analysis, one can establish the following \textit{sufficient} conditions: 
\begin{itemize}
	\item If there exists $\alpha>0$, s.t. the solution $s\mapsto \hat{u}^{\sigma}(s)$ is an $L^2(\setR^2)$-analytic function in $\setCp_{\alpha}$, and satisfies the following bound
	\begin{align*}
		\|\hat{u}^{\sigma}(s)\|_{L^2(\setR^2)}\leq C\max(1, (\Re s)^{-m})(1+|s|)^p\|\hat{F}\|_{L^2(\setR^2)}, \quad C>0, \quad m,p\geq 0, 
	\end{align*}
then the corresponding time-domain problem is well-posed. See, e.g., \cite{BecacheKachanovska:17} for more details. 
\item If the above holds true with $\alpha=0$, the respective time-domain problem is also stable. 
\end{itemize}
More details on this type of analysis can be found in the monograph by F.-J. Sayas \cite{tdbie_bib}; it was used in the PML context in e.g. \cite{chen},\cite{halpern_rauch}, \cite{BecacheKachanovska:17},  \cite{duru_gabriel_kreiss}, \cite{bkwave}, \cite{bkw}, and also in the time-domain BIE community, cf., e.g., \cite{banjai_lubich_sayas}. 
\subsection{Radial PMLs are well-posed}
\label{sec:rad_pml_well_posed}
Let us consider the problem \eqref{eq:lpl_u} in the variational form and establish the following bound:
\begin{align}
\label{eq:usigmabound}
\|\hat{u}^{\sigma}(s)\|_{H^1(\setR^2)}\leq C (1+|s|)^{p}\|\hat{F}\|_{H^{-1}(\setR^2)}, \quad s\in \setCp_{\alpha}.  
\end{align}
This is done in the lemma below.
\begin{lem}
\label{lem:coercivity}
  Under Assumption \ref{assum:sigma}, there exists $\alpha>0$ depending on $\|\sigma\|_{L^\infty(\setR^+)}, \|\shpml\|_{L^\infty(\setR^+)}, \bA$,  such that $\Re \sfpml^s(u, su)\geq \|su\|^2_{L^2(\setR^2)}+\|\nabla u\|^2_{L^2(\setR^2)}$ for all $s\in\setCp_{\alpha}$.
\end{lem}
\begin{proof}
The main idea of the proof is to use the fact that, as $|s|\rightarrow +\infty$, the sesquilinear form $a_{\sigma}(\cdot,\cdot)$ becomes close to the sesquilinear form without the PML (i.e., $a_0(\cdot,\cdot)$).
First of all, remark that
\begin{align}
	\label{eq:zero_order_term}
	\begin{split}
\Re (s^2\langle \dtpml \dpml u, su\rangle)&=|s|^2\Re s\|u\|^2_{L^2(\setR^2)}+|s|^2\|(\shpml+\sigma)^{1/2}u\|^2_{L^2(\setR^2)}+\Re s\|(\shpml\sigma)^{1/2}u\|^2_{L^2(\setR^2)}.
\end{split}
\end{align}
For the remaining term, we use the definition of the matrix $\bA_{\sigma}=\Js^{-1}\bA\Js^{-\top}\operatorname{det}\Js$ as per \eqref{eq:js}, \eqref{eq:basigma_main}. In particular, for all $s\in\setCp_1$,
with a constant $C$ depending only on $\|\sigma\|_{\infty}, \|\tilde{\sigma}\|_{\infty}, \, \bA$, we have that 
\begin{align}
	\label{eq:upp_bound}
	\left|\langle (\bA_{\sigma}-\bA)\nabla u, s\nabla u\rangle\right|\leq C\|\nabla u\|^2_{L^2(\setR^2)},
\end{align}
because $\|\bA_{\sigma}-\bA\|_2$ is bounded for $s$ with $\Re s>1$ and in particular behaves as $O(|s|^{-1})$ as $|s|\to+\infty$ (since $\Js=\operatorname{Id}+O(|s|^{-1})$, for $|s|\rightarrow +\infty$). 
  Therefore,
\begin{align*}
\Re \langle \bA_{\sigma}\nabla u, s\nabla  u\rangle
&=\Re \langle (\bA_{\sigma}-\mathbf{A})\nabla u, s\nabla u\rangle+\Re \langle \mathbf{A}\nabla u, s\nabla u\rangle\\
&
  \geq \Re s\evmax^{-1}\|\nabla u\|^2_{L^2(\setR^2)} -C\|\nabla u\|^2_{L^2(\setR^2)},
\end{align*}
  where in the last inequality we used \eqref{eq:upp_bound}.  Finally, taking $\alpha=\max(1,\evmax(C+1))$, and combining the above lower bound with the identity \eqref{eq:zero_order_term}, we arrive at the conclusion in the statement of the lemma. 
\end{proof}
With the Lax-Milgram lemma, we conclude that the bound \eqref{eq:usigmabound} holds true. This, together with the arguments of the extended version  \cite[Proposition 3.3]{bkreport} of \cite{BecacheKachanovska:17}, implies the following result.
\begin{cor}
	\label{cor:wp}
	The problem \eqref{eq:rad_time} is well-posed. 
\end{cor}
%
%
%
%
%
%
%
%
%
%

\subsection{Radial PMLs are stable if the source is located far away from the absorbing layer}
\label{sec:rad_pml_stable_source}
Let us remind that all over this section we consider $\sigma$ that satisfies Assumption \ref{assump:piecewise_const}. This assumption is not necessary, but simplifies some of the computations. Moreover, we will use the principal definition of the square root ($\Re\sqrt{z}>0$, $z\in \mathbb{C}\setminus (-\infty, 0]$).

\subsubsection{The fundamental solution of the PML problem \eqref{eq:rad_time} in the Laplace domain and  its properties}
\label{sec:fund_sol_pml}
Because the PML problem \eqref{eq:rad_time} is well-posed, its solution can be expressed with the help of the associated fundamental solution. Analytic expressions of such fundamental solutions in the time domain are difficult to obtain, cf. \cite{diaz_joly} for the Cartesian PMLs in two dimensions, and therefore we are going to work purely in the Laplace domain. First of all, as expected, the fundamental solution for the PML problem coincides with the fundamental solution of the original wave equation with the PML change of variables applied to it. 
\begin{prop}[Fundamental solution of the PML problem]
\label{prop:fs}
Let $\hat{F}\in L^2(\setR^2)$, and let
  $s\in\setCp_{\alpha}$,
where $\alpha>0$ is sufficiently large. Then the unique solution $\hat{u}^{\sigma}\in H^1(\setR^2)$ to \eqref{eq:lpl_u} is given by
\begin{align}
\label{eq:usigma_g}
\hat{u}^{\sigma}(s,\bx):=\int_{\setR^2}G_{\sigma}(s;\bx,\by)\hat{F}(\by)d\by, \quad \bx\in \setR^2, 
\end{align}
  where $G_{\sigma}$ is a fundamental solution of the PML problem \eqref{eq:lpl_u} defined by
\begin{align}
\label{eq:def_d_sigma}
  G_{\sigma}(s; \bx, \by)&:=\frac{1}{2\pi \sqrt{\operatorname{det}\bA}}K_0\left(s\sqrt{(\bx_{\sigma}(s)-\by_{\sigma}(s))^{\top}\bB(\bx_{\sigma}(s)-\by_{\sigma}(s))}\right),&\bB&=\bA^{-1},
\end{align}
  and $\by_{\sigma}$ is defined similarly to $\bx_{\sigma}$ (see \eqref{eq:pmldefns}). 
  In the above $K_0$ is the McDonald function, defined as in \cite[10.25]{nist}, with the branch cut $\setRmz:=\{t\in\setR:t\leq 0\}$.
\end{prop}
\begin{proof}
See Appendix \ref{appendix:Gsigma}. 
\end{proof}
To show how the properties of the fundamental solution in the Laplace domain translate to time-domain bounds on $u^{\sigma}(t,\bx)$, we will study the behaviour of $s\mapsto G_{\sigma}(s; \bx,\by)$ in the complex plane.
\begin{rmk}
  As we will see further, the use of the fundamental solution enables us to show stronger stability results than the ones provided by the study of the PML sesquilinear form in Section \ref{sec:rad_pml_well_posed}.  
\end{rmk}

	First of all, we are interested in the analyticity of the argument of $K_0$ in Proposition \ref{prop:fs}.
	Clearly, the case $\bx, \, \by\in \domint$ reduces to the analysis of the fundamental solution without the PML. The case $\bx,\by\in \domext$ is less interesting, since in applications the source $\hat{F}$ in \eqref{eq:usigma_g} is supported inside $\domint$. Therefore, we consider the case of $\bx\in {\domext}$ and $\by\in {\domint}$. 
	
	To proceed, let us introduce the following two sets: 
		\begin{align}
			\label{eq:dominst_st}
			\begin{split}
				\Omega_{\operatorname{inst}}&:=\{\by\in\domint:\, \exists\, \bx\in \domext, \text{ s.t. } s\mapsto G_{\sigma}(s; \bx,\by) \text{ is not analytic in }\mathbb{C}^+\},\\
				\Omega_{\operatorname{st}}&:=\domint\setminus {\Omega_{\operatorname{inst}}}.
			\end{split}
		\end{align}
		If the source term $\hat{F}$ is s.t. $\operatorname{supp}\hat{F}(s,.)\cap \Omega_{\operatorname{inst}}\neq \emptyset$ for some 'well-chosen' $s\in \mathbb{C}^+$, then we can expect that the corresponding time-domain solution is unstable. On the other hand, the sources $\hat{F}(s,.)$ whose support is inside $\Omega_{\operatorname{st}}$ for all $s\in\mathbb{C}^+$, do not necessarily yield stable time-domain solutions. Nonetheless, as we will see later, if the support of $\hat{F}(s,.)$ is further confined to a pre-defined subdomain of $\Omega_{\operatorname{st}}$, it is possible to prove the stability of the time-domain solution.
	
	%
	%
	%
	The analyticity of $G_{\sigma}$ is linked to the analyticity of the argument of $K_0$ in Proposition \ref{prop:fs}, namely
	the function 
	\begin{align}
		\label{eq:defh}
		\hpml(s;\bx,\by):=(\bx_{\sigma}(\bx)-\by)^{\top}\bB(\bx_{\sigma}(\bx)-\by).
	\end{align}
	Remark that for $s\in\setCp$ s.t.\ $\hpml(s; \bx,\by)\notin \setRmz$, it  holds that $s\sqrt{\hpml(s;\bx,\by)}\notin (-\infty, 0]$.  Indeed, for this it is sufficient to check that $\arg \left(s\sqrt{\hpml(s; \bx, \by)}\right)\in \left(-\pi, \pi\right)$, for all $s\in \setCp$. This follows from $\arg \sqrt{\hpml(s; \bx, \by)}\in \left(-\pi/2, \pi/2\right)$. Therefore, the fundamental solution is well-defined in these points, and we have the following equivalent definition:
	\begin{align}
		\Omega_{\operatorname{inst}}&=\{\by\in\domint:\, \exists\, \bx\in \domext, \text{ s.t. } h_{\sigma}(s; \bx,\by)\leq 0 \text{ for some }s\in \mathbb{C}^+\}.
	\end{align}
	%
	%
	%
	The first principal result of this section reads. 
	\begin{lem}
		\label{lem:lower_bound}
		Let $\delta>0$, $r_*=\radpml/\mu_*-\delta$ and $\mu_*$ like in \eqref{eq:cstar_def}. Then $B(0, r_*)\subset \Omega_{\operatorname{st}}$. Moreover, 
		%
		%
		there exists $C>0$, which depends only on $\radpml$, the matrix $\mathbf{B}$ and $\delta$, s.t.\ for all $s\in \setCp$, we have that,  
		\begin{align*}
			&\Re (s\sqrt{\hpml(s; \bx, \by)})\geq C\Re s\|\bx\|, \quad\text{ for all } \bx\in \domext, \, \by\in B(0,r_*).
		\end{align*}
	\end{lem}
	The proof of this statement relies on several other auxiliary lemmas, and can be found in the end of this section. 
	%
	This result plays a key part in the proof of the stability of the problem, provided that the support of the source is sufficiently separated from the absorbing layer (i.e., Theorem \ref{theorem:summary}). 
	
	The next result of this section, namely Lemma \ref{lem:reh_negative}, provides the first step that allows to characterize the set $\Omega_{\operatorname{inst}}$. Using this result, we will show that for anisotropic media, necessarily, $\Omega_{\operatorname{inst}}\neq \emptyset$, cf. Remark \ref{rmk:negative}.  In Section \ref{sec:remedy} we will use the geometric conditions of Lemma \ref{lem:reh_negative} to construct a stable PML. 
	
\begin{figure}
	\centering
	\includegraphics[width=0.3\textwidth, trim ={0 15 0 0},clip]{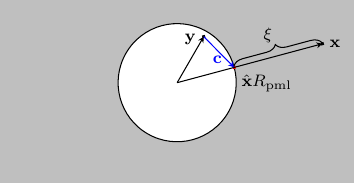}
	\caption{An illustration to the notation \eqref{eq:notation1}. The PML medium is marked in gray and the physical medium in white.}
	\label{fig:notation}
\end{figure}

To proceed with our analysis, let us fix $\bx\in {\domext}$ and $\by\in {\domint}$, and introduce (recall) the following definitions: 
\begin{align}
	\label{eq:notation1}
	&\dpml(s)=1+\frac{\sigma_c}{s}, \quad \bc(\hat{\bx},\by):=\radpml\hat{\bx}-\by, \quad \xi(\|\bx\|):=\|\bx\|-\radpml.
\end{align}
In the above and what follows, we denote by $\hat{\bx}:=\frac{\bx}{\|\bx\|}$. Remark that for $\bx$ inside the PMLs $\dpml(s,\|\bx\|)$ is constant, therefore we use the notation $\dpml(s)$ instead. 
With these definitions, we have in particular that 
\begin{align*}
	\bx_{\sigma}=\left(\|\bx\|+\frac{\sigma_c}{s}(\|\bx\|-\radpml)\right)\hat{\bx}
	=
	\left(1+\frac{\sigma_c}{s}\right)\xi(\|\bx\|)\hat{\bx}+\radpml\hat{\bx}.
\end{align*}
Therefore, the quantity $\hpml(s; \bx, \by)$ that we wish to study can be rewritten as follows:
\begin{align}
	\label{eq:defh2}
	\hpml(s; \bx, \by)&=(\bx_{\sigma}-\by)^{\top}\bB(\bx_{\sigma}-\by)=(\bc+\hat{\bx}\xi\dpml)\bB(\bc+\hat{\bx}\xi\dpml)\\
	\nonumber
	&\phantom{:}=\bc^{\top}\bB\bc+\dpml^2(s)\xi^2\hat{\bx}^{\top}\bB\hat{\bx}+2\dpml(s)\xi\bc^{\top}\bB\hat{\bx}. 
\end{align}
With $\hat{\bc}=\bc\|\bc\|^{-1}$, we define the following quantities:	
\begin{align}
	\label{eq:notation2}
	\gamma_{11}(\hat{\bx},\by):=\hat{\bc}^{\top}\bB\hat{\bc},
	\quad \gamma_{22}(\hat{\bx}):=\hat{\bx}^{\top}\bB\hat{\bx},
	\quad\text{ and }\quad\gamma_{12}(\hat{\bx},\by):=\hat{\bx}^{\top}\bB\hat{\bc}. 
\end{align}
This allows us to rewrite $\hpml(s; \bx,\by)$ as follows:
\begin{align*}
	\hpml=\|\bc\|^2\gamma_{11}+2\dpml\xi\|\bc\|\gamma_{12}+\dpml^2\xi^2\gamma_{22}. 
\end{align*}
An illustration to the notation \eqref{eq:notation1} is given in Figure \ref{fig:notation}.
We are interested in the behaviour of $\hpml$ for 
$s\in\setCp$
and $\xi\geq 0$, in particular in the points where $\hpml(s;\bx,\by)\in \setR_{0}^{-}$, so that the argument of the fundamental solution is not analytic.
We then have the following result.
\begin{lem}
	\label{lem:reh_negative}
	Let $\bx\in \domext$ and $\by\in \domint$ be fixed. Then $\hpml(s; \bx, \by)\leq 0$ if and only if $s\in \setCp$ is s.t.\ the following two conditions hold true:
	\begin{align}
		\label{eq:lem_main}
		(a)\, \Re \dpml(s)=-\frac{\|\radpml\hat{\bx}-\by\|}{\|\bx\|-\radpml}\frac{\gamma_{12}(\hat{\bx},\by)}{\gamma_{22}(\hat{\bx})};
		\quad (b)\, \cos^2(\arg \dpml(s))\leq \frac{\gamma_{12}^2(\hat{\bx},\by)}{\gamma_{11}(\hat{\bx},\by)\gamma_{22}(\hat{\bx})}. 
	\end{align}
\end{lem}
\begin{proof}
	Computing the real and imaginary part of $\hpml$ explicitly yields with
	{$d_{\sigma,r}:=\Re\dpml$, $d_{\sigma,i}:=\Im\dpml$} {and the notation from \eqref{eq:notation2}}
	{%
		\begin{align}
			\label{eq:h_re_im}
			\hpml&=\|\bc\|^2\gamma_{11}+2d_{\sigma,r}\xi\|\bc\|\gamma_{12}+\left(d_{\sigma,r}^2-d_{\sigma,i}^2\right)\xi^2\gamma_{22}+2id_{\sigma,i}\xi\left(\|\bc\|\gamma_{12}+d_{\sigma,r}\xi \gamma_{22}\right).
	\end{align}}
	Thus
	$\Im \hpml(s; \bx, \by)=0$ if and only if $d_{\sigma,i}=0$ or (a) is satisfied, i.e.,
	\begin{align}
		\label{eq:lem_main_a}
		d_{\sigma,r}(s)=-\frac{\|\bc\|}{\xi}\frac{\gamma_{12}}{\gamma_{22}}.
	\end{align}
	The identity $d_{\sigma,i}=0$ is equivalent to $s\in \setRp$ which implies $\hpml(s; \bx, \by)>0$.
	This leaves us with (a). For $\hpml(s; \bx, \by)\leq 0$, it is further necessary that
	\begin{align*}
		\|\bc\|^2\gamma_{11}+2d_{\sigma,r} \xi\|\bc\|\gamma_{12}+(d_{\sigma,r}^2-d_{\sigma,i}^2)\xi^2\gamma_{22}\leq 0.
	\end{align*}
	Re-expressing $\xi$ via $d_{\sigma,r}$ from \eqref{eq:lem_main_a} (or \eqref{eq:lem_main}(a)) yields 
	\begin{align*}
		\gamma_{11}-2\frac{\gamma_{12}^2}{\gamma_{22}}+d_{\sigma,r}^{-2}(d_{\sigma,r}^2-d_{\sigma,i}^2)\frac{\gamma_{12}^2}{\gamma_{22}}\leq 0.
	\end{align*}
	This can be rewritten as 
	\begin{align}
		\label{eq:di_dr}
		\frac{d_{\sigma,i}^2}{d_{\sigma,r}^2}+1\geq \frac{\gamma_{11}\gamma_{22}}{\gamma_{12}^2}.
	\end{align}
	The above is equivalent to the condition \eqref{eq:lem_main}(b).
\end{proof}
To understand the statement of Lemma \ref{lem:reh_negative}, remark that the condition $(b)$ is a condition on a 'source' point $\by$ and on the direction $\hat{\bx}\in \mathbb{S}^2$. The condition (a) is more involved, and involves the distance between the point $\bx$ inside the PML and the interface between the PML and the physical media.

\begin{figure}
	\centering
	\includegraphics[width=0.4\textwidth]{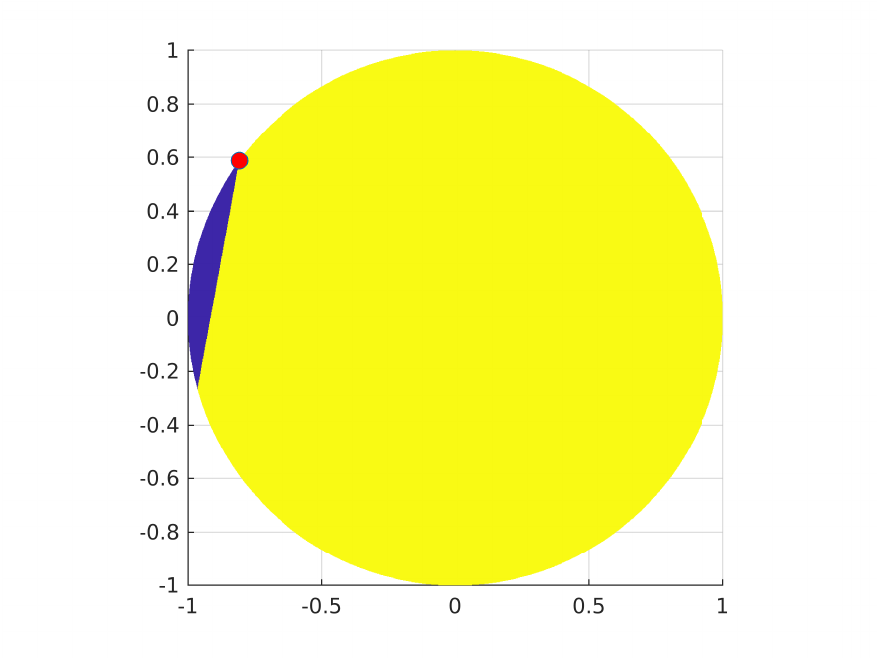} 
	\includegraphics[width=0.4\textwidth]{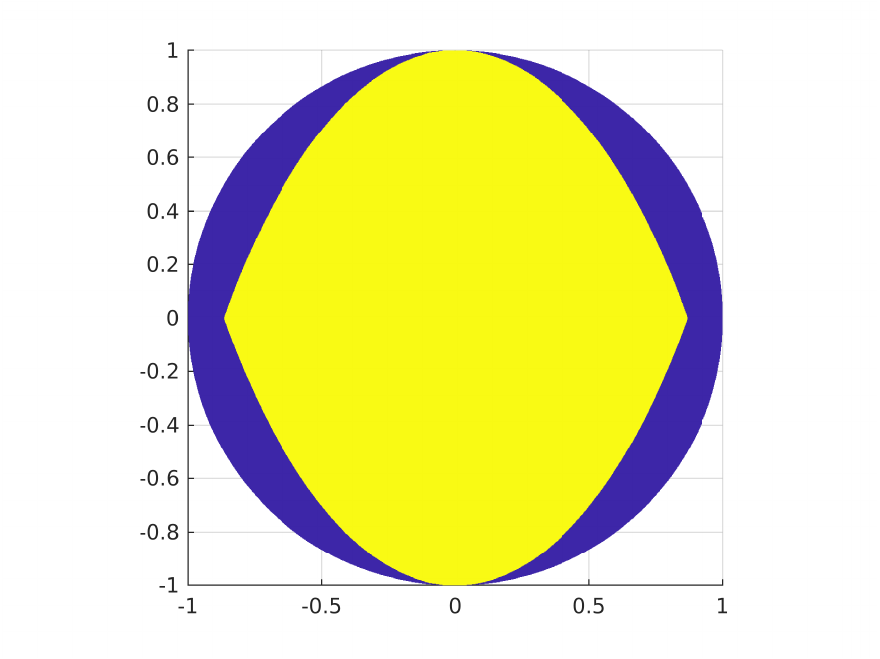} 	
	\caption{Left: the quantity $\operatorname{sign}\gamma_{12}(\hat{\bx},\by)$ as defined in \eqref{eq:notation2} depending on $\by\in B(0,1)$, for a fixed $\bx=(\cos0.8\pi,\sin0.8\pi)$, $\radpml=1$. In dark blue we mark $\by$, s.t.\ $\operatorname{sign}\gamma_{12}(\hat{\bx},\by)=-1$, while in yellow $\by$ with  $\operatorname{sign}\gamma_{12}(\hat{\bx},\by)=1$. The point $\bx$ is marked in red. Right: the sets $\Omega_{\operatorname{st}}$ (in yellow) and $\Omega_{\operatorname{inst}}$ (in blue). 
		%
		In all the experiments, the matrix $\bB=\operatorname{diag}(1,1/4)$. }
	\label{fig:illustr1}
\end{figure}

Let us now study the validity of the conditions \eqref{eq:lem_main} $(a)$, $(b)$ for arbitrary points $\bx\in \domext$ and $\by\in \domint$. It appears that the condition $(b)$ always holds true for some $s\in \setCp$ and almost all directions $\hat{\bx}$.
\begin{lem}
	\label{lem:b_always_true}
	Let $\by\in \domint$ and $\hat{\bx}\in \mathbb{S}^2$ be s.t. $\gamma_{12}(\hat{\bx},\by)\neq 0$. Then there exists $s\in \setCp$, s.t.\ the condition (b) in \eqref{eq:lem_main} holds true. 
\end{lem}
\begin{proof}
	First, the right-hand side of (b) in \eqref{eq:lem_main} is strictly positive.  
	Next, remark that $s\mapsto \dpml(s)=1+\frac{\sigma_c}{s}$ is a bijection from $\setCp$ to $\setCp_1$. Thus, $s\mapsto\arg \dpml(s)$ takes all the values in $\left(-\frac{\pi}{2}, \frac{\pi}{2}\right)$, hence the conclusion. 
\end{proof}
The validity of the condition (a), however, depends on whether the medium is isotropic or anisotropic. In particular, for isotropic media, the condition (a) never holds true, and we have that $\hpml(s;\bx,\by)\notin \setR_{0}^{-}$.  
\begin{lem}
	\label{lem:isotropic_positive}
	For the isotropic medium ($\bB=a\operatorname{Id}$, $a>0$), for all $\hat{\bx}\in \mathbb{S}^2$ and $\by\in\domint$, it holds that $\gamma_{12}=\gamma_{12}(\hat{\bx}, \by)>0$. Therefore, for all $s\in \setCp$, $\hpml(s; \bx,\by)\notin \setR_{0}^{-}$, and thus  $\Omega_{\operatorname{st}}=\Omega_{\operatorname{int}}$.
\end{lem}
\begin{proof}
	Since for isotropic media $\bB=a\operatorname{Id}$, $a>0$, we have that
	\begin{align*}
		\operatorname{sign}\gamma_{12}(\hat{\bx},\by)=\operatorname{sign}\left(\hat{\bx}^{\top}(\radpml\hat{\bx}-\by)\right)=\operatorname{sign}\left(\hat{\bx}^{\top}\bc\right).
	\end{align*}
	The condition $\|\by\|<\radpml$ is equivalent to $\|\radpml\hat{\bx}-\bc\|^2<\radpml^2$, which, after straightforward computations yields $\|\bc\|^2-2\radpml \hat{\bx}^\top \bc<0$. Thus $\hat{\bx}^\top\bc>\|\bc\|^2/(2\radpml)>0$ and hence $\gamma_{12}(\hat{\bx},\by)>0$.
	Since the left hand-side of \eqref{eq:lem_main}(a) is always positive, we conclude with Lemma \ref{lem:reh_negative}.
\end{proof}
For an anisotropic medium, the above result for the sign of  $\gamma_{12}(\hat{\bx},\by)$ does not hold true. Indeed, 
\begin{align*}
	\operatorname{sign}\gamma_{12}(\hat{\bx},\by)=\operatorname{sign}\left(\hat{\bx}^{\top}\bB\bc\right),
\end{align*}
and the positivity of $\hat{\bx}^{\top}\bc$ does not imply that the above quantity is positive. 
This is illustrated in Figure \ref{fig:illustr1}.
\begin{rmk}
	\label{remark:existence_xy}
	It is possible to prove that when eigenvalues of $\bB$ differ, i.e. $\evmin\neq \evmax$, then there always exists $\by\in \domint$ and a direction $\hat{\bx}_0\in \mathbb{S}^2$, s.t. $\gamma_{12}(\hat{\bx}_0,\by)<0$. 
\end{rmk}
The condition on the sign of $\gamma_{12}$ is of utmost importance, and will help us to characterize the set $\Omega_{\operatorname{inst}}$ in Lemma \ref{lem:omegainst}. Its proof relies on the following two observations. 
\begin{lem}
	\label{lem:a_true_when_gamma_negative}
	Assume that there exist $\hat{\bx}_0\in \mathbb{S}^2$ and $\by\in \domint$, s.t.\ $\gamma_{12}(\hat{\bx}_0,\by)<0$. Then for any $s\in \setCp$, the condition (a) of \eqref{eq:lem_main} holds true for $\bx\in \domext$ defined by $\bx=\rho\hat{\bx}_0$ with 
	\begin{align*}
		\rho=\radpml-\|\radpml \hat{\bx}_0-\by\|\frac{\gamma_{12}(\hat{\bx}_0,\by)}{\Re\dpml(s)\gamma_{22}(\hat{\bx}_0)}\radpml>\radpml.
	\end{align*} 
\end{lem}
The proof of this lemma is left to the reader. 
It also implies the following positive result.
\begin{cor}
	\label{cor:main_h}
	Let $\by\in \domint$. If $\gamma_{12}(\hat{\bz},\by)\geq 0$ for all $\hat{\bz}\in \mathbb{S}^2$, then $\hpml(s;\bx,\by)\notin\setR_{0}^{-}$ for all $s\in \setCp$. Therefore, $\by\in \Omega_{\operatorname{st}}$. 
\end{cor}
%
Now we are in the position to prove the following lemma. 
	\begin{lem}
		\label{lem:omegainst}
		The set $\Omega_{\operatorname{inst}}$ admits the following characterization:
		\begin{align}
			\label{eq:omegainst}
			\Omega_{\operatorname{inst}}=\{\by\in \domint: \, \text{ there exists }\hat{\bx}_0\in\mathbb{S}^2, \, \text{ s.t. }\gamma_{12}(\hat{\bx}_0, \by)<0\}.
		\end{align}
	\end{lem}
	\begin{proof}[Proof of Lemma \ref{lem:omegainst}]
		The inclusion $\Omega_{\operatorname{inst}}\subseteq\ldots$ follows from Corollary \ref{cor:main_h}. 
		To prove the opposite inclusion, it suffices to show that as soon as $\gamma_{12}(\hat{\bx}_0, \by)<0$ for some $\hat{\bx}_0\in \mathbb{S}^2$, then $h_{\sigma}(s; \bx,\by)\in \mathbb{R}_0^{-}$ for some $s\in \mathbb{C}^+$, $\bx\in\domext$. 
		
		Fix the point $\by\in \domint $ and the direction $\hat{\bx}_0$, s.t. $\gamma_{12}(\hat{\bx}_0,\by)<0$. Next choose $s\in \mathbb{C}^+$ like in Lemma \ref{lem:b_always_true} and set $\bx=\rho\hat{\bx}_0$, with $\rho$ chosen as in Lemma \ref{lem:a_true_when_gamma_negative}. It remains to apply Lemma \ref{lem:reh_negative} for $h_{\sigma}(s; \bx, \by)$.
	\end{proof}
\begin{rmk}
		\label{rmk:negative} 
		With Remark \ref{remark:existence_xy} Lemma \ref{lem:omegainst} implies that $\Omega_{\operatorname{inst}}\neq \emptyset$ whenever the medium is anisotropic (i.e. when the eigenvalues of $\mathbf{B}$ differ).
	\end{rmk}
In Fig.\ \ref{fig:illustr1}, right, we plot the sets $\Omega_{\operatorname{st}}$ and $\Omega_{\operatorname{inst}}$. It indicates that $\Omega_{\operatorname{st}}$ contains a  ball around the origin.
This is quantified below.
\begin{lem}
	\label{lem:stability_condition1}
	
	Let $\mu_*$ be like in \eqref{eq:cstar_def}. 
	For all $\hat{\bz}\in \mathbb{S}^2$ and $\by\in \domint$ s.t.\
	$\|\by\|<\radpml/\mu_*$, 
	the function $\mathbb{S}^2\ni\hat{\bz}\mapsto \gamma_{12}(\hat{\bz},\by)$ is strictly positive. Therefore, cf.,~Lemma~\ref{lem:omegainst},  $B(0,{\radpml/\mu_*})\subset\Omega_{\operatorname{st}}$. 
	%
	%
\end{lem}
\begin{proof}
	By Corollary \ref{cor:main_h}, to have that $\hpml(s; \bx, \by)\in \mathbb{C}\setminus \setR_{0}^{-}$ it is sufficient to ensure that $\by$ is s.t. $\mathbb{S}^2\ni\hat{\bz}\mapsto \gamma_{12}(\hat{\bz},\by)$ is strictly positive. We have that 
	\begin{align*}
		\operatorname{sign}\gamma_{12}(\hat{\bz},\by)=\operatorname{sign}\left(\radpml \hat{\bz}^{\top}\bB\hat{
			\bz}-\|\by\|\hat{\bz}^{\top} \bB\hat{\by}\right).
	\end{align*}
	By straightforward calculations carried out in Lemma \ref{lem:upperbound_cxx_cxy}, for any $\hat{\bz}_1,\, \hat{\bz}_2\in \mathbb{S}^2$, 
	\begin{align*}
		\left|\frac{\hat{\bz}^{\top}_1 \bB\hat{\bz}_2}{\hat{\bz}^{\top}_1\bB\hat{\bz}_1}\right|\leq \mu_*,
	\end{align*}
	where $\mu_*$ is defined as in the statement of the present lemma. Then it follows immediately that $\mathbb{S}^2\ni\hat{\bz}\mapsto \gamma_{12}(\hat{\bz},\by)>0.$
	%
\end{proof}
{Finally, we have all the necessary ingredients to prove Lemma \ref{lem:lower_bound}, which gives a refined bound on $\hpml(s; \bx, \by)$. }
\begin{proof}[Proof of Lemma \ref{lem:lower_bound}]
	We start with $\Re (s\sqrt{\hpml})=\Re s\Re \sqrt{\hpml}-\Im s\Im \sqrt{\hpml}$.
	Let us first prove that  $\Im\sqrt{\hpml}$ is of the same sign as $-\Im s$. For this let us recall that by \eqref{eq:h_re_im} 
	\begin{align*}
		\Im \hpml(s; \bx, \by)=2\Im \dpml(s)\left(\xi\|\bx-\by\|\gamma_{12}+\Re \dpml(s)\xi^2 \gamma_{22}\right), 
	\end{align*}
	and as $\gamma_{12}>0$ (cf. Lemma \ref{lem:stability_condition1}), $\gamma_{22}>0$ and $\Re \dpml(s)>0$, we obtain that $$\operatorname{sign}\Im \hpml(s; \bx, \by)=\operatorname{sign}\Im \dpml(s)=\sign\Im (\spml/s)=-\sign\Im s \text{ if }\xi\neq 0,$$ and
	$\Im \hpml(s; \bx,\by)=0$ otherwise. 
	Therefore, 
	\begin{align}
		\label{eq:princ_ident}
		\Re(s\sqrt{\hpml(s; \bx, \by)})\geq \Re s\Re \sqrt{\hpml(s; \bx, \by)}. 
	\end{align}
Let us now remark that $\hpml$ depends on $s$ via $\bx_{\sigma}$ only. Let us thus define $H\colon \mathbb{C}^2\times \overline{B(0, r_*)}$ via $H(\mathbf{z}, \mathbf{y}):=(\mathbf{z}-\by)^{\top}\bB(\mathbf{z}-\by)$, so that $H(\bx_{\sigma}, \by)=\hpml(s; \bx, \by)$. 

	For $(\bx_{\sigma}, \by)$ belonging to a compact set $K\times \overline{B(0, r_*)}$, $\Re {\sqrt{H(\bx_{\sigma},\by)}}>c_K>0$, since $H(\bx_{\sigma},\by)=\hpml(s; \bx,\by)$ does not belong to $\mathbb{R}_0^{-}$. 
	For large $\|\bx_{\sigma}\|\rightarrow +\infty$, we have that, uniformly in $\|\by\|\leq \radpml/\mu_*-\delta$,  $$\sqrt{H( \bx_{\sigma},\by)}=\sqrt{\bx_{\sigma}^{\top}\mathbf{B}\bx_{\sigma}}\left(1+o(1)\right)=\left(1+\frac{\sigma_c}{s}\right)\sqrt{\bx^{\top}\mathbf{B}\bx}\left(1+o(1)\right).$$
	We then have 
	\begin{align*}
		\Re \left(1+\frac{\sigma_c}{s}\right)\sqrt{\bx^{\top}\mathbf{B}\bx}\geq\sqrt{\evmin}\|\bx\|, 
	\end{align*}
	which yields the result in the statement of the lemma.
\end{proof}

\subsubsection{Proof of Theorem \ref{theorem:summary}.(\ref{enum:thm_summary})}
\label{sec:rad_pml_stable_source_proof}
This section is dedicated to proving that radial PMLs for anisotropic media are stable, as long as the support of the source is sufficiently separated from the PML.
We start with the following auxiliary result. 
\begin{lem}[The source far from the PML]
\label{lemma:source_far}
Let $r_*:=\radpml/\mu_*-\delta$, where $\mu_*$ is defined in \eqref{eq:cstar_def} and $\delta>0$. We have the following bound for all $s\in \setCp$:
\begin{align*}
&\int_{\mathbb{R}^2}\int_{B(0, r_*)}|G_{\sigma}(s;\bx, \by)|^2 d\by\,d\bx \leq C \max\left(1, \left(\frac{1}{\Re s}\right)^2\right),  
\end{align*}
 where the constant $C$ depends only on $\radpml$, $\delta$,  and the matrix $\mathbf{B}$.
\end{lem}
\begin{proof}
We split the integral 	
	\begin{align*}
		&\int_{\mathbb{R}^2}\int_{B(0, r_*)}|G_{\sigma}(s;\bx, \by)|^2 d\by\,d\bx =	I_0+I_{\sigma}, \\
		& I_0=	\int_{\domint}\int_{B(0, r_*)}|G_{\sigma}(s;\bx, \by)|^2 d\by\,d\bx
		, \quad I_{\sigma}=\int_{\domext}\int_{B(0, r_*)}|G_{\sigma}(s;\bx, \by)|^2 d\by\,d\bx.
	\end{align*}
Next, let us bound the above two integrals. We start with a less classical one, namely $I_{\sigma}$, and next argue that the corresponding estimate can be extended to $I_0$. 

\textit{Step 1. A bound on $I_{\sigma}$. } Recall an explicit expression for $G_{\sigma}$ given in Proposition \ref{prop:fs}. We start by bounding $G_{\sigma}$. 
	First of all, remark that the function $z\mapsto K_0(z)$ is analytic in $\mathbb{C}\setminus \setR^{-}$. 
By \cite[10.32.9]{nist}, it holds that 
	\begin{align*}
		K_0(z)=\int_0^{\infty}\mathrm{e}^{-z\cosh t}dt,\quad z\in \setCp, \text{ therefore }|K_0(z)|\leq |K_0(\Re z)| \text{ in }\setCp.
	\end{align*}
 Moreover, from the above representation it follows that the function $\setR^+\ni\lambda\mapsto K_0(\lambda)$ is decreasing. Therefore, by the above considerations and Lemma \ref{lem:lower_bound}, we have that $$|G_{\sigma}(s; \bx,\by)|\leq \frac{1}{2\pi \sqrt{\operatorname{det}\bA}}|K_0(\Re (s\sqrt{\hpml(s; \bx,\by)}))|\leq \frac{1}{2\pi \sqrt{\operatorname{det}\bA}}|K_0(C\Re s\|\bx\|)|.$$ 
This allows to rewrite
\begin{align*}
	I_{\sigma}&\leq \frac{1}{4\pi^2\operatorname{det}\bA}\int_{\domext}\int_{B(0,r_*)}|K_0(C\Re s\|\bx\|)|^2 d\bx d\by\leq \frac{ r_*^2}{4\pi\operatorname{det}\bA}\int_{\domext}|K_0(C\Re s\|\bx\|)|^2d\bx\\
	&\leq \frac{ r_*^2}{4\pi\operatorname{det}\bA}\int_{\mathbb{R}^2}|K_0(C\Re s\|\bx\|)|^2d\bx=\frac{ r_*^2}{4\pi(\Re s)^2C^2\operatorname{det}\bA}\int_{\mathbb{R}^2}|K_0(\|\bx\|)|^2d\bx.
\end{align*}
By \cite[10.30]{nist}, for any $a>0$ (where the constant $C_a>0$ in the bound depends on $a>0$ as well), it holds that 
	\begin{align}
		\label{eq:K0Z}
		|K_0(z)|\leq C_a \left\{
		\begin{array}{ll}
			\max(1,|\log|z||) & |z|<a, \\
			\mathrm{e}^{-\Re z}, & |z|\geq a.
		\end{array}
		\right.
	\end{align}
The above bound shows that $\bx\mapsto K_0(\|\bx\|)\in L^2(\mathbb{R}^2)$. This yields the validity of the bound in the statement of the lemma for $I_{\sigma}$. 

\textit{Step 2. A bound on $I_0$. }
First of all, we remark that $I_0\leq\|G_{\sigma}(s; ., .)\|^2_{L^2(\domint\times\domint)}$. 
Now recall that for $\bx,\by\in \domint$, it holds that $G_{\sigma}(s; \bx,\by)=\frac{1}{2\pi\sqrt{\operatorname{det}\bA}}K_0(s\|\bx-\by\|_{\bB})$, where $\|\bx\|_{\bB}^2=\bx^{\top}\bB\bx$.  Like before, we conclude that $|G_{\sigma}(s; \bx,\by)|\leq \frac{1}{2\pi\sqrt{\operatorname{det}\bA}}|K_0(\Re s\|\bx-\by\|_{\bB})|$. 

Using the fact that $\domint=B(0, \radpml)$, denoting $\domint^{\Re s}:=B(0, \radpml\Re s)$ and by re-scaling, we can rewrite 
\begin{align*}
	I_0&\leq (\Re s)^{-4}\frac{1}{4\pi^2\det\bA}\iint_{\domint^{\Re s}\times \domint^{\Re s}}|K_0(\|\bx-\by\|_{\bB})|^2d\bx\,d\by\\
	&\leq (\Re s)^{-4}\frac{1}{4\pi^2\operatorname{det}\bA}\int_{\domint^{2\Re s}}\int_{\domint^{\Re s}}|K_0(\|\bp\|_{\bB})|^2d\by d\bp\leq (\Re s)^{-2}C(\bB, \radpml)\int_{\domint^{2\Re s}}|K_0(\|\bp\|_{\bB})|^2d\bp\\
	&\leq (\Re s)^{-2}C(\bB, \radpml)\int_{\mathbb{R}^2}|K_0(\|\bp\|^2_{\bB})|d\bp. 
\end{align*}
We apply the bound \eqref{eq:K0Z} which allows to conclude that $I_0\leq \tilde{C}(\bB, \radpml)(\Re s)^{-2}$. 
\end{proof}
We immediately obtain from the above the following corollary, which is also the principal result of this section. 
\begin{cor}
	\label{cor:main_stability}
	Consider the PML problem \eqref{eq:rad_time} with the source terms $(f, f_v, \boldsymbol{f}_p, \boldsymbol{f}_q)\in L^2(\setR^+; L^2(\setR^2))$ and vanishing initial conditions. Let the source term $F$, defined after \eqref{eq:lpl_u}, be s.t.\ $F\in H^1(\setR^+; L^2(\setR^2))$, and $F(0)=0$. Let,  additionally, for all $t>0$, $\operatorname{supp}F(t)\subset B(0, r_*)$, with $r_*$ being like in Lemma \ref{lemma:source_far}. 

	Then, for all $T>0$, \eqref{eq:rad_time} admits a unique solution $u^{\sigma}\in H^1(0,T; L^2(\setR^2))$. This solution satisfies the following bound, with the constant $C>0$ depending on $r_*$, $\radpml$ and the matrix $\mathbf{B}$, but independent of $\spml_c$:
	\begin{align*}
		&{\color{black}\|u^{\sigma}\|_{L^{\infty}(0, T; L^2(\mathbb{R}^2))}\leq CT^{1/2}(1+T)\|F\|_{H^1(0, T; L^2(\setR^2))}. }
        \end{align*}
\end{cor}
\begin{proof}
	First of all, as shown in Section \ref{sec:rad_pml_well_posed}, the solution to the PML problem is unique, and its Fourier-Laplace transform is given by, cf. Proposition \ref{prop:fs}, 
	\begin{align}
		\label{eq:representation_usigma}
		\hat{u}^{\sigma}(s,\bx)=\int_{\setR^2}G_{\sigma}(s, \bx, \by)\hat{F}(s,\by)d\by. 
	\end{align}
	We remark that for $\|\bx\|<\radpml$, $\hat{u}^{\sigma}(s,\bx)=\hat{u}(s,\bx)$ (i.e., equal to the solution without the PML), and therefore is $L^2(\domint)$-valued analytic function.
%
As for the $\domext$, let us fix $\bx\in \domext, \, \by\in B(0, r_*)$. Remark that $s\mapsto G_{\sigma}(s; \bx,\by)$ can be extended by analyticity to $\setCp$. 

To deduce the analyticity of $s\mapsto \hat{u}^{\sigma}(s)\in L^2(\setR^2\setminus B(0, \radpml))$ in $\setCp$ it suffices to prove the analyticity of the scalar valued-function $s\mapsto \langle \hat{u}^{\sigma}(s),v\rangle_{L^2(\setR^2\setminus B(0, \radpml))}$ in $\setCp$ for each $v\in L^2(\setR^2\setminus B(0, \radpml))$ \cite[Thm.~1.6.1]{GohbergLeiterer}.
The latter follows from the definition of $\hat{u}^{\sigma}(s)$ in terms of $G_\sigma(s;\bx,\by)$, the analyticity properties of $G_\sigma(s;\bx,\by)$ in $\setCp$ and Lebesgue's dominated convergence theorem (i.e., we can exchange differentiation and integration).

With the Cauchy-Schwarz inequality and Lemma \ref{lemma:source_far}, for all $s\in \setCp$, it holds that  
	\begin{align*}
          \|\hat{u}^{\sigma}\|_{L^2(\setR^2)}\leq C\|\hat{F}\|_{L^2(\setR^2)}\max(1, (\Re s)^{-1}), \quad C>0.
	\end{align*}
It remains to use the Plancherel inequality and the causality argument to obtain the desired stability bound, see \cite[Theorem 4.1]{bkwave}, \cite[Theorem 2.11]{bkw}.
\end{proof}

\subsection{The origin of the instabilities: presence of the essential spectrum for small frequencies}\label{subsubsec:ess-spec}

The goal of this section is to prove Theorem \ref{theorem:fredhomlness}. 
Its proof relies on the following lemma, which shows that the principal symbol of the PDE associated to the frequency-domain PML problem \eqref{eq:rad_lpl_v0} may vanish for some $s_0\in \setCp$ and $\vecx_0\in\domext$. As a corollary, the underlying operator $\oppml$ is no longer Fredholm. 
To prove this result, it is sufficient to show the existence of a singular sequence, i.e., $(v_n)_{n\in\mathbb{N}}\subset H^1(\setR^2)$, s.t. $\|v_n\|_{H^1(\setR^2)}=1$, $(v_n)$ does not have a convergent subsequence, and $\|\oppml(s)v_n\|_{}\rightarrow 0$, see \cite[Thm.\ 3.3]{Halla:19Diss}.
In order to construct such a sequence, we will use (truncated) plane waves $\mathrm{e}^{in\boldsymbol{\xi}_0\cdot\bx}$, with a well-chosen phase $\boldsymbol{\xi}_0$. The choice of the phase is given by the following lemma.

\begin{lem}\label{lem:s0xi0}
Assume that $\spml$ is continuous for $r>\radpml$ and  non-decreasing (and not necessarily constant). Let additionally $\evmax\neq \evmin$. 
Then there exist $\vecx_0\in\domext$, $s_0\in\setCp$ and $\vecxi_0\in\setR^2\setminus\{0\}$ such that $\vecxi_0^\top \Apml(s_0,\vecx_0) \vecxi_0=0$. 
\end{lem}
\begin{proof}
Remark that it is sufficient to show the existence of ${\vecxi}\neq 0$ s.t.\ 
${\vecxi}^\top \bA_{\sigma}^{\phi}(s_0,\vecx_0) {\vecxi}=0$, as $\bA_{\sigma}=\bR_{\phi}\bA_{\sigma}^{\phi}\bR_{\phi}^{\top}$, see Section \ref{sec:pml_intro} for respective definitions of the matrices (and in particular \eqref{eq:rad_lpl_v0}). 
 
Let $\vecx=r(\cos\phi,\sin\phi)^{\top}\in\domext$ be an arbitrary vector (such that $\spml(r), \shpml(r)>0$), and $s\in \setCp$,  both quantities to be fixed later. Following the path of \cite[p.~2721]{Halla:22PMLani}, we express 
$\Apml^{\phi}(s,\vecx)$ via the entries of the matrix $\bA^{\phi}=(A_{ij}^{\phi})_{i,j}$:
\begin{align*}
{\bA}_{\sigma}^{\phi}(s,\bx)=
\begin{pmatrix}
e^{-i\tau} c_d^{-1} A_{11}^{\phi} & A_{12}^{\phi} \\ A_{12}^{\phi} & e^{i\tau} c_d A_{22}^{\phi}	
\end{pmatrix},
\quad \tau:=\arg(\dpml(s,r)/\dtpml(s,r)),\quad 
c_d := |\dpml(s,r)/\dtpml(s,r)|.
\end{align*}
%
%
Recall that $\bA=\diag(\evmax^{-1},\evmin^{-1})$, and that under the condition $\evmax\neq\evmin$, for $\phi \neq \pi n/2, \, n\in \mathbb{Z}$, it holds that $A_{12}^{\phi}\neq 0$.
From now on we restrict ourselves to this range of $\phi$.
Next we plug $\vecxi=(\xi_{1}, \xi_{2})^{\top}\in \setR^2$ into the equation $\vecxi^\top \Apml^{\phi}(s,\bx) \vecxi=0$. The imaginary part of the left-hand side vanishes if $\xi^2_1 = c_d^{2} \xi_{2}^2{A_{22}^{\phi}}/{A_{11}^{\phi}} $.
Plugging in
$\vecxi=(c_d^{1/2}\sqrt{A_{22}^{\phi}}, -c_d^{-1/2}\sqrt{A_{11}^{\phi}}\operatorname{sign}A_{12}^{\phi})^\top$
into $\vecxi^\top \Apml^{\phi}(s,\vecx) \vecxi=0$ yields the identity
	\begin{align}
	\label{eq:cos_condition}
	\cos\tau=|A_{12}^{\phi}|/(A_{11}^{\phi}A_{22}^{\phi})^{1/2}=:a_0.
	\end{align} 
It remains to show that there exist $\bx$ and $s\in \setCp$, s.t.\ the above holds true. This will be the values $\bx_0$ and $s_0$ as in the statement of the lemma. First, as $\bA^{\phi}$ is positive definite, and by the restriction on $\phi$, $0<a_0<1$.
Next, we remark that 
\begin{align*}
	\operatorname{arg}\tau=\operatorname{arg}z(\eta, \omega, r), \qquad  z(\eta,\omega,r)=\eta^2+\omega^2+(\widetilde{\sigma}+\sigma)\eta+\sigma\widetilde{\sigma}+i(\widetilde{\sigma}-\sigma)\omega,
	\end{align*}
where $z$ depends on $r$ through $\sigma,\, \tilde{\sigma}$. 
Let us consider the above expression on the lines $\omega=\operatorname{const}$. For $\eta=0$, $z=\omega^2+i(\widetilde{\sigma}-\sigma)\omega$. Moreover, since $\sigma$ is non-decreasing, we have that $\widetilde{\sigma}(r)\leq  \frac{r-\radpml}{r}\sigma(r)$, and thus $\widetilde{\sigma}-\sigma\neq 0$ for any $r>\radpml$. From this it follows that there exists $\omega=\omega_0$ (sufficiently small) and $r=r_0$, s.t.\ $\cos\operatorname{arg}z(0,\omega_0,r_0)<a_0.$ By the continuity argument, this inequality holds true for sufficiently small $\eta>0$. As $\eta\rightarrow +\infty$, $\cos\operatorname{arg}z(\eta,\omega_0,r_0)\rightarrow 1$. By the mean-value theorem, we conclude that there exists $\eta_0>0$, s.t.\ $\cos\operatorname{arg}z(\eta_0,\omega_0,r_0)=a_0$. Thus, in the statement of the lemma we fix $s_0=\eta_0+i\omega_0$ and choose $r_0,\, \phi_0,\, \vecxi_0$ as discussed above.

\end{proof}
\begin{rmk}
Note that the condition $\vecxi_0^\top \Apml(s_0,\vecx_0)\vecxi_0=0$ with $\vecxi_0\in\setR^2\setminus\{0\}$ is stronger than $\vecxi_0^* \Apml(\vecx_0) \vecxi_0=0$, $\vecxi_0\in\setC^2\setminus\{0\}$.
\end{rmk}

With the help of the above lemma, we can prove Theorem \ref{theorem:fredhomlness}.
\begin{proof}[Proof of Theorem \ref{theorem:fredhomlness}]
Assume that $\oppml$ is Fredholm. 
We are going to derive a contradiction.
  Let $\bx_0$, $s_0$ and $\vecxi_0$ be like in Lemma \ref{lem:s0xi0}. 
For $\delta>0$ let 
\begin{align*}
\bA_{\sigma,\delta}(s,\vecx):=\left\{
	\begin{array}{ll}
	\Apml(s,\vecx), & \|\vecx-\vecx_0\|>\delta, \\
	\Apml(s,\vecx_0), & \|\vecx-\vecx_0\|\leq \delta. 
	\end{array}
	\right.
\end{align*}
%
Let $\sfpmld^s(\cdot,\cdot)$ be defined as $\sfpml^s(\cdot,\cdot)$ with $\Apml$ being replaced by $\Apmld$. 
Since $\Apml(s_0,\vecx)$ is continuous at $\vecx_0$, for a fixed $\epsilon>0$, we can find $\delta>0$ such that $\|\Apml(s_0,.)-\Apmld(s_0,.)\|_{L^\infty(\setR^2)}<\epsilon$.
Hence there exists $\delta>0$ such that the Riesz representation $\oppmld(s_0)$ of $\sfpmld^{s_0}(\cdot,\cdot)$ is Fredholm.
We are going to construct a singular sequence for $\oppmld(s_0)$ in a similar fashion as in the proof of \cite[Thm.~6.2.1]{Agranovich:15}, and thus arrive at a contradiction. 

Let $\chi$ be a $C^2$ cut-off function with $\chi(\vecx_0)=1$, $\supp\chi\subset B_{\delta}(\vecx_0)$. 
Consider the sequence of functions $u_n(\vecx):=\chi(\vecx)\exp(in \vecxi_0\cdot\vecx)$, $n\in\setN$, for which $\|u_n\|_{H^1(\setR^2)}\approx n$.
It is not difficult to check that $u_n/\|u_n\|_{H^1}\rightharpoonup 0$ in $L^2(\setR^2)$, therefore, $u_n/\|u_n\|_{H^1}$ has no convergent in $H^1$ subsequence. 
  Next, let us verify that there exists a constant $C>0$, s.t. $|\sfpmld^{s_0}(u_n,\testf{u})|\leq C\|\testf{u}\|_{H^1(\setR^2)}$ uniformly in $n$.
We compute
\begin{align*}
\innerprod{\Apmld(s_0,\vecx)\nabla u_n}{\nabla \testf{u}}_{L^2}= \innerprod{\Apmld(s_0,\vecx_0)\nabla u_n}{\nabla \testf{u}}_{L^2}
&=\underbrace{\innerprod{\exp(in \vecxi_0\cdot\vecx)\Apml(s_0,\vecx_0)\nabla \chi}{\nabla \testf{u}}_{L^2}}_{I_1^{(n)}} \\
&+\underbrace{\innerprod{in \chi \exp(in \vecxi_0\cdot\vecx)\Apml(s_0,\vecx_0)\vecxi_0}{\nabla \testf{u}}_{L^2}}_{I_2^{(n)}}.
\end{align*}
Clearly, $|I_1^{(n)}|\leq C(\chi,\Apml(s_0,\vecx_0))\|\testf{u}\|_{H^1(\setR^2)}$. It remains to show that the same holds true for  $I_2^{(n)}$. We rewrite it by using the integration by parts
\begin{align*}
I_2^{(n)}
=-\innerprod{(in)^2\chi \exp(in \vecxi_0\cdot\vecx) \vecxi_0^\top \Apml(s_0,\vecx_0)\vecxi_0}{\testf{u}}_{L^2}-\innerprod{in \exp(in \vecxi_0\cdot\vecx) \nabla\chi^\top\Apml(s_0,\vecx_0)\vecxi_0}{\testf{u}}_{L^2}. 
\end{align*}
By Lemma \ref{lem:s0xi0}, the first term in the right-hand side of the above vanishes. To deal with the second term, we rewrite $in\exp(in\vecxi_0\cdot\vecx)\vecxi_0=\nabla \exp(in\vecxi_0\cdot\vecx)$ and integrate by parts again, which yields the following identity:
\begin{align*}
I_2^{(n)}&=-\innerprod{\nabla \exp(in \vecxi_0\cdot\vecx)}{\testf{u} \Apml(s_0,\vecx_0)^*\nabla\chi}_{L^2}\\
&=\innerprod{\exp(in \vecxi_0\cdot\vecx)}{(\nabla \testf{u})^\top \Apml(s_0,\vecx_0)^*\nabla\chi
+\testf{u}\div(\Apml(s_0,\vecx_0)^*\nabla\chi)}_{L^2}.
\end{align*}
The absolute value of the latter expression can be bounded from above by $C(\chi,\Apml(\vecx_0))\|\testf{u}\|_{H^1}$.
Hence $u_n/\|u_n\|_{H^1}$ is a singular sequence for $\oppmld^{s_0}$ and thus $\oppmld(s_0)$ is not Fredholm.
This is a contradiction to the Fredholmness of $\oppml(s_0)$. 
This implies that the essential spectrum of $s\mapsto\oppml(s)$ is non-empty in $\setCp$. 
\end{proof}

\begin{rmk}
Alternatively, one can prove Theorem \ref{theorem:fredhomlness} by constructing a singular sequence based on the localized fundamental solution around the points $\bx$ where its argument vanishes, see Lemma \ref{lem:reh_negative}. 
Yet another approach in the isotropic time-harmonic context is presented in  \cite[Section 5.3.1]{WessDiss} based on a decomposition into polar coordinates and a singular sequence based on spherical Hankel functions with increasing index.
\end{rmk}
\begin{rmk}
  We make the following observation regarding the proof of Theorem \ref{theorem:fredhomlness} and, in particular, Lemma \ref{lem:s0xi0}. As seen from the proof of these results, the support of the functions of the singular sequence is localized in the absorbing layer. It is thus quite natural that they cannot be 'triggered' by the data with a sufficient distance to the absorbing layer, cf. Theorem \ref{theorem:summary}.
	In addition, with increasing index $n\in\setN$ the functions $u_n$ are more and more oscillating, and thus are harder to approximate by the finite element functions. Hence again, it is not surprising that coarse discretizations do not capture the part of the spectrum in $\setCp$, cf., the numerical experiments in Section \ref{sec:first_experiments}. 
\end{rmk}

Now we have proven Theorems \ref{theorem:summary} and \ref{theorem:fredhomlness}. These results show that PMLs for anisotropic problems  in general fail when applied for computational purposes.
Thus in the next section we will seek a remedy for this issue.

\section{A possible remedy: shifted change of variables?}
\label{sec:remedy}
In Section \ref{subsubsec:ess-spec} we have observed that $\oppml$ can loose its Fredholmness when $\cos\arg(\dtpml/\dpml)$ is not bounded far enough away from zero. Indeed, in the isotropic case ($\evmax=\evmin$) the condition \eqref{eq:cos_condition}, sufficient for the loss of Fredholmness, becomes $\cos\arg(\dpml/\dtpml)=0$ (which never holds true for $s\in \setCp$). In the anisotropic case ($\evmin\neq\evmax$) \eqref{eq:cos_condition} can be written as
\begin{align}
	\label{eq:anis_condition}
	\cos\arg(\dpml/\dtpml)=\frac{\left|(\evmax-\evmin)\sin\phi\cos\phi\right|}{(\evmax\cos^2\phi+\evmin\sin^2\phi)^{1/2}(\evmin\cos^2\phi+\evmax\sin^2\phi)^{1/2}}.
\end{align}
Recall that this condition can be fulfilled because $\dpml(s,r)/\dtpml(s,r)\sim (1+\frac{\sigma}{s})$ in the vicinity of the PML interface inside the PML (i.e. when $r=\radpml+\varepsilon$, with $\varepsilon$ being sufficiently small). Then it is always possible to choose $s\in \setCp$, s.t.\  $\arg (1+\frac{\sigma}{s})=\arctan \frac{\sigma \Im s}{|s|^2+\sigma \Re s}$ can take an assigned value $(-\frac{\pi}{2}, \frac{\pi}{2})$, which ensures the validity of the identity \eqref{eq:anis_condition} for a range of $\phi$. 
It is then natural to look for an alternative PML change of variables, e.g., $$r\mapsto r+\frac{\int_{\radpml}^r\sigma(\tilde{r})d\tilde{r}}{s\zeta(s)},$$ for some function $\zeta$ s.t.\ $s\zeta(s)$
is bounded in the vicinity of $0$. Indeed, $\zeta$ should satisfy other properties (in particular  analyticity in $\setCp$) to ensure that the corresponding PML problem remains well-posed. We do not address the question on working out sufficient conditions for general $\zeta(s)$ here, but concentrate on one simple specific choice of $\zeta$.
In the frequency-domain work \cite[p.~2723]{Halla:22PMLani}, it was suggested to choose $\zeta(s)=1+\frac{\gamma}{s}$, for some $\gamma>0$, so that $s\zeta(s)=s+\gamma$. This change of variables amounts to replacing $\dpml(s), \dtpml(s)$ by $\dpml(s+\gamma)$, $\dtpml(s+\gamma)$, respectively. 

  \begin{rmk}
    We note that a similar approach has already been used to stabilize PMLs for isotropic dispersive materials (cf., e.g., \cite{BecacheJolyVinoles,BecacheKachanovska:17, bkw}). In \cite{BecacheKachanovska:17,bkw} $\zeta$ is chosen to also have mathematical properties of physical media (cf., \cite[Section 2.1.1]{bkw}). These assumptions are also fulfilled for the choice $\zeta(s)=1+\gamma s^{-1}$. However, the sufficient assumption on $s\zeta(s)^{-1}$ to ensure convergence of the PML error as $\sigma_c\rightarrow +\infty$ (cf., \cite[Assumption 2.10]{bkw}) is not fulfilled by our choice.
  \end{rmk}
The parameter $\gamma$ is then chosen to ensure that $\cos\operatorname{arg}\frac{\dpml}{\dtpml}$ exceeds the quantity in the right hand side of \eqref{eq:anis_condition} for all $\phi$ and all $s\in \setCp$. Then it is possible to ensure that $\oppml(s)$ is Fredholm for all $s\in\setCp$, in the spirit of \cite[Thm.~3.6]{Halla:22PMLani}. 
%
%
Indeed, since $\cos \operatorname{arg}(z_1/z_2)>\cos\operatorname{arg}z_1$ when $\operatorname{arg}z_1, \operatorname{arg}z_2, \operatorname{arg}z_1-\operatorname{arg}z_2\in (0, \pi/2)$ (resp.\ in $(-\pi/2, 0)$), we have that
\begin{align*}
 	\cos\arg(\dpml(s+\gamma)/\dtpml(s+\gamma))&\geq \cos\arg(\dpml(s+\gamma)) = \cos\arg \Big(1+\frac{\spml_c}{s+\gamma}\Big)\\
	&\geq \cos\arctan \Im\Big(1+\frac{\spml_c}{s+\gamma}\Big) 
	= \cos\arctan \bigg( \frac{\spml_c}{\gamma}\Im\Big( \frac{1}{s\gamma^{-1}+1}\Big) \bigg)
\end{align*}
for all $s\in\setCp$.
Since
\begin{align*}
	\Big|\Im \frac{1}{s\gamma^{-1}+1}\Big|=\frac{\gamma^{-1}|\Im(s)|}{(\Re(s\gamma^{-1})+1)^2+\gamma^{-2}(\Im s)^2}
	\leq \frac{\gamma^{-1}|\Im(s)|}{1+\gamma^{-2}(\Im s)^2} \leq \frac{1}{2}
\end{align*}
for all $s\in\setCp$, we deduce that
\begin{align*}
	\inf_{s\in\setCp} \cos\arg\left(\frac{\dpml(s+\gamma)}{\dtpml(s+\gamma)}\right) &\geq 
	\cos\arctan \Big( \frac{\spml_c}{2\gamma} \Big)
	= \frac{1}{\sqrt{1+(\spml_c/(2\gamma))^2}}.
\end{align*}
Thus by choosing $\spml_c/2\gamma$ small enough we can ensure that $\inf_{s\in\setCp} \cos\arg(\dpml(s+\gamma)/\dtpml(s+\gamma))$ is larger than the right hand side of \eqref{eq:anis_condition} for all $\phi$,
which yields that $\oppml(s)$ is Fredholm for all $s\in\setCp$ \cite[Thm.~3.6]{Halla:22PMLani}.

In order to obtain a more precise stability condition on $\gamma$, we will examine the behaviour of the fundamental solution of the respective problem, in the spirit of Section \ref{sec:fund_sol_pml}.
\begin{rmk}
The idea of using the complex frequency shift in order to stabilize the PMLs for anisotropic media is not new, and, to the best of our knowledge, first appeared in the works of K.\ Duru and G.\ Kreiss \cite{DuruKreiss:12}, \cite{KreissDuru:13}.
However, \emph{complex-frequency shifted} PMLs were introduced already much earlier, e.g., in \cite{RodenGedney:00}. See the discussion in the introduction of \cite{BecachePetropoulosGedney:04} for more details.
\end{rmk}

\subsection{A stability condition on $\gamma$}
Let us introduce the fundamental solution $G_{\sigma,\gamma}$ corresponding to the PML change of variables with $\psi(s)=1+\frac{\gamma}{s}$, analogously to \eqref{eq:def_d_sigma}:
\begin{align}
	\nonumber
	&G_{\sigma,\gamma}(s, \bx, \by):=\frac{1}{2\pi \sqrt{\operatorname{det}\bA}}K_0\left(s\sqrt{(\bx_{\sigma,\gamma}(s)-\by_{\sigma,\gamma}(s))^{\top}\bB(\bx_{\sigma,\gamma}(s)-\by_{\sigma,\gamma}(s))}\right), \quad \bB=\bA^{-1},\\
	\label{eq:def_d_sigma_gamma}
	&	\bx_{\sigma,\gamma}(s):=\hat{\bx}\big(\|\bx\|+\frac{1}{s+\gamma}\int_{\radpml}^{\|\bx\|}\sigma(r')dr'\big)=\bx\left(1+\frac{\sigma_c}{s+\gamma}\frac{\|\bx\|-\radpml}{\|\bx\|}\right)=\bx\, \dtpml(s+\gamma, \|\bx\|), 
\end{align}
and $\by_{\sigma,\gamma}$ is defined similarly to $\bx_{\sigma,\gamma}$. The proof of the counterpart of Proposition \ref{prop:fs} with $\gamma$-shifted change of variables follows the same steps as the proof of Proposition \ref{prop:fs} in Appendix \ref{appendix:Gsigma}, and thus we leave it to the reader. 

To study the behaviour of the fundamental solution, we can proceed as in Section \ref{sec:fund_sol_pml}. Remark that instead of $\hpml(s; \bx,\by)$ as defined in \eqref{eq:defh}, we now study $\hpml(s+\gamma;
\bx,\by)$. Our goal is to choose $\gamma$ so that Lemma \ref{lem:lower_bound} holds true for any source term, and not only for the data supported far enough from the interface.

For this we will require in particular that $\hpml(s+\gamma;\bx,\by)\notin (-\infty, 0]$ for all $s\in \setCp$, \textbf{all} $\by\in \domint$ and $\bx\in \dompml$.  By Lemma \ref{lem:reh_negative}, this requires that one of the following holds true:
\begin{align}
\label{eq:dsgamma}
(a)\,	\Re \dpml(s+\gamma)\neq-\frac{\|\radpml\hat{\bx}-\by\|}{\|\bx\|-\radpml}\frac{\gamma_{12}(\hat{\bx},\by)}{\gamma_{22}(\hat{\bx})}; \quad (b)\,\cos^2(\arg \dpml(s+\gamma))> \frac{\gamma_{12}^2(\hat{\bx},\by)}{\gamma_{11}(\hat{\bx},\by)\gamma_{22}(\hat{\bx})}. 
\end{align} 
Let us define the following ratio which will play an important role in the analysis that follows:
\begin{align}
	\label{eq:nudef}
	\nu:=\frac{\sigma_c}{\gamma}.
\end{align} 
Observe that the image of $\setCp$ under the mapping $\dpml(\cdot+\gamma)$ depends on this ratio only (rather than on $\gamma$, $\sigma_c$ separately). It is the circle centered in the point $\bs_{\nu}=\left(1+\frac{\nu}{2}, 0\right)^\top\in \setR^2$ and of radius $\nu/2$; indeed, 
\begin{align}
	\label{eq:def_boule}
	\begin{split}
B_{\nu/2}(\boldsymbol{s}_{\nu})
	&=\{z\in\mathbb{C}\colon (\Re z-1-\nu/2)^2+(\Im z)^2<\nu^2/4\}
	=\{1+z\colon z\in\mathbb{C}, \, (\Re z-\nu/2)^2+(\Im z)^2<\nu^2/4\}\\
	&=\{1+z\colon z\in\mathbb{C}, \; |z|^2<\nu\Re z\}
	=\{1+z\colon z\in\mathbb{C}, \; \Re (z^{-1}-\nu^{-1})>0\}\\
	&=\{1+z\colon \; z=\frac{1}{s+\nu^{-1}}, \quad s\in \mathbb{C}^+ \}
	=\{1+z\colon \; z=\frac{\sigma_c}{s+\gamma}, \quad s\in \setCp \}.
	\end{split}
\end{align}
Let us now study under which condition on $\nu$ \eqref{eq:dsgamma}(a) or (b) holds true. First of all, the result below shows that the condition \eqref{eq:dsgamma}(a) necessarily fails for some $(\by,\bx)\in\domint\times\dompml$ (see Remark \ref{remark:existence_xy}).
\begin{lem}
Let $\by\in \domint$, $\hat{\bz}\in \mathbb{S}^2$ be s.t. $\gamma_{12}(\hat{\bz},\by)<0$. Let $\nu, \, \gamma>0$. Then, for $s=\gamma$ and $\bx=\rho \hat{\bz}\in \domext$, where $\rho$ is given by
\begin{align*}
	\rho=\radpml-\left(1+\frac{\nu}{2}\right)^{-1}\|\radpml\hat{\bz}-\by\|\frac{\gamma_{12}(\hat{\bz},\by)}{\gamma_{22}(\hat{\bz}, \by)}>\radpml,
\end{align*}
the condition \eqref{eq:dsgamma}(a) \textbf{does not} hold true.
\end{lem}
Again, the proof of this result is left to the reader.
Therefore, in our choice of $\nu$, we will enforce \eqref{eq:dsgamma}(b), for all $\by\in \domint$ and $\hat{\bx}\in \mathbb{S}^2$. Let us define
\begin{align*}
  \beta^2:=\max_{\hat{\bx}\in\mathbb{S}^2,\, \by\in \Omega_{-}(\hat{\bx})}\frac{\gamma_{12}^2(\hat{\bx},\by)}{\gamma_{11}(\hat{\bx},\by)\gamma_{22}(\hat{\bx})},\quad
  \Omega_{-}(\hat{\bx}):=\{\by\in \domint: \gamma_{12}(\hat{\bx},\by)<0\}.
\end{align*}
Then we would like that
\begin{align}
	\label{eq:bcond}
\cos^2(\arg \dpml(s+\gamma))>\beta^2, \text{ for all }s\in \setCp.
\end{align} 
By  Lemma \ref{lemma:fx}, $\beta<1$ and equals
$\beta=\frac{\evmax-\evmin}{\evmax+\evmin}.$
It remains to rewrite \eqref{eq:bcond} as follows, cf.\ \eqref{eq:di_dr}:
\begin{align}
	\label{eq:tancond}
	\tan^2\arg \dpml(s+\gamma)<\frac{1}{\beta^2}-1, \quad s\in \setCp.
\end{align}
In view of the definition of the image of $\dpml(\cdot+\gamma)$, given in \eqref{eq:def_boule},  we have that \footnote{Because the maximum of $\tan \arg  z$ is achieved on the boundary of the circle $\mathcal{B}_{\nu/2}(\bs_{\nu})$, it suffices to find the supremum of the function $f(x)=\frac{\nu^2/4-(x-1-\nu/2)^2}{x^2}$ on the interval $ (1,
1+\nu)$.} $$\operatorname{max}_{z\in \mathcal{B}_{\nu/2}(\bs_{\nu})}\tan^2\arg z=\frac{\nu^2}{4(\nu+1)}.$$ 
Therefore, \eqref{eq:tancond} is equivalent to 
\begin{align}
	\label{eq:stab_cond}
 \frac{\nu^2}{4\nu+4}<\frac{1}{\beta^2}-1\iff \nu<\nu_*:=2\sqrt{\frac{1}{\beta^2}-1}\left(\sqrt{\frac{1}{\beta^2}-1}+\frac{1}{\beta}\right).
\end{align}
Remark that as $\beta\rightarrow 1$ (which amounts to $\evmax/\evmin\rightarrow +\infty$), the above upper bound requires that $\nu\rightarrow 0$.

Our next step would be to show that the condition \eqref{eq:stab_cond} ensures the stability of the PML system. However, this condition per se is not sufficient. Indeed, a complete stability result would also require the analysis of the problem when the sources are located \textbf{inside} the perfectly matched layers. Otherwise we can expect that the instabilities will manifest themselves on the discrete level, just like they did in the case when the source was supported far away from the PMLs, and the corresponding continuous problem was stable. Indeed, recall the numerical experiments of Section  \ref{sec:two_different_types_behaviour}, where the parameters are chosen so that Theorem \ref{theorem:summary} yields the stability of the continuous solution, and yet the discrete solution is unstable.

On the other hand, since such an analysis is fairly technical, and, moreover, it will not allow us to relax \footnote{We will discuss later, in Section \ref{sec:non_convergence}, why this requirement is disadvantageous for the PMLs.} the requirement $\sigma_c/\gamma<\nu_*$, we proceed as follows. We will prove a somewhat weaker stability result, where $\nu_*$ is replaced by a sufficiently small number. Before proceeding, we formulate the time-domain system that we are going to study.

{%
\subsection{The frequency-shifted first-order system}
\label{sec:gamma_system}
Up to this point all our reasoning was solely based on the frequency-domain fundamental solution. The corresponding time-domain system (cf., \eqref{eq:rad_time}) needs yet to be derived. We do this along the lines of the paragraphs preceding \eqref{eq:rad_time} to obtain the frequency shifted first-order in time system \eqref{eq:rad_time_gamma}.
Note that, to this end we merely have to replace the arguments of $\dpml,\tilde\dpml$ by $s+\gamma$ in \eqref{eq:sdd} to obtain
\begin{align*}
		s\tilde \dpml(s+\gamma)\dpml(s+\gamma) &= s +\sigma+\tilde\sigma +\frac{\sigma\tilde\sigma-\gamma(\sigma+\tilde\sigma)}{s+\gamma}-\frac{\gamma\sigma\tilde\sigma}{(s+ \gamma)^2}, \\
  s\frac{\dpml(s+\gamma)}{\tilde \dpml(s+\gamma)} &= s + \sigma-\tilde\sigma -\frac{(\sigma-\tilde\sigma)(\gamma+\tilde\sigma)}{s+\gamma+\tilde \sigma},\\
  s\frac{\tilde \dpml(s+\gamma)}{\dpml(s+\gamma)} &= s - (\sigma-\tilde\sigma) +\frac{(\sigma-\tilde\sigma)(\gamma+\sigma)}{s+ \gamma+\sigma}.
\end{align*}
With $\Jsg(s,\bx)=\Js(s+\gamma,\bx)$, the frequency-shifted counterparts to \eqref{eq:detJ_s} and \eqref{eq:detJJJ_s} are given by 
\begin{align*}
  s \operatorname{det}\Jsg\hat{u}^{\sigma}={}&(s+\sigma+\tilde\sigma)\hat u^\sigma+ \frac{\sigma\tilde\sigma-\gamma(\sigma+\tilde\sigma)}{s+\gamma}\hat u^\sigma-\frac{\gamma\sigma\tilde\sigma}{(s+\gamma)^2}\hat u^\sigma,\\
    s\det\Jsg^{-1}\Jsg^{\top}\bA^{-1}\Jsg\hat\bp^{\sigma} ={}&s\bA^{-1}\hat \bp^{\sigma}+\left(\sigma-\tilde\sigma -\frac{(\sigma-\tilde\sigma)(\gamma+\tilde\sigma)}{s+\tilde \sigma+\gamma}\right)\Px \bA^{-1}\Px\hat\bp^{\sigma},\\
    &+\left(- (\sigma-\tilde\sigma) +\frac{(\sigma-\tilde\sigma)(\gamma+\sigma)}{s+ \sigma+\gamma}\right)\Pxp \bA^{-1}\Pxp\hat\bp^{\sigma}.
\end{align*}
We introduce the auxiliary frequency-domain unknowns
\begin{align*}
  \hat v &:= \frac{1}{\gamma+s}\hat u^\sigma,&
  \hat w &:= \frac{1}{\gamma+s}\hat v,\\
  \hat\bq &:= \frac{\tilde\sigma+\gamma}{s+\gamma+\tilde \sigma}\Px\hat\bp^{\sigma} +\frac{\sigma+\gamma}{s+ \gamma+\sigma}\Pxp \hat\bp^{\sigma},
\end{align*}

This leads to the first order system
\begin{align}
  \begin{split}
    \partial_t u^\sigma &=-(\sigma+\tilde\sigma) u^\sigma + (\gamma(\sigma+\tilde\sigma)-\sigma\tilde\sigma)  v+\gamma\sigma\tilde\sigma w+ \div\bp^{\sigma}+f,\\
  \partial_t v &= u^\sigma-\gamma v, \\
  \partial_t w &= v-\gamma w, \\
  \bA^{-1} \partial_t\bp^{\sigma}&=(\sigma-\tilde\sigma)\left(\Pxp \bA^{-1}\Pxp\bp^{\sigma}-\Px \bA^{-1}\Px\bp^{\sigma}-\Pxp \bA^{-1}\Pxp\bq+\Px \bA^{-1}\Px\bq\right) +\nabla u^\sigma,\\
    \partial_t \bq&=(\tilde\sigma+\gamma)\left(\Px\bp^{\sigma}-\Px\bq\right)+(\sigma+\gamma)\left(\Pxp\bp^{\sigma}-\Pxp \bq\right).
  \end{split}
  \label{eq:rad_time_gamma}
\end{align}
}%
The above problem, equipped with initial conditions, source terms and posed in $\setR^2$, can be proven to be well-posed, following the same arguments as in Section \ref{sec:rad_pml_well_posed}. We will not repeat them here, but rather concentrate on the stability question, which we again study by examining the fundamental solution. 
\begin{rmk}
Remark that compared to the system \eqref{eq:rad_time}, the above system requires the introduction of an extra auxiliary scalar unknown in the PML layer. 	
\end{rmk}
\subsection{Stability of the obtained PML system in the free space}
The principal result of this section is a counterpart of Corollary \ref{cor:main_stability}, and reads. 
\begin{thm}
		\label{thm:stability_shifted_pml}
There exists $\nu_0>0$, s.t. for any $\sigma_c, \gamma>0$, s.t. $\sigma_c/\gamma<\nu_0$, the following holds true. 
Let $F\in H^1(\setR^+; L^2(\setR^2))$ be s.t. $F(0)=0$.  
Then $$\hat{u}^{\sigma}(s,\bx)=\int_{\setR^2}G_{\sigma,\gamma}(s,\bx,\by)\hat{F}(s,\by)d\by$$ is a Fourier-Laplace transform of a  $\mathcal{C}(\setR^+_0; L^2(\setR^2))\cap H^1_{loc}(\setR^+; L^2(\setR^2))$-function $u^{\sigma}$, which additionally satisfies the following stability bound for all $T>0$: 
	\begin{align*}
	\|u^{\sigma}\|_{L^{\infty}(0, T; L^2(\setR^2))}\leq C_1(\bB, \radpml)\mathrm{e}^{C_2(\bB)\sigma_c\radpml}T^{1/2}\max(1, T^2)\|F\|_{H^1(0, T; L^{2}(\setR^2))}, 
	\end{align*}
with some non-negative constants $C_1(\bB, \radpml)$ and $C_2(\bB)$.
\end{thm}
The proof of this theorem can be found in the end of the present section.
\begin{rmk}
We believe that the dependence on $\radpml$ in the above estimate can be waived. As for the dependence on $\sigma_c$, we did not manage to get rid of it.
\end{rmk}
\begin{rmk}
	At a first glance, an exponential increase of the above stability estimate with respect to $\sigma_c$ seems to indicate that it is impossible to ensure convergence of the complex frequency-shifted PMLs as $\sigma_c\rightarrow +\infty$. This is confirmed by the numerical experiments in the section that follows (Section \ref{sec:non_convergence}). Instead, in Section \ref{sec:numerics}, we propose an alternative, truncation-less perfectly matched layer, whose convergence does not rely on increasing $\sigma_c$. Thus, the non-uniformity of the stability estimate in $\sigma_c$ does not seem to pose any problems for the methods that we suggest.  
\end{rmk}

%
%
Our considerations rely on the analysis of the fundamental solution. Remark that the argument of $K_0$ depends now on $(\bx_{\sigma,\gamma}-\by_{\sigma,\gamma})^{\top}\bB(\bx_{\sigma,\gamma}-\by_{\sigma,\gamma})\equiv\hpml(s+\gamma; \bx,\by)$, with $\hpml$ as defined in \eqref{eq:defh}. The key result in the proof of Theorem \ref{thm:stability_shifted_pml} is the following lemma, which is a counterpart of Proposition  \ref{prop:cs_pml_orig_appendix}.
 \begin{prop}
	\label{prop:cs_pml_orig}
	There exist constants $\nu_0,\, C_{\pm},\, c_{\pm} >0$ depending on $\bB$ only, s.t.\ for all $\sigma_c,\gamma>0\colon\, \sigma_c/\gamma<\nu_0$, and for all $s\in \mathbb{C}^+$, {$h_{\sigma}(s+\gamma; \bx,\by)\notin \mathbb{R}_0^{-}$ for $\bx\neq \by$, and }
	the following bounds hold true.

        1. For $(\bx,\by)\in \setR^2\times \setR^2$,
	\begin{align}
		\label{eq:basic_bound}
		C_+|s|\|\bx-\by\|\geq |s\sqrt{\hpml(s+\gamma;\bx,\by)}|\geq C_{-}|s|\|\bx-\by\|.
	\end{align}
	2. There exists $C_{\bB}$ depending on the matrix $\bB$ 
        s.t. for all $(\bx,\by)\in {\overline{\domext}}\times \setR^2$ with  $\|\bx-\by\|>\rho:=C_{\bB}\radpml$, it holds that 
		\begin{align}
			\label{eq:res3}
			\Re \left(s\sqrt{\hpml(s+\gamma;\bx,\by)}\right)\geq c_{+}\Re s\|\bx-\by\|.
		\end{align} 
                3. For $(\bx,\by)\in \overline{\domext}\times\setR^2$, 
		\begin{align}
			\label{eq:reh}
			\Re (s\sqrt{\hpml(s+\gamma; \bx,\by)})\geq
			(c_{+}\Re s-c_{-}\sigma_c)\|\bx-\by\|.
		\end{align}
		Moreover, $s\mapsto G_{\sigma,\gamma}(s; \bx,\by)$ is analytic in $\mathbb{C}^+$ for all $\bx\neq \by$. 
	\end{prop}
\begin{rmk}
	In the above, we believe that the dependence of $\rho$ on $\radpml$ can be waived, and is due to the non-optimality of the estimates.
\end{rmk}
For the proof of the above proposition please see Appendix \ref{appendix:h_behave}, where this proposition is re-stated as Proposition \ref{prop:cs_pml_orig_appendix} and Corollary \ref{cor:Gsigma}. The above lemma allows to obtain Laplace-domain bounds on the solution of the shifted PML problem.
\begin{prop}
	\label{prop:lpl_domain}
Assume that $g\in L^2(\setR^2)$. Then, with $\nu_0>0$ like in Lemma \ref{prop:cs_pml_orig}, and for all $\sigma_c, \gamma>0$, s.t. {$\sigma_c/\gamma<\nu_0$}, the function $v(s, \bx)=\int_{\setR^2}G_{\sigma,\gamma}(s; \bx,\by)g(\by)d\by$ satisfies the following bound:
	\begin{align*}
	\|v(s)\|_{L^2(\setR^2)}\leq C_1(\bB, \radpml)\mathrm{e}^{C_2(\bB)\sigma_c \radpml}\max(1, (\Re s)^{-2})\|g\|_{L^2(\setR^2)},
	\end{align*}
with some non-negative constants $C_1(\bB, \radpml)$ and $C_2(\bB)$.
\end{prop}
\begin{proof}
	Let us consider 
	\begin{align*}
		\|v(s)\|_{L^2(\setR^2)}^2=\int_{\setR^2}\left|
		\int_{\setR^2}G_{\sigma,\gamma}(s; \bx,\by)g(\by)d\by\right|^2 d\bx.
	\end{align*}
	We make use of the explicit expression of the fundamental solution \eqref{eq:def_d_sigma_gamma} and the bound \eqref{eq:K0Z}, which we repeat for the convenience of the reader below. For any $a>0$, there exists $C(a)>0$, s.t., for $z\in \mathbb{C}\setminus (-\infty, 0]$, 
	\begin{align}
		\label{eq:K0Z2}
		|K_0(z)|\leq C(a)\left\{
		\begin{array}{ll}
			\max(1,|\log|z||), & |z|\leq a, \\
			\mathrm{e}^{-\Re z}, & |z|>a.
		\end{array}
		\right.
	\end{align}
	For a fixed $\bx\in \setR^2$, we split, based on the elements of Lemma \ref{prop:cs_pml_orig},
	\begin{align}
		\label{eq:splitting}
		&\setR^2=\mathcal{O}_1(\bx)\cup\mathcal{O}_2(\bx)\cup\mathcal{O}_3(\bx), \\
		\nonumber
		&\mathcal{O}_1(\bx):=\{\by\in \setR^2\colon\, |s|\|\bx-\by\|\leq 1\}, \\
		\nonumber
		&\mathcal{O}_2(\bx):=\{\by\in \setR^2\colon\, |s|\|\bx-\by\|> 1, \quad \|\bx-\by\|\leq\rho\}, \\
		\nonumber
		&\mathcal{O}_3(\bx):=\{\by\in \setR^2\colon\, |s|\|\bx-\by\|> 1, \quad \|\bx-\by\|>\rho\}.
	\end{align}
	Let us now consider $G_{\sigma,\gamma}(s; \bx,\by)$ in each of these regions. 
	In $\mathcal{O}_1$, by the upper bound in \eqref{eq:basic_bound}, $|s\sqrt{\hpml(s+\gamma; \bx,\by)}|$ $\leq$ $C_{+}$, and therefore, by \eqref{eq:K0Z2}, 
	\begin{align}
		\label{eq:c1}
		|G_{\sigma,\gamma}(s; \bx,\by)|\leq C(C_+)\max(1, |s\sqrt{\hpml(s+\gamma; \bx,\by)}|)\leq C_1\left(1+\left|\log(|s|\|\bx-\by\|)\right|\right),
	\end{align}
for some $C_1>0$, and 
	where in the last bound we exploited the lower bound in \eqref{eq:basic_bound}.
	In $\mathcal{O}_2$,  using the lower bound \eqref{eq:basic_bound} and the definition of $\mathcal{O}_2$,  $\left|s\sqrt{\hpml(s+\gamma;\bx,\by)}\right|>C_{-}$, therefore, we employ \eqref{eq:K0Z2} to show that 
	\begin{align*}
		|G_{\sigma,\gamma}(s; \bx,\by)|\leq C(C_{-})\mathrm{e}^{-\Re (s\sqrt{\hpml(s+\gamma; \bx,\by)})}.
	\end{align*}
	We then use \eqref{eq:reh} to bound further
	\begin{align*}
		|G_{\sigma,\gamma}(s; \bx,\by)|\leq C(C_{-})\mathrm{e}^{-c_+\Re s\|\bx-\by\|+c_{-}\sigma_c\|\bx-\by\|}\leq C_2\mathrm{e}^{c_{-}\sigma_c\rho}, \quad C_2>0.
	\end{align*}
	Finally, in $\mathcal{O}_3$, we have again that $\left|s\sqrt{\hpml(s+\gamma;\bx,\by)}\right|\geq C_{-}$, and thus use the bound of \eqref{eq:K0Z2}, combined with the bound \eqref{eq:res3}, which yields 
	\begin{align*}
		|G_{\sigma,\gamma}(s; \bx,\by)|\leq C_3\mathrm{e}^{-c_{+}\Re s\|\bx-\by\|}, \quad C_3=C(C_{-})>0.
	\end{align*}
	The splitting \eqref{eq:splitting} and the above bounds allow us to rewrite 
	\begin{align}
		\nonumber
		&\|v(s)\|_{L^2(\setR^2)}^2\leq 3\sum\limits_{i=1}^3S_i, \quad S_i=\int_{\setR^2}\left|
		\int_{\mathcal{O}_i(\bx)}G_{\sigma,\gamma}(s; \bx,\by)g(\by)d\by\right|^2 d\bx\\
		\label{eq:bsi}
		&\leq \int_{\setR^2}\left|
		\int_{\setR^2}H_{i}(\bx-\by)|g(\by)|d\by\right|^2 d\bx,\\
		\nonumber
		&H_{1}(\bx)=C_1\left(1-\log\left(|s|\|\bx\|\right)\right)\mathbbm{1}_{|s|\|\bx\|\leq 1},\quad 
		H_2(\bx)=C_2\mathrm{e}^{c_{-}\sigma_c\rho}\mathbbm{1}_{1< |s|\|\bx\|\leq |s|\rho},\\
		\nonumber 
		& H_{3}(\bx)=C_3\mathrm{e}^{-\Re s\|\bx\|}\mathbbm{1}_{|s|\|\bx\|>\max(1,|s|\rho)}.
	\end{align}
	Next we recognize in the bound \eqref{eq:bsi} for $S_i$ the $L^2$-norm of a convolution product. With Young's inequality for convolutions we obtain $S_i\lesssim \|H_i\|_{L^1}^2\|g\|_{L^2}^2$. We end up with the following inequality:
	\begin{align*}
		\|v(s)\|^2_{L^2(\setR^2)}\leq  \|g\|^2_{L^2(\setR^2)}\sum\limits_{i=1}^3 \|H_i\|^2_{L^1(\setR^2)}.	
	\end{align*}
	It remains to estimate the $L^1$-norms of $H_i$, $1\leq i\leq 3$. This yields
	\begin{align*}
		\|H_1\|_{L^1(\setR^2)}&=C_1\int_{B_{\frac{1}{|s|}}}(1-\log\left(|s|\|\bx\|\right))d\bx\lesssim \left(|s|^{-2}-|s|^{-2}\int_{B_1}\log\|\bx\|d\bx\right)\\
		&\lesssim \max(1, (\Re s)^{-2}),\\
		\|H_2\|_{L^1(\setR^2)}&\leq \int_{B_{\rho}}\mathrm{e}^{c_{-}\sigma_c\rho}d\bx\lesssim \rho^2\mathrm{e}^{c_{-}\sigma_c\rho}, \\
		\|H_3\|_{L^1(\setR^2)}&\lesssim \int_{\mathbb{R}^2}\mathrm{e}^{-c_{-}\Re s\|\bx\|}d\bx\lesssim \max((\Re s)^{-2}, 1).
	\end{align*}
	Recalling the definition of $\rho$ from Lemma \ref{prop:cs_pml_orig} allows to obtain the desired estimate. 
\end{proof}

We have the necessary ingredients to prove Theorem \ref{thm:stability_shifted_pml}. 
\begin{proof}[Proof of Theorem \ref{thm:stability_shifted_pml}]
\textbf{Analyticity of $\hat{u}^{\sigma}$. }As discussed in the proof of Corollary \ref{cor:main_stability}, to prove $L^2$-analyticity of $\hat{u}^{\sigma}(s,\bx)$, it is sufficient to show that for all $\varphi\in L^2(\mathbb{R}^2)$, the function 
\begin{align*}
	s\mapsto f_{\varphi}(s):=(\hat{u}^{\sigma}(s,.), \varphi)
\end{align*}
is analytic (in $\mathbb{C}^+$). To see that this indeed holds true it is sufficient to check that on each compact $K\subset \mathbb{C}^+$, $f_{\varphi}$ is a uniform limit of analytic functions
\begin{align*}
f_{\varphi}^{\varepsilon}(s):=\int_{\mathbb{R}^2\times \mathbb{R}^2: \|\bx-\by\|>\varepsilon}G_{\sigma}(s; \bx, \by)\hat{F}(s,\by)\varphi(\bx)d\bx\,d\by, 
\end{align*}
which will allow to conclude by \cite[Thm.~1.6.1]{GohbergLeiterer}. 

\textit{Analyticity of $f_{\varphi}^{\varepsilon}$. }To verify that $f_{\varphi}^{\varepsilon}(s)$ is analytic, we remark that $s\mapsto \hat{F}(s, .)$ is $L^2(\mathbb{R}^2)$-analytic in $\mathbb{C}^+$ (as the Laplace transform of $F$); on the other hand, $G_{\sigma}(s; \bx,\by)$ is analytic in $\mathbb{C}^+$ for each $(\bx,\by)\in \{\mathbb{R}^2\times \mathbb{R}^2: \|\bx-\by\|>\varepsilon\}$, as follows from Corollary \ref{cor:Gsigma}. Moreover, with \cite[10.29.3]{nist},  $$\partial_s G_{\sigma,\gamma}(s; \bx,\by)=-\frac{1}{2\pi \sqrt{\det\mathbf{A}}}K_1(s\sqrt{h_{\sigma}(s+\gamma; \bx,\by)})\left(\sqrt{h_{\sigma}(s+\gamma;\bx,\by)}+\frac{s}{\sqrt{h_{\sigma}(s+\gamma; \bx,\by)}}\partial_s h_{\sigma}(s+\gamma; \bx,\by)\right).$$ 

For any $a>0$, there exists $C(a)>0$, s.t., for $z\in \mathbb{C}^+$, 
\begin{align*}
	|K_1(z)|\leq C(a)\left\{
	\begin{array}{ll}
	|z|^{-1}, & |z|\leq a,\\
	\mathrm{e}^{-\Re z}, & |z|>a, 	
	\end{array}
	\right.
\end{align*} 
cf. \cite[10.30, 10.40]{nist}. Thus, to bound $\partial_s G_{\sigma, \gamma}$ we can pursue the same strategy as in the proof of Proposition \ref{prop:lpl_domain} (remark that we consider the case when $\|\bx-\by\|\geq \varepsilon$), cf. e.g. the bound \eqref{eq:K0Z2}, and next argue that $|\partial_s G_{\sigma,\gamma}(s; \bx,\by)\hat{F}(s,\by) \varphi(\bx)|+| G_{\sigma,\gamma}(s; \bx,\by)\partial_s\hat{F}(s,\by) \varphi(\bx)|$ is bounded by an $L^1$-function, in order to apply the Lebesgue's dominated convergence theorem (cf. the proof of Proposition \ref{prop:lpl_domain}). This would prove that for all $\varepsilon>0$, $f_{\varphi}^{\varepsilon}$ is analytic in $\mathbb{C}^+$. 

\textit{Uniform convergence. }It remains then to argue that for any $K\subset \mathbb{C}^+$, $\sup_{s\in K}|f_{\varphi}(s)-f_{\varphi}^{\varepsilon}(s)|\rightarrow 0$, as $\varepsilon\rightarrow 0$. For this we consider 
\begin{align*}
f_{\varphi}(s)-f_{\varphi}^{\varepsilon}(s)=\int_{\mathbb{R}^2\times \mathbb{R}^2: \|\bx-\by\|<\varepsilon}G_{\sigma,\gamma}(s; \bx, \by)\hat{F}(s,\by)\varphi(\bx)d\bx\,d\by.
\end{align*}
Let us start by bounding $|G_{\sigma,\gamma}(s; \bx,\by)|$. 
For $s$ belonging to a compact subset of $K$, we use the bound \eqref{eq:c1}, which yields 
\begin{align*}
	|G_{\sigma,\gamma}(s; \bx,\by)|\lesssim C_K |\log\|\bx-\by\| |, \quad \|\bx-\by\|<\varepsilon.
\end{align*}
Therefore, for all $s\in K$, 
\begin{align*}
	|f_{\varphi}(s)-f_{\varphi}^{\varepsilon}(s)|&\leq C_K\int_{\mathbb{R}^2\times \mathbb{R}^2}|\log\|\bx-\by\||\, \mathbbm{1}_{\|\bx-\by\|<\varepsilon}|\hat{F}(s,\by)|\, |\varphi(\bx)|d\bx\,d\by\\
	&\leq \|\varphi\|_{L^2(\mathbb{R}^2)}\|\int_{\mathbb{R}^2}\mathcal{G}^{\varepsilon}(\|\bx-\by\|)|\hat{F}(s,\by)|d\by\|_{L^2(\mathbb{R}^2)},
\end{align*}
with $\mathcal{G}^{\varepsilon}(\|\bx\|)=\log\|\bx\|\mathbbm{1}_{\|\bx\|<\varepsilon}$. With the Young inequality for convolutions , we conclude that 
\begin{align*}
\sup_{s\in K}|f_{\varphi}-f_{\varphi}^{\varepsilon}|&\leq\sup_{s\in K}\|\hat{F}(s)\|_{L^2(\mathbb{R}^2)}\|\varphi\|_{L^2(\mathbb{R}^2)}\|\mathcal{G}^{\varepsilon}\|_{L^1(\mathbb{R}^2)}\rightarrow 0, \quad \varepsilon\rightarrow 0.
\end{align*}		
\textbf{Proof of the time-domain bound. }Use the bound of Proposition \ref{prop:lpl_domain}, and the Plancherel theorem in the spirit of \cite{bkwave, bkw}. 
\end{proof}
\subsection{Investigation of convergence}
\label{sec:non_convergence}
Recall that the convergence of the PMLs is quantified by two parameters: the absorption parameter $\sigma$ and the width of the PML layer $L$, see \cite{diaz_joly,bkw,bkwave}. In the classical, B\'erenger's case, both in waveguides and in the free space,  the error between the PML solution and the exact solution on the time interval $(0,T)$ decreases at least exponentially fast as $\sigma_c\rightarrow +\infty$ and $L$ is fixed, or as $L\rightarrow +\infty$ and $\sigma_c$ is fixed  \cite{diaz_joly, bkwave, bkw}. 
Since from the computational viewpoint varying $L$ can be quite expensive,  much effort in recent years was dedicated to the choice of the optimal profiles of the damping functions $\sigma$ \cite{collino_monk_optimizing, axel:optimizing}.
In this section we present the analysis for $\sigma$ being piecewise constant (cf. Assumption \ref{assump:piecewise_const}).

Unfortunately, \textbf{when $\gamma$ is chosen proportionally to $\sigma_c$, it seems impossible to ensure the PML convergence as $\sigma_c\rightarrow +\infty$ and the PML width $L$ is kept constant}.

In the section that follows we present an analytic argument which will show that the best error convergence we can hope for is $O(\sigma_c^{-1/2})$, and in the next section we will present a numerical experiment where the error stagnates as $\sigma_c\rightarrow +\infty$, and no convergence is observed.
\begin{rmk}
		\label{rmk:freq_domain}
	The lack of convergence is easy to understand when considering the time-harmonic regime ($s=-i\omega$). Indeed, in one dimension, under the PML change of variables $x\mapsto x+\frac{\int_{\radpml}^{x}\sigma(x)dx}{-i\omega+\gamma}$, the time-harmonic outgoing plane wave $u_{pw}(t,x)=\mathrm{e}^{-i\omega t+ikx}$ is transformed into an evanescent wave
	\begin{align*}
		u_{pw}^{\sigma}(t,x)=u_{pw}(t,x)\lambda_{\operatorname{att}}(x,\omega,k), \quad \lambda_{\operatorname{att}}(x,\omega,k)=\mathrm{e}^{-\frac{k \sigma_c (x-\radpml)}{\omega+i\gamma}}.
	\end{align*}
For $\gamma=\nu\sigma_c$, the attenuation factor is then given by $|\lambda_{\operatorname{att}}(x,\omega,k)|=\mathrm{e}^{-\frac{k\omega\sigma_c(x-\radpml)}{\omega^2+\nu^2\sigma_c^2}}$. In particular, for $\omega=k$, we have that 
\begin{align}
	\label{eq:lambdaatt}
		|\lambda_{\operatorname{att}}(x,k,k)|\sim \mathrm{e}^{-\frac{k^2}{\nu^2\sigma_c}(x-\radpml)}, \text{ as }\sigma_c\rightarrow +\infty.
\end{align}
The above shows that the convergence can be assured only by increasing the width of the perfectly matched layer, proportionally to $\sigma_c$. 
\end{rmk}

\subsubsection{Analysis}
In the 1-dimensional case it is possible to obtain an explicit expression of the PML solution. Indeed, consider the following model problem:
\begin{align*}
	&(\partial_t^2-\partial_x^2)u=0,\quad  x>0, \\
	&\left. u\right|_{x=0}=g(t), \quad +\text{ vanishing i.c.}
\end{align*}
%
We apply the $\gamma$-shifted PMLs stemming from the change of variables with $\sigma_c, \gamma>0$, 
\begin{align*}
x_{\sigma}(x)=\left\{
\begin{array}{ll}
x+\frac{\sigma_c(x-\radpml)}{s+\gamma}, & \radpml\leq x<\radpml+L,\\
x, & x<\radpml.
\end{array}
\right.	
\end{align*}
We truncate the PML with the Dirichlet boundary conditions at $x=\radpml+L$. 
The time-domain error between the exact solution and the solution to the time-domain counterpart of the resulting PML problem can be expressed in terms of the series of modified Bessel functions (see Appendix \ref{appendix:error}). We have the following explicit expression of the error (where the series is finite for each final time $t>0$):
\begin{align}
	\label{eq:esigmat}
	\begin{split}
		e^{\sigma}(t,x)&=u^{\sigma}(t,x)-u(t,x)=\mathrm{e}^{-\gamma t}\sum\limits_{\ell=1}^{\infty}\mathrm{e}^{-2\sigma_c L\ell}T_{\ell}g+\sum\limits_{\ell=1}^{\infty}\mathrm{e}^{-2\sigma_c L\ell}\alpha^{1/2}_{\ell}\int_0^{t}\mathrm{e}^{-\gamma \tau} I_1(2\alpha_{\ell}\tau^{1/2})\tau^{-1/2}T_{\ell}g(t-\tau,x)d\tau.
	\end{split}
\end{align}
In the above $I_1$ is a modified Bessel function, and the operators $T_{\ell}$ and the coefficient $\alpha_{\ell}$ are defined by 
\begin{align*}
	T_{\ell}g(t,x)=g(t- x-2 (\radpml+L)\ell)-g(t+ x-2 (\radpml+L)\ell), 	\quad \alpha_{\ell}=\sqrt{2L\gamma\sigma_c \ell}.
\end{align*}
We are interested in asymptotics of the above expression as $\sigma_c\rightarrow +\infty$, and the ratio $\sigma_c/\gamma$ is kept constant equal to $\nu$. Assume w.l.o.g. the following: 
\begin{align}
	\label{eq:assumptions1}
	&(a)\; \nu\geq 1; \quad (b)\;\sqrt{2L\nu}\in (\radpml/8, 3\radpml/8);\\
		\label{eq:assumptions2}
 &(c)\; g(t)\geq 0 \text{ on }(0, \radpml/2) \text{ and }g(t)>C_g>0 \text{ on }(\radpml/8, 3\radpml/8).
\end{align}

Measuring the error at the time $t=t_*=2\radpml+2L$ at the point $x=x_*=\radpml/2$ yields an explicit expression 
\begin{align*}
	e^{\sigma}(t_*, x_*)&=\mathrm{e}^{-\gamma t_*}\mathrm{e}^{-2\sigma_c L}g\left(\frac{\radpml}{2}\right)+\int_0^{\radpml/2}k(\tau)g\left(\frac{\radpml}{2}-\tau\right)d\tau,\\
	k(\tau) &=\mathrm{e}^{-2\sigma_c L}\alpha_1^{1/2}\mathrm{e}^{-\gamma \tau}I_1(2\alpha_1\tau^{1/2})\tau^{-1/2}, \quad \alpha_1=\sigma_c\sqrt{2L\nu^{-1}}.
\end{align*}
With \eqref{eq:assumptions2}, we can write 
\begin{align}
	\label{eq:esigma}
		e^{\sigma}(t_*, x_*)\geq C_g \int_{\radpml/8}^{3\radpml/8}k(\tau)d\tau.
\end{align}
 Using the asymptotic behaviour of modified Bessel functions $I_1(z)\sim (2\pi z)^{-1/2}\mathrm{e}^z$, as $z\rightarrow +\infty$, cf. \cite[10.40.1]{nist}, we arrive at the following expression for $\tau$ bounded away from the origin: 
\begin{align*}
	k(\tau)\sim (4\pi)^{-1}\mathrm{e}^{-\sigma_c\nu^{-1}(\sqrt{2L\nu}-\sqrt{\tau})^2}\tau^{-3/4}, \quad \sigma_c\rightarrow +\infty.
\end{align*}
Plugging the above in \eqref{eq:esigma} yields,
\begin{align*}
		e^{\sigma}(t_*, x_*)\gtrsim  \int_{\radpml/8}^{3\radpml/8}\mathrm{e}^{-\sigma_c\nu^{-1}(\sqrt{2L\nu}-\sqrt{\tau})^2}\tau^{-3/4}d\tau\gtrsim \int_{\radpml/8}^{3\radpml/8}\mathrm{e}^{-\sigma_c\nu^{-1}(\sqrt{2L\nu}-p)^2}dp,
\end{align*}
where the hidden constant depends on $\radpml$. Finally, by assumption \eqref{eq:assumptions1}, making $\sigma_c\rightarrow +\infty$, we can bound from below
\begin{align*}
e^{\sigma}(t_*,x_*)\gtrsim \int_{-\infty}^{\infty}\mathrm{e}^{-\sigma_c \nu^{-1} p^2}dp=O(\sigma_c^{-1/2}), \quad \sigma_c\rightarrow +\infty.
\end{align*}
%
%


%
The above example indicates that, in general, we can expect only algebraic convergence of the complex-scaled PMLs with respect to $\sigma_c\rightarrow +\infty$. Such a slow convergence in $\sigma_c^{-1/2}$ is of limited practical interest: it had been observed \cite{collino_monk_optimizing,baffet_grote_et_al} that using large values of $\sigma_c$ in practice requires the use of fine discretizations or high-order finite elements to prevent numerical reflections from the interface between the absorbing layer and the physical media.  Moreover, the numerical experiment of the following section shows that even this estimate is quite optimistic, and the error may stagnate for $\sigma_c$ being sufficiently large.

\subsubsection{Numerical experiments and conclusions}
\label{sec:numerics_cfs}
\begin{figure}
	\centering
	\includegraphics*[width=0.4\textwidth]{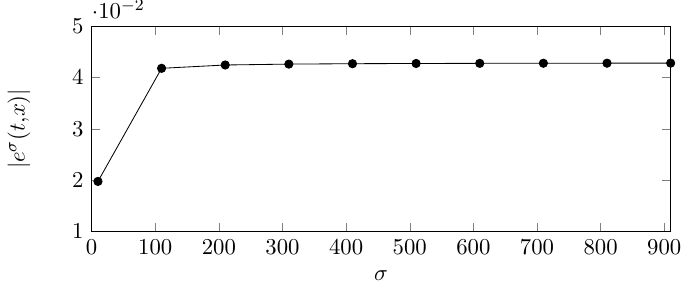}
	\caption{Illustration to the experiment of Section \ref{sec:numerics_cfs}. Dependence of $|e^{\sigma_c}(t,x)|$ on $\sigma_c$. }
	\label{fig:error_cfs}
\end{figure}

To support the considerations of the previous section, let us plot the dependence of the PML error on $\sigma_c$, when keeping $\gamma=\nu^{-1}\sigma_c$ with $\nu$ being constant.  We compute this error numerically using its expression in the Laplace domain \eqref{eq:error_lpl_1d} and the convolution quadrature method \cite{MR0923707}. 
We fix $L=1$, $\radpml=0.2$, $x=0.1$, $\nu=1$, the final time $t=10$, choose $g(t)=\mathrm{e}^{-t^2}t^2$ and compute $|e^{\sigma_c}(t,x)|$ depending on $\sigma_c$. The corresponding value is shown in Figure \ref{fig:error_cfs}.
%
The results of Figure \ref{fig:error_cfs} suggest that the corresponding PMLs do not converge as $\sigma_c\rightarrow +\infty$, and the only way the convergence can be assured is by taking $L\rightarrow +\infty$. Unfortunately, in this case the PMLs are hardly more efficient than a classical super-cell method: indeed, if one is interested in the solution inside the box $(-a,a)\times (-a,a)$, due to the finite velocity $v=1$ of the wave propagation, for a fixed time $t>0$, one can always compute the solution inside the larger box $(-a-L, a+L)\times (-a-L, a+L)$, with $L=t/2$ and homogeneous Dirichlet boundary conditions at the boundary of the box. And thus the improvement provided by the PMLs would be only marginal. 

Therefore, in the section that follows we suggest two alternative strategies. The first one is based on the use of the infinite elements (Hardy-space methods) for the $\gamma$-shifted PML system.
The second one is based on mapping of the exterior domain to a bounded domain, which is applied after the complex scaling.
The advantage of both approaches is that no PML truncation is done, and thus at the continuous level $\gamma$-shifted PMLs become exact.


\section{Two different numerical methods}
\label{sec:numerics}
In the following we present two strategies to discretize the system \eqref{eq:rad_time_gamma} without truncating the exterior domain, to preserve stability of the system and still obtain a converging method.
Both of the methods are Galerkin methods based on the following semi-discrete weak formulation of \eqref{eq:rad_time_gamma} to find $u^\sigma,v,w,\bp^{\sigma},\bq\in C^1([0,T],\scalsp)^3\times C^1([0,T],\vecsp)^2$ such that
  \begin{align}
    \begin{split}
      \langle\partial_t u^\sigma,u^\dagger\rangle &=\langle(-(\sigma+\tilde\sigma) u^\sigma +(\gamma(\sigma+\tilde\sigma)-\sigma\tilde\sigma)v+\gamma\sigma\tilde\sigma w,u^\dagger\rangle-\langle\bp^{\sigma},\nabla u^\dagger\rangle+\langle f, u^\dagger\rangle,\\
      \langle\partial_t v,v^\dagger\rangle &= \langle u^\sigma-\gamma v,v^\dagger\rangle,\\
      \langle\partial_t w,w^\dagger\rangle &= \langle v-\gamma w,w^\dagger\rangle,\\
      \langle\bA^{-1} \partial_t\bp^{\sigma},\bp^{\dagger}\rangle&=\langle(\sigma-\tilde\sigma)\bA^{-1}(\bp^{\sigma}_\tangsub-\bq_\tangsub),\bp_\tangsub^\dagger\rangle-\langle(\sigma-\tilde\sigma)\bA^{-1}(\bp^{\sigma}_\radsub-\bq_\radsub),\bp_\radsub^\dagger\rangle+\langle\nabla u^\sigma,\bp^{\dagger}\rangle,\\
      \langle\partial_t \bq,\bq^\dagger\rangle&=\langle(\tilde\sigma+\gamma)\left(\bp^{\sigma}_\radsub-\bq_\radsub\right),\bq_\radsub^\dagger\rangle+\langle(\sigma+\gamma)\left(\bp^{\sigma}_\tangsub-\bq_\tangsub\right),\bq_\tangsub^\dagger\rangle,
    \end{split}
    \label{eq:rad_time_weak}
  \end{align}
  for all $u^\dagger,v^\dagger,w^\dagger\in\scalsp,\bp^\dagger,\bq^\dagger\in\vecsp$
  and some (discrete) spaces $\scalsp\subset H^1(\dom), \vecsp\subset \left(L^2(\dom)\right)^2,$
  where the radial and tangential components of vector unknowns are denoted by subscripts $\radsub,\tangsub$ (i.e., $\vecs f_\radsub:=\Px \vecs f,\quad\vecs f_\tangsub := \Pxp \vecs f$, with $\Px$ being the projection onto the space spanned by $\bx$). The above system is equipped with vanishing initial conditions. 
  Subsequently we apply the Crank-Nicholson time-stepping method to the above semi-discrete system.
\subsection{Method 1: infinite elements}
\label{sec:numerics_ie}
As a first approach to construct the discrete spaces $\scalsp,\vecsp$ we use infinite elements.
To be more specific, we use Hardy space infinite elements (HSIEs) \cite{hsm,Halla:16}.
Note that so far HSIEs have been primarily used for time-harmonic problems, with the exception of \cite{RuprechtSchaedleSchmidt:13}, \cite[10.2.4]{WessDiss}, \cite{td_ie}.
Since HSIEs are constructed by means of a rather technical apparatus working in the Laplace domain, we choose to work here with an equivalent, but more accessible presentation as in \cite{NW22,WessDiss}.
That is we directly specify the infinite elements as discrete subspaces of $H^1(\domext)$ and $\left(L^2(\domext)\right)^2 $.
In $\domint$ we use conforming finite element spaces $\scalspint\subset H^1(\domint), \vecspint\subset\big(L^2(\domint)\big)^2$. In our particular implementation, we use Lagrange $\mathbb{P}_k$ finite elements for discretizing $u^{\sigma}$, and discontinuous Lagrange elements $\left(\mathbb{P}_{k-1}\operatorname{-DG}\right)^2$ for discretizing $\bp^{\sigma}$. In $\domext$ we use tensor products of the the trace spaces of the interior spaces and (scalar) radial spaces, i.e.,
\begin{align*}
  \scalspext&:=\scalspextrad\otimes\scalspexttang=\spa\{\vecx=r\hat\vecx\mapsto \tilde v(r) \hat v(\radpml\hat\vecx)\colon\; \tilde v\in\scalspextrad, \;\hat v\in\scalspexttang\},\\
  \vecspext&:=\vecspextrad \otimes \vecspexttang = \spa\{\vecx=r\hat\vecx\mapsto \tilde q(r) \hat \bq(\radpml\hat\vecx)\colon\; \tilde q\in\vecspextrad, \;\hat \bq\in\vecspexttang\},
\end{align*}
where
\begin{align*}
  \scalspexttang &:= \{v|_\bdpml,v\in \scalspint\},&
  \vecspexttang &:= \{\bq|_\bdpml,\bq\in \vecspint\}.
\end{align*}
In the latter case, as the elements of $\vecspint$ are piecewise-polynomials, the trace $\bq|_\bdpml$ is well-defined on each mesh element adjacent to $\bdpml$.
The spaces for discretizing the radial part are then defined as follows.
To capture the combination of the exponential decay and the oscillatory behaviour of frequency-domain PML solutions, cf. Remark \ref{rmk:freq_domain}, we use products of exponentially decaying functions and polynomials
\begin{align}
  \label{eq:radspace}
  	\scalspextrad(\iepar_0,\iepar_1,N):=\vecspextrad(\iepar_0,\iepar_1,N):=\spa\left\{\exp(-\iepar_0 r)p(r),\, \exp(-\iepar_1 r)p(r): p\in\mathcal P^N\right\}, 
\end{align}
with $N\in \mathbb{N}, \, \eta_0\neq \eta_1>0.$ 
In the above $\mathcal P^N$ denotes the space of polynomials of degree $N$, and thus the respective spaces are $2(N+1)$-dimensional.  For a particular case $\iepar_0=\iepar_1$, we define the $2(N+1)$-dimensional space
  \begin{align}	\scalspextrad(\iepar_0,\iepar_0,N):=\vecspextrad(\iepar_0,\iepar_0,N):=\spa\left\{\exp(-\iepar_0 r)p(r)\,\colon\, p\in\mathcal P^{2N+1}\right\}.
  \label{eq:radspace_onepole}
  	\end{align}
The space $\scalspextrad(\iepar_0,\iepar_0,N)=\vecspextrad(\iepar_0,\iepar_0,N)$ corresponds to the classical Hardy infinite element space introduced in \cite{hsm}. At the same time, the space   $\scalspextrad(\iepar_0,\iepar_1,N)=\vecspextrad(\iepar_0,\iepar_1,N)$ relates to the two-pole Hardy space, as introduced in \cite{HallaHohageNannenSchoeberl:16} (specifically to the two-scale version from \cite{HallaNannen:18}). 

The main difficulties in the practical use of the above spaces are the choice of a 'good' basis (it should be well-conditioned, yield sparse discretization matrices, and allow for a 'sparse' coupling between the interior and the exterior, cf. \cite[p.70]{WessDiss}) and further evaluation of the matrix elements. 
 
For the classical method a convenient set of basis functions was constructed in \cite{hsm} in the spacial Laplace domain. Afterwards in \cite{NW22} the basis was translated  into the physical space resulting in a subset of generalized Laguerre functions. The construction of such basis functions for the two-scale method is more subtle, and was done in the spacial Laplace domain in \cite{HallaHohageNannenSchoeberl:16}.
The numerical implementation of the above method relies on the knowledge of the mass and stiffness matrices ($\langle \tilde v_n,\tilde v_m\rangle_{L^2(\setRp)}$, $\langle r\tilde v_n,\tilde v_m\rangle_{L^2(\setRp)}$, $\langle \partial_r \tilde v_n,\partial_r \tilde v_m\rangle_{L^2(\setRp)}$, etc.), which can be computed semi-analytically. The corresponding expressions and computational procedures are stated for the convenience of the reader in Appendix \ref{app:implementation}.

%

\subsubsection*{Numerical experiments}
\label{sec:ie_experiments}

\begin{figure}[!b]
  \centering
  \begin{subfigure}[t]{0.45\textwidth}
  \includegraphics[width=\textwidth]{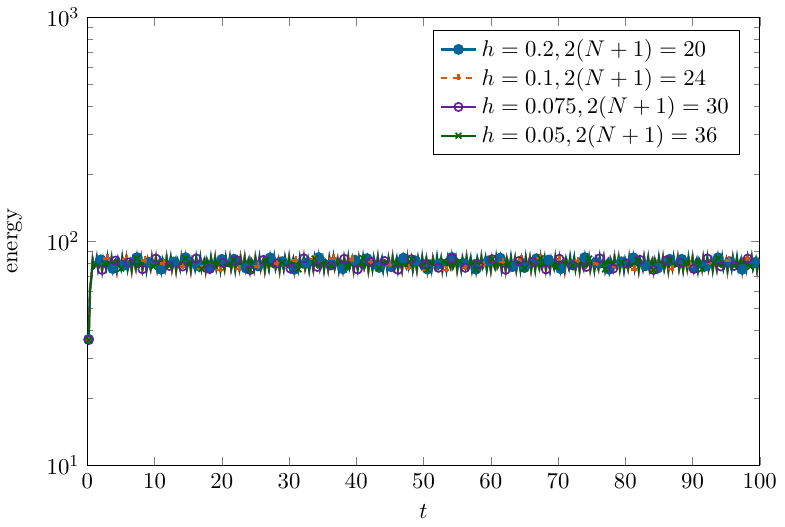}
  \caption{Long-time stability of infinite element discretizations for different mesh sizes and number of infinite elements. The energy curves are hard to distinguish, since they lie exactly on top of each other.}
  \label{fig:ie_gamma_stability}
  \end{subfigure}\hfill
  \begin{subfigure}[t]{0.45\textwidth}
  \includegraphics[width=\textwidth]{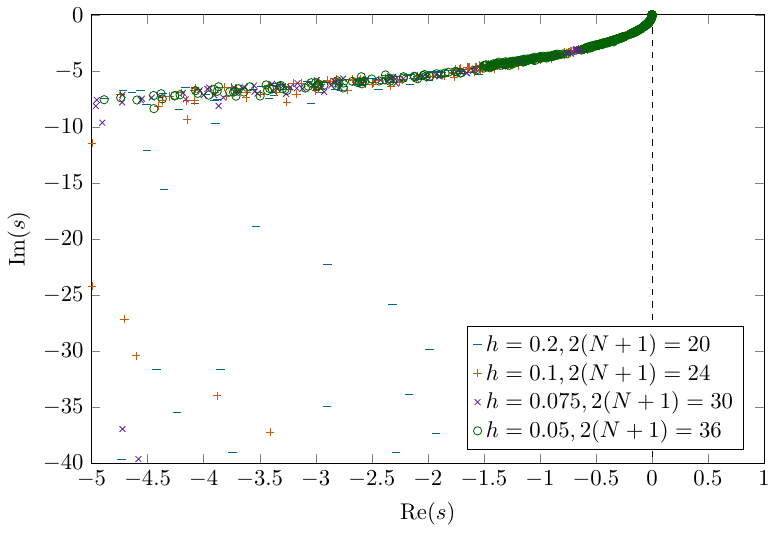}
  \caption{Spectra of the time harmonic problems corresponding to Figure \ref{fig:ie_gamma_stability}.}
  \label{fig:ie_gamma_stability_res}
  \end{subfigure}
  \caption{Stability of the frequency shifted infinite element discretizations.}
  \label{fig:ie_gamma_stability_all}
\end{figure}

\paragraph{Stability}
To illustrate numerically the stability of the system we discretize the same problem as in Section \ref{sec:first_experiments}, this time adding the $\gamma$-shift and using infinite elements (two-pole method with $\eta_0=1$ and $\eta_1=20$). 
In the interior  we use the elements with the order $k=4$ and vary mesh sizes. The number of infinite elements in the exterior is then varied accordingly (as defined by $\mathcal{V}_{\bdpml}$, $\vecsp_{\bdpml}$). 
%
{Because the statement of Theorem \ref{thm:stability_shifted_pml} is not quantitative (i.e. there is no explicit expression of $\nu_0$ given there), we choose $\nu$ satisfying \eqref{eq:stab_cond}.}
In particular, we take $\sigma_c=20$ and $\gamma=10$, which, together with  $\beta=\frac{\evmax-\evmin}{\evmax+\evmin}=4/5$,  leads to (cf. \eqref{eq:stab_cond})
\begin{align*}
  \nu=\frac{\sigma_c}{\gamma}=2<\nu_*=3.
\end{align*}
The dependence of the energy of the solution on time is shown in Figure \ref{fig:ie_gamma_stability}. We observe a long-time stability both for coarse and fine discretizations.
To underline the stability of the discrete system, Figure \ref{fig:ie_gamma_stability_res} shows the discrete spectra of the corresponding time-harmonic problems. Contrary to Figure \ref{fig:coarse_stability_res} even for finer discretizations the discrete resonances do not enter the positive complex half plane as indicated by the analysis in the previous sections.
\paragraph{Convergence}
\begin{figure}[t]
  \centering
  \begin{subfigure}{0.49\textwidth}
    \includegraphics[width=0.8\textwidth]{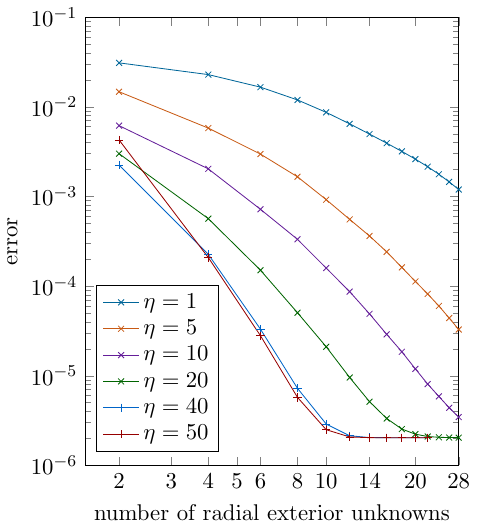}
    \caption{Logarithmic plot of the errors}
  \label{fig:convergence_ie_loglog}
  \end{subfigure}
  \begin{subfigure}{0.49\textwidth}
  \includegraphics[width=0.8\textwidth]{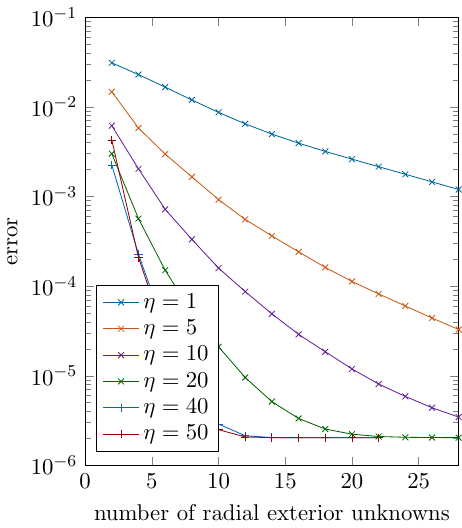}
    \caption{Semi-logarithmic plot of the errors}
  \label{fig:convergence_ie_semilog}
  \end{subfigure}
  \caption{Convergence of two-scale complex scaled infinite elements with respect to the number of the infinite elements in the basis $\mathcal{V}_{\operatorname{rad}}$.}
  \label{fig:convergence_ie}
\end{figure}
The theory from \cite{Halla:16,NW22} predicts super-algebraic convergence of the two-scale Hardy-space method w.r.t. the number of the infinite elements for time-harmonic problems (for a fixed frequency). The goal of this section is to verify numerically whether the same convergence rate can be achieved in the time-domain regime, in our setting. 
We choose $\domint=B_1\setminus\overline{B_{0.5}}:=\{\bx\in\setR^2:\|\bx\|<1\}$ and an anisotropy with $\bA=\diag(1/8,1)$. We use a piecewise constant scaling function $\sigma$ (i.e., satisfying Assumption \ref{assump:piecewise_const}), with $\sigma_c=20$, $\gamma=10$. As a source we choose
\begin{align*}
  f(t,\bx)=2400\sin(10t)\left(\exp\left(-160\|\bx-(0.7,0)^\top\|^2\right)+\exp\left(-160\|\bx+(0.7,0)^\top\|^2\right)\right),
\end{align*}
As before, we merely simulate one quarter of the domain.
The interior discretization is done by the finite elements of order $k=6$ and a mesh-size of $h=0.075$. For the infinite elements,  
in all our experiments, we set $\eta_0=1$, and $\eta_1=\eta$, and vary $\eta$. We compare the numerical solution to the reference solution. The latter is computed by using homogeneous boundary conditions on a domain large enough such that the wave is not yet reflected back to the interior and finite elements of order $k=8$.

Figure \ref{fig:convergence_ie} displays convergence of the $L^2$-errors (in space and time in $\domint$) of the primal variable $u^{\sigma}(t,\bx)$ for end-time $T=0.5$ for various choices of $\eta$.  The error from the interior- and time-discretization is roughly $2\cdot 10^{-6}$, thus the observed, converging dominating error is exclusively the error of the exterior discretization.

Like in the time-harmonic regime, we observe the super-algebraic convergence of the error with respect to the number of the infinite elements. 
Moreover, as we see, the two-scale Hardy finite element space drastically outperforms the one-pole method, though in both cases the super-algebraic convergence can be observed. This is most likely due to a poor quality of approximation of the solutions vanishing outside of a ball of radius $R=R(T)$ where $T$ is the final simulation time, due to the finite speed of wave propagation. Indeed, such solutions can be approximated better by $\mathrm{e}^{-\eta x}$, with sufficiently large $\eta$. Choosing optimal parameters for the performance of such infinite element methods, however, remains an open question. 	


\FloatBarrier
\subsection{Method 2: truncation-free PMLs}
\label{sec:numerics_infpml}
Let us now present an alternative method. We use the new change of variables, cf. \eqref{eq:radpml}:
\begin{align}
	\label{eq:pml_sigma_gamma}
\tilde{r}_{\sigma,\gamma}:=
\left\{
\begin{array}{ll}
f_{\lpml}(r)+\frac{1}{s+\gamma}\int_{\radpml}^{f_{\lpml}(r)}\sigma(r')dr', & r>\radpml, \\
r, & r\leq \radpml.
\end{array}
\right.
\end{align}
with 
\begin{align*}
  f_{\lpml}(r)=\frac{\radpml \lpml}{\radpml+\lpml-r}.
\end{align*}
The change of variables defined by $f_{\lpml}$ is a bijective mapping from $[\radpml,\radpml+\lpml)$ to $[\radpml,\infty)$ (and thus maps $\domext$ into a bounded layer $\dompml$ of width $L$). It therefore yields a truncation-free version of the PML. 
Note that this approach follows the ideas for the frequency domain of \cite{HugoninLalanne:05,yang_wang_gao} and \cite[Sect.\ 4.5.1]{Halla:19Diss}, and is different to the one used in \cite{BermudezHervellaNPrietoRodriguez:08} where merely the imaginary part of the scaling is unbounded.
To avoid singular integrals due to the change of variables defined above we use an approach similar to \cite{yang_wang_gao}. We apply the change of variables to the weak form \eqref{eq:rad_time_weak}  while  the test functions $v^\dagger$ are replaced by $(\radpml+\lpml-r)^3v^\dagger$ (differing slightly from \cite[3.3.2]{yang_wang_gao} where also the trial function is scaled).

\subsubsection*{Numerical experiments}
\label{sec:infpml_experiments}
\begin{figure}[!b]
  \centering
  \begin{subfigure}[t]{0.45\textwidth}
  \includegraphics[width=\textwidth]{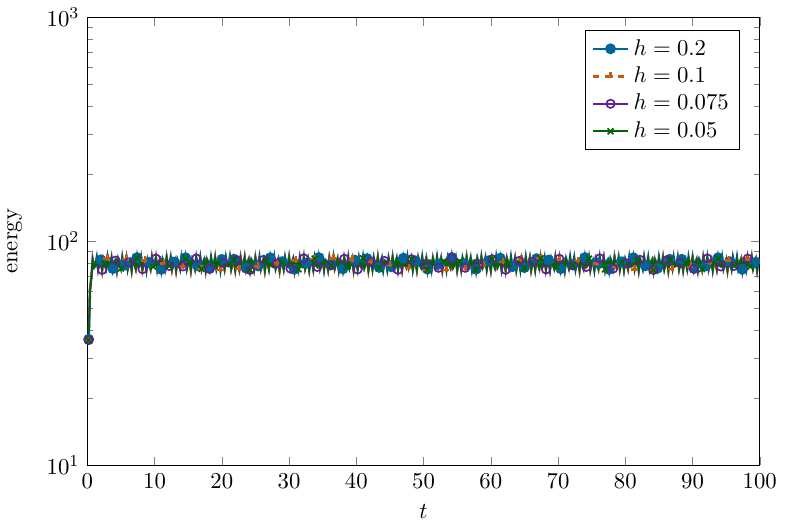}
  \caption{Long-time stability of truncation-free PML discretizations for different mesh sizes. The energy curves are hard to distinguish, since they lie exactly on top of each other.}
  \label{fig:infpml_gamma_stability}
  \end{subfigure}\hfill
  \begin{subfigure}[t]{0.45\textwidth}
  \includegraphics[width=\textwidth]{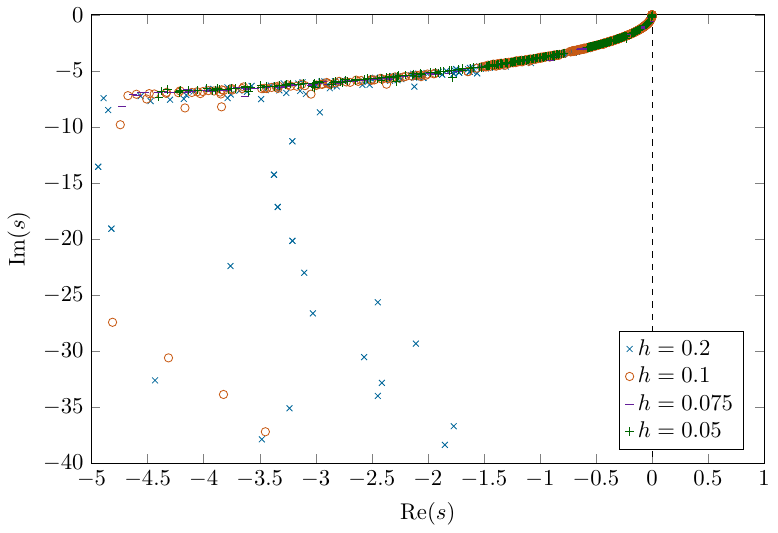}
  \caption{Spectra of the time harmonic problems corresponding to Figure \ref{fig:infpml_gamma_stability}.}
  \label{fig:infpml_gamma_stability_res}
  \end{subfigure}
  \caption{Stability of the frequency shifted infinite PMLs.}
  \label{fig:infpml_gamma_stability_all}
\end{figure}

\begin{figure}[h]
  \centering
  \includegraphics[width=0.5\textwidth]{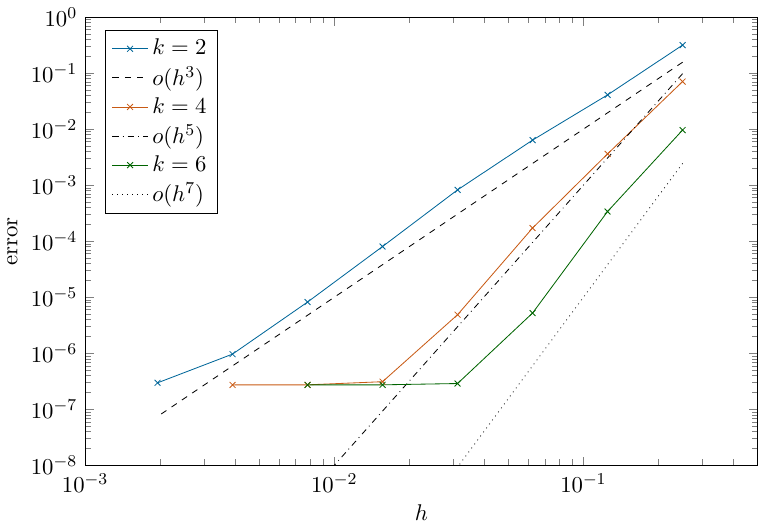}
  \caption{Convergence of the truncation-free PMLs.}
  \label{fig:infpml_convergence}
\end{figure}

In all the experiments we  set $L=1$. 
\paragraph{Stability}
We conduct the same experiments as in Sections \ref{sec:first_experiments}, \ref{sec:ie_experiments}, but using the approach described above.
Similarly to the infinite element approach, we observe long-time stability in the time domain experiments in Figure \ref{fig:infpml_gamma_stability}. This is again confirmed by the fact that, at least for the discretizations considered, the respective resonance problems have no spectrum in $\setCp$, as illustrated in Figure \ref{fig:infpml_gamma_stability_res}.
%
%
\paragraph{Convergence}
Contrary to the case of the infinite elements it is not straightforward to separate the interior from the exterior discretization error, due to the fact that the finite element meshes in the interior and exterior domains are not independent. Thus we use a uniform mesh size in the whole domain to study the convergence of the $L^2$-error (in time and space) of the primal variable $u^\sigma(t,x)$ in the interior domain (cf. Figure \ref{fig:infpml_convergence}). The remaining parameters are identical to the ones in the convergence experiments in Section \ref{sec:ie_experiments}.
For truncation-free PMLs we expect to observe two different kinds of errors (apart from the discretization error of the interior domain):
\begin{itemize}
  \item{the discretization error of the absorbing layer which is expected to be of order $h^{k+1}$ (in the $L^2$-norm), and}
  \item{an additional error due to the fact that the outmost elements with homogeneous Dirichlet boundary conditions in the absorbing layer are mapped to an infinite domain, which is expected to cause some sort of artificial reflections.}
\end{itemize}
In the experiments in Figure \ref{fig:infpml_convergence} we use an end time $T=0.5$, and thus we observe the first kind of error and the convergence of order $h^{k+1}$. The plateau in Figure \ref{fig:infpml_convergence} occurring at errors of magnitude $~2\cdot 10^{-7}$  is the time-discretization error. In the section that follows we perform experiments on longer times.

\subsection{Comparison of the two methods}
\label{sec:long_time}
\begin{figure}
  \centering
  \begin{subfigure}[t]{0.45\textwidth}
  \includegraphics[width=\textwidth]{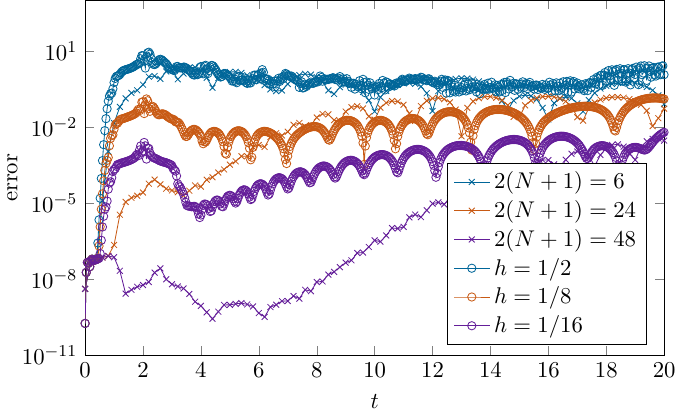}
  \caption{The error of infinite PMLs and infinite elements for various discretizations with respect to time.}
  \label{fig:convergence_bessel_time}
  \end{subfigure}\hfill
  \begin{subfigure}[t]{0.45\textwidth}
  \includegraphics[width=\textwidth]{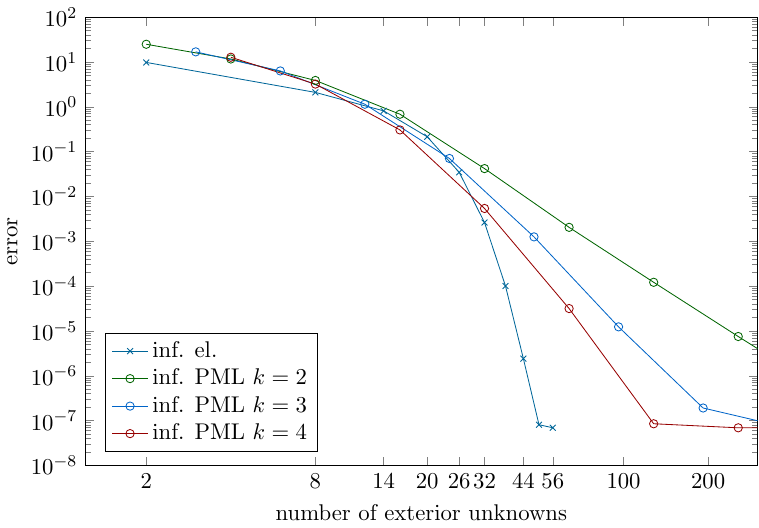}
  \caption{Comparison of the error of the two methods with respect to exterior degrees of freedom.}
  \label{fig:convergence_comparison_longtime_loglog}
  \end{subfigure}
  \caption{Comparison of the error of truncation-free PMLs and infinite elements.}
  \label{fig:bessel_plots}
\end{figure}

The goal of this short section is to compare the methods presented in Sections \ref{sec:numerics_infpml} and \ref{sec:numerics_ie} numerically. 
Since for larger end times $T$ we were not able to compute the reference solution with a very high accuracy, we will compare the performance of these methods for a model isotropic problem.  
In this setting, we consider the problem \eqref{eq:main_problem} with $\bA=\operatorname{Id}$, vanishing source $f=0$, and the initial conditions $\left.u(t, r,\phi)\right|_{t=0} = 120\exp(-50r)$, $\left. \partial_t u(t, r,\phi)\right|_{t=0}=0$, so that the solution does not depend on the variable $\phi$. This allows to reduce the problem to a single spatial dimension.  We next solve this new problem by applying two methods described above. 
We use the following parameters for the methods: $\sigma_c=20$, $\gamma=10$, $\iepar_0=1, \, \iepar_1=20$ and $L=1$.
We choose a very fine interior discretization, so that the interior error is dominated by the error from the exterior discretization (and until the wave reaches the absorbing layer it is dominated by the time stepping error; the time step is chosen as $\tau=10^{-5}$ which  gives an error of $\approx 10^{-7}$).
In Figure \ref{fig:convergence_bessel_time}, we study the dependence of the $L^2$-spatial errors on time for both methods.
The number of infinite elements and the exterior discretization of the truncation-free PMLs with order $k=3$ is chosen so that the number of unknowns discretizing the exterior is equal in both methods.

We observe that the error starts growing when the wave hits the absorbing layer at $T=1$. After an initial settling phase the errors grow. Nonetheless, on smaller time intervals the error of the infinite elements is by far smaller. At the end of the experiments the errors of the two methods are comparable.
Figure \ref{fig:convergence_bessel_time} shows the error of the two methods for various discretizations with respect to time.
Figure \ref{fig:convergence_comparison_longtime_loglog} shows comparison of the two methods for the separated equation at the end time $T=10$ with respect to the number of degrees of freedom in the exterior domain. For the chosen examples the infinite element method is more efficient than the 
truncation-free PML.
\FloatBarrier 
\section*{Used notation}
\noindent
\begin{tabularx}{\textwidth}{p{2cm}p{8cm}l}
\textbf{symbol}&\textbf{description}&\textbf{definition}\\
\hline
\multicolumn{3}{c}{\textit{sets of numbers}} \\
\hline
$\setRpz$&non-negative real numbers& $\{t\in\setR:t\geq 0\}$\\
$\setRp$&positive real numbers &$\{t\in\setR:t> 0\}$\\
$\setCp$&complex numbers with positive real part&$\{s\in\setC:\Re s>0\}$\\
$\setCp_\alpha$&complex numbers with real part bounded from below&$\{s\in\setC:\Re s>\alpha\}$\\
$B_R$&open ball of radius $R$&$\{\bx\in\setR^2\text{ or }\setC:\|\bx\|< R\}$\\
\end{tabularx}\\
\begin{tabularx}{\textwidth}{p{2cm}p{8cm}l}
\hline
\multicolumn{3}{c}{\textit{problem parameters}} \\
\hline
$\bA$&material anisotropy&mostly $\bA = \diag(\evmax^{-1},\evmin^{-1})$\\
$\bB$&&$\bA^{-1}$\\
$\bA^\phi$&material anisotropy in polar coordinates&\\
$\evmin,\evmax$&minimal and maximal eigenvalues of $\bB$&\\
$\mu^*$&stability radius&{$\frac{1+\mu}{2\sqrt{\mu}}=\frac{\evmax+\evmin}{2\sqrt{\evmax\evmin}}$,
$\mu=\frac{\evmin}{\evmax}$}\\
\end{tabularx}\\
\begin{tabularx}{\textwidth}{p{2cm}p{8cm}l}
\hline
\multicolumn{3}{c}{\textit{domains}} \\
\hline
$\domint$&interior domain& for radial PMLs $B_{\radpml}$\\
$\dompml^\infty$&(untruncated) exterior domain& $\setR^2\setminus\overline{\domint}$\\
$\dompml$&truncated exterior domain& for radial PMLs $B_{\radpml+L}\setminus\overline{\domint}$\\
\end{tabularx}\\
\begin{tabularx}{\textwidth}{p{2cm}p{8cm}l}
\hline
\multicolumn{3}{c}{\textit{PML quantities}} \\
\hline
$\bA^\phi_\sigma$&PML material in polar coordinates&\eqref{eq:rad_lpl_v0}\\
$\bA_\sigma$&PML material&\eqref{eq:basigma_main}\\
$L$&PML thickness&\\
$\radpml$&starting radius of the complex layer&\\
$\sigma$&damping function\\
$\sigma_c$&damping constant\\
$\gamma$&complex frequency shift parameter&\\
$\nu$&$\sigma_c/\gamma$&\eqref{eq:nudef}\\
$\nu_*$ & $2\sqrt{\beta^{-2}-1}(\sqrt{\beta^{-2}-1}+\beta^{-1})$,
$\beta=\frac{\evmax-\evmin}{\evmax+\evmin}$&\eqref{eq:stab_cond}\\
$r_\sigma(s,r)$&complex scaled radius&\eqref{eq:radpml}\\
$\dtpml,\dpml,\tilde\sigma$&PML {auxiliary functions}&\eqref{eq:pmldefns}\\
$\vecxs(s,\bx)$&complex scaled variable&\eqref{eq:pmldefns}\\
$u^\sigma$&complex scaled solution&\\
$\Js$&Jacobian of the scaling&$D_{\bx} \vecxs$\\
$G_\sigma$&fundamental solution of the PML system&Proposition \ref{prop:fs}\\
$h_\sigma$&$(\vecxs-\by)^\top\bB(\vecxs-\by)$&\eqref{eq:defh}\\
$\gamma_{11},\gamma_{12},\gamma_{22}$&$\hat\bc^\top\bB\hat\bc$, 
$\hat\bx^\top\bB\hat\bc$, $\hat\bx^\top\bB\hat\bx$ &\eqref{eq:notation2}\\
$\bc$&$\radpml\hat\vecx-\by$&\eqref{eq:notation1}\\
$\xi$ &$\|\bx\|-\radpml$&\eqref{eq:notation1}\\
\end{tabularx}\\
\begin{tabularx}{\textwidth}{p{2cm}p{8cm}l}
\hline
\multicolumn{3}{c}{\textit{various definitions}} \\
\hline
$\Pi_\parallel(\vecx),\Pi_\perp(\vecx)$&projections onto the radial/tangential direction&\\
$\hat u$&Laplace transform in time&\\
$\hat\bx$&unit vector $\bx/\|\bx\|$&
\end{tabularx}\\
\FloatBarrier
\bibliography{bibliography}
\appendix
\section{Auxiliary lemmas quantifying the behaviour of $\hpml$}
\label{appendix:h_behave}
We will assume everywhere in this section that $\sigma$ satisfies Assumption \ref{assump:piecewise_const}, i.e., is piecewise-constant and equal to $\sigma_c>0$ inside the perfectly matched layer. We additionally use the notation $\|\bx\|_{\bB}:=\sqrt{\bx^\top\bB\bx}$, for an arbitrary vector $\bx\in \setR^2$.
Recall the definition \eqref{eq:defh}, which we extend to $\by\in \domext$:
\begin{align*}
	\hpml(s; \bx,\by)=(\bx_{\sigma}-\by_{\sigma})^\top \bB(\bx_{\sigma}-\by_{\sigma}), \quad \bx,\by\in \setR^2.
\end{align*}
In what follows, we will need more explicit expressions and bounds on $\hpml$ in different regions of $\setR^2\times \setR^2$. We will first state and prove these bounds for $\bx,\by\in \overline{\domext}\times \overline{\domext}$, next for $\bx\in \overline{\domext}$ and $\by\in \domext$, and, finally, gather all the bounds in Section \ref{sec:summary_of_the_bounds}
\subsection{Bounds on $\hpml(s; \bx,\by)$ for $(\bx,\by)\in {\overline{\domext}}\times {\overline{\domext}}$}
\label{section:bounds_classical_dompml_dompml}
\subsubsection{Bounds for classical PMLs}
As a preparation for the proofs of the following results we begin by rewriting
\begin{align}
	\label{eq:ppml}
	\bx_{\sigma}-\by_{\sigma}=(\bx-\by)(1+s^{-1}\sigma_c)-s^{-1}\sigma_c(\radpml\hat{\bx}-\radpml\hat{\by})=\mathbf{p}(1+s^{-1}\sigma_c)-s^{-1}\sigma_c \ppml,
\end{align} 
with $\bp=\bx-\by$ and $\ppml=(\radpml\hat{\bx}-\radpml\hat{\by})$. We then have that 
\begin{align}
	\label{eq:exprhpml}
	\hpml(s; \bx,\by)=\left(1+\frac{\sigma_c}{s}\right)^2\bp^{\top}\bB\bp-2\frac{\sigma_c}{s}\left(1+\frac{\sigma_c}{s}\right)\bp^{\top} \bB\ppml+\frac{\sigma_c^2}{s^2} \ppml^{\top}\bB\ppml.
\end{align}
We start with the following bound. 
\begin{lem}
	\label{lem:A1_part1}
	Let $\sigma_c>0$. Then there exist $\gamma_{\bB}, \, c_{\bB}>0$, which depend on $\bB$, s.t.\ for all  $s\in\setCp_{\alpha}$, with  $\alpha=\gamma_{\bB}\sigma_c$, 
	\begin{align}
		\label{eq:first_bound_part1}
		&\Re \hpml(s; \bx,\by)\geq c_{\bB}\|\bx-\by\|^2,\qquad \text{ for all }\bx,\by\in \overline{\domext}. 
	\end{align}
\end{lem}
\begin{proof}[Proof of Lemma \ref{lem:A1_part1}]
We make use of \eqref{eq:exprhpml}, \eqref{eq:ppml}. 
	For all $s\in\setCp$ we have that $\Re(1+s^{-1}\sigma_c)^2\geq \Re(1+s^{-2}\sigma_c^2)$, and therefore, in particular, for $s\in \setCp_{\alpha}$, it holds that 
	\begin{align*}
		\Re \hpml(s; \bx,\by)\geq \left(1-\frac{\sigma_c^2}{\alpha^2}\right)\|\bp\|^2_{\bB}-\frac{2\sigma_c}{\alpha}\left(1+\frac{\sigma_c}{\alpha}\right)\left|\bp^{\top}\bB\ppml\right|-\frac{\sigma^2_c}{\alpha^2}\|\ppml\|^2_{\bB}.
	\end{align*}
	We have that $\|\bp\|\geq \|\ppml\|$.\footnote{ Indeed, denoting by $b:=\hat{\bx}^\top\hat{\by}$, we remark that the minimum of the linear in $b\in[-1,1]$ function
		\begin{align*}
			\|\bp\|^2-\|\ppml\|^2=\|\bx\|^2+\|\by\|^2-2\|\bx\|\|\by\|b-2\radpml^2(1-b)
		\end{align*}
		is achieved at $b=1$ ($\|\bx\|,\|\by\|>\radpml$), and equal to $|\|\bx\|-\|\by\||^2\geq 0$.
	}
	Therefore, 
	\begin{align}
		\label{eq:important_bound1}
                \|\bp\|_{\bB}\geq\sqrt{\evmin}\|\bp\|\geq\sqrt{\evmin}\|\ppml\|\geq \sqrt{\evmin/\evmax}\|\ppml\|_{\bB}.
	\end{align}
	We conclude that for all $\alpha>\gamma_{\bB}\sigma_c$, with $\gamma_{\bB}>0$ sufficiently large, all the sign-indefinite and negative terms are controlled by $\|\bp\|_{\bB}^2$. Hence the bound in the statement of the lemma.
\end{proof}
An upper bound on $\hpml$, stated below, is then obtained similarly, and thus we leave its proof to the reader. 
\begin{lem}
	\label{lem:A2_part1}
	Let $\sigma_c>0$. Then there exists a constant $C_{\bB}>0$, which depends on $\bB$, s.t. for all  $s\in\setCp$, 
	\begin{align*}
		|\hpml(s; \bx,\by)|\leq C_{\bB}\left(1+\frac{\sigma_c}{|s|}\right)^2\|\bx-\by\|^2, \quad \text{for all }\bx,\by\in \overline{\domext}.
	\end{align*}
\end{lem}
We also have another lower bound on $\hpml$ which will be of use later. 
\begin{lem}
	\label{lem:A3_part1}
	Let $\sigma_c>0$. Then, there exist $\tilde{\gamma}_{\bB}, \, \tilde{c}_{\bB}>0$, which depend on $\bB$,  s.t.\ for all  $s\in\setCp_{\tilde{\alpha}}$, with  $\tilde{\alpha}=\tilde{\gamma}_{\bB}\sigma_c$,
	\begin{align}
		\label{eq:res4_standard_pml_part1}
		\Re(s \sqrt{\hpml(s;\bx, \by)})\geq \tilde{c}_{\bB}\Re s\|\bx-\by\|,\, \text{ whenever }\bx,\by\in \overline{\domext}.
	\end{align}
\end{lem} 
\begin{proof} 
We use the identity $$\Re (s\sqrt{\hpml})=\Re s\Re \sqrt{\hpml}-\Im s\Im \sqrt{\hpml}.$$

Recall the explicit expression for the complex root:
\begin{align}
	\sqrt{z}=\sqrt{\frac{|z|+\Re z}{2}}+i\operatorname{sign}\Im z\sqrt{\frac{|z|-\Re z}{2}}.
	\label{eq:sqrt}
\end{align}
A lower bound for $\Re \sqrt{\hpml}>\sqrt{\Re\hpml}$ can be obtained from Lemma \ref{lem:A1_part1}. This yields 
\begin{align}
	\label{eq:bound0}
	\Re (s\sqrt{\hpml(s; \bx,\by)})=\Re s\Re \sqrt{\hpml(s; \bx,\by)}-\Im s\Im \sqrt{\hpml(s; \bx,\by)}\geq c_{\bB}^{1/2}\Re s \|\bx-\by\|-|\Im s\Im \sqrt{\hpml(s; \bx,\by)}|.
\end{align}
To bound $\Im \sqrt{\hpml}$, we use that, for $z\in \mathbb{C}\setminus \setR^{-}_0$, cf.\ \eqref{eq:sqrt},
\begin{align*}
	\left|\Im \sqrt{z}\right|=\sqrt{\frac{|z|-\Re z}{2}}.
\end{align*}
We next use the inequality $\sqrt{a^2+b^2}-a\leq \frac{b^2}{a}$, valid for $a>0$ ({remark that this is more optimal than $\sqrt{a^2+b^2}-a\leq |b|$ for small $|b|/a$}), with $a=\Re{\hpml}$ and $b=\Im\hpml$. This inequality implies that
$
\left|\Im \sqrt{\hpml}\right|\leq \frac{|\Im \hpml|}{\sqrt{\Re \hpml}}.
$
With Lemma \ref{lem:A1_part1}, for $s\in \setCp_{\alpha}$, it holds that
\begin{align}
	\label{eq:imhpml}
	\left|\Im \sqrt{\hpml}\right|\leq \frac{|\Im \hpml|}{c_{\bB}^{1/2}\|\bp\|}.	
\end{align}
Finally, it remains to bound $|\Im \hpml|$.  From the expression \eqref{eq:exprhpml}, it follows that 
\begin{align}
\nonumber
	\Im \hpml(s; \bx,\by)&=\left(2\sigma_c\Im \frac{1}{s}+\sigma_c^2\Im\frac{1}{s^2}\right)\|\bp\|^2_{\bB}-\left(2\sigma_c\Im \frac{1}{s}+2\sigma_c^2\Im \frac{1}{s^2}\right)\bp^{\top}\bB\ppml+\sigma_c^2\Im\frac{1}{s^2}\|\ppml\|^2_{\bB}\\
		\label{eq:expr}
	&=-2\sigma_c \frac{\Im s}{|s|^2}(\|\bp\|_{\bB}^2-\bp^{\top}\bB\ppml)-2\sigma_c^2  \frac{\Re s\Im s}{|s|^4}\|\bp-\ppml\|^2_{\bB}.
\end{align}
  With the inequality  \eqref{eq:important_bound1} $\|\bp\|_{\bB}\geq \sqrt{\evmin/\evmax}\|\ppml\|_{\bB}$, and using $\frac{|\Im s|}{|s|}<1, \frac{\Re s}{|s|}<1$, the above yields, with some constant $\tilde{C}_{\bB}>0$,
\begin{align}
	\label{eq:ims}
	|\Im \hpml(s; \bx,\by)|\leq \tilde{C}_{\bB} \frac{\sigma_c}{|s|}\left(1+\frac{\sigma_c}{|s|}\right)\|\bp\|^2.
\end{align}
  Inserting the above inequality into \eqref{eq:imhpml}, and next inserting \eqref{eq:imhpml} into \eqref{eq:bound0}, we obtain,
\begin{align*}
	\Re (s\sqrt{\hpml(s; \bx,\by)})&\geq c_{\bB}\Re s\|\bx-\by\|-\frac{\tilde{C}_{\bB}\sigma_c |\Im s|}{c_{\bB}^{1/2}|s|}\left(1+\frac{\sigma_c}{|s|}\right)\|\bx-\by\|\\
	&\geq  c_{\bB}\Re s\left(1-\frac{\tilde{C}_{\bB}}{c_{\bB}^{3/2}}\frac{\sigma_c}{\Re s}\left(1+\frac{\sigma_c}{\Re s}\right)\right)\|\bx-\by\|.
\end{align*}
Choosing $\Re s>\tilde{\gamma}_{\bB}\sigma_c$ with sufficiently large $\tilde{\gamma}_{\bB}$ depending on $\bB$ only yields the desired estimate in the statement of the lemma. 
\end{proof}
\subsubsection{Bounds on complex-shifted PMLs}
For the complex-shifted PMLs, the counterpart of \eqref{eq:exprhpml} reads
\begin{align}
	\label{eq:exprhpml2}
	\hpml(s+\gamma; \bx,\by)=\left(1+\frac{\sigma_c}{s+\gamma}\right)^2\bp^{\top}\bB\bp-2\frac{\sigma_c}{s+\gamma}\left(1+\frac{\sigma_c}{s+\gamma}\right)\bp^{\top} \bB\ppml+\frac{\sigma_c^2}{(s+\gamma)^2} \ppml^{\top}\bB\ppml.
\end{align}
We start with the following auxiliary result.
\begin{lem}
	\label{lem:aux_res_part1}
	There exist constants $\nu_0>0$ and $c_{\bB}>0$ depending on $\bB$, s.t.\ for all  $\sigma_c<\nu_0\gamma$ the inequality  
	\begin{align}
		\label{eq:rehpml}
		\Re \hpml(s+\gamma; \bx,\by)>c_{\bB}\|\bx-\by\|^2
	\end{align}
	holds true for all $s\in \setCp$ and $\bx,\by\in \overline{\domext}$.
\end{lem}
\begin{proof}
	By Lemma \ref{lem:A1_part1}, one has $\Re\hpml(s; \bx,\by)>c_{\bB}\|\bx-\by\|^2$, whenever $\Re s>\gamma_{\bB}\sigma_c$. In other words, $\Re\hpml(s+\gamma; \bx,\by)>c_{\bB}\|\bx-\by\|^2$ whenever $\Re s+\gamma>\gamma_{\bB}\sigma_c$. If $\gamma>\gamma_{\bB}\sigma_c$ (i.e. $\nu_0=\gamma_{\bB}^{-1}$), the estimate $\Re\hpml(s+\gamma; \bx,\by)>c_{\bB}\|\bx-\by\|^2$ holds for all $s\in \mathbb{C}^+$.
\end{proof}
With the above lemma, we can prove the following proposition. 
 \begin{prop}
	\label{prop:cs_pml_part1}
	There exist constants $\nu_0,\, C_{\pm},\, c_{\pm} >0$ depending on $\bB$ only, s.t.\ for all $\sigma_c,\gamma>0\colon\, \sigma_c<\nu_0\gamma$, for any $s\in \mathbb{C}^+$, any $\bx,\by\in \overline{\domext}$, the following bounds hold true:
	\begin{align}
		\label{eq:basic_lower_bound_appendix1}
		&|s\sqrt{\hpml(s+\gamma;\bx,\by)}|\geq C_{-}|s|\|\bx-\by\|,\\
		\label{eq:basic_upper_bound_appendix1}
		& |s\sqrt{\hpml(s+\gamma;\bx,\by)}|\leq C_+|s|\|\bx-\by\|,\\
		\label{eq:basic_lower_bound_appendix1_v2}
		&\Re(s\sqrt{\hpml(s+\gamma;\bx,\by)})\geq (c_+\Re s-c_{-}\sigma_c)\|\bx-\by\|.
	\end{align}
If, additionally, $\|\bx-\by\|>\rho_*=2\radpml\sqrt{\evmax/\evmin}$, then 
\begin{align}
		\label{eq:basic_lower_bound_appendix1_v3}
		\Re \left(s\sqrt{\hpml(s+\gamma;\bx,\by)}\right)\geq c_{+}\Re s\|\bx-\by\|.
\end{align}
\end{prop}
The value of $\rho_*$ in the above is not optimal, but we do not pursue the goal of optimizing it. 
\begin{proof}
The bound \eqref{eq:basic_lower_bound_appendix1} follows from Lemma \ref{lem:aux_res_part1}, combined with $|s\sqrt{h_{\sigma}}|=|s|\sqrt{|h_{\sigma}|}\geq |s|\sqrt{\Re h_{\sigma}}$. The upper bound \eqref{eq:basic_upper_bound_appendix1} can be obtained by straightforward computations from \eqref{eq:exprhpml2}, thus we leave the proof to the reader. 
Let us now prove the bound \eqref{eq:basic_lower_bound_appendix1_v2}. 
We have that 
\begin{align}
	\label{eq:hpml_lb}
	\Re (s\sqrt{\hpml(s+\gamma;\bx,\by)})=\Re s\Re \sqrt{\hpml(s+\gamma; \bx,\by)}-\Im s\Im\sqrt{\hpml(s+\gamma; \bx, \by)}. 
\end{align}
By Lemma \ref{lem:aux_res_part1},  and using the bound $\Re\sqrt{h_{\sigma}}\geq \sqrt{\Re h_{\sigma}}$, we have that 
\begin{align}
	\label{eq:lower_bound_reocc}
	\operatorname{Re}\sqrt{\hpml(s+\gamma; \bx,\by)}>c_{\bB}^{1/2}\|\bx-\by\|.
\end{align} Therefore, it remains to bound $\left|\Im s\Im\sqrt{\hpml(s+\gamma; \bx, \by)}\right|$. We proceed like in the proof of Lemma \ref{lem:A3}; in particular, we have, with \eqref{eq:ims} that 
\begin{align*}
	|\Im \hpml(s+\gamma; \bx,\by)|\leq \tilde{C}_{\bB}\frac{\sigma_c}{|s+\gamma|}\left(1+\frac{\sigma_c}{|s+\gamma|}\right)\|\bx-\by\|.
\end{align*} 
Inserting the two above bounds into \eqref{eq:hpml_lb} yields 
\begin{align*}
		\Re (s\sqrt{\hpml(s+\gamma;\bx,\by)})\geq c_{\bB}^{1/2}\Re s\|\bx-\by\|-\tilde{C}_{\bB}\frac{\sigma_c|\Im s|}{|s+\gamma|}\left(1+\frac{\sigma_c}{|s+\gamma|}\right)\|\bx-\by\|.
\end{align*}
With $|\Im s|<|s+\gamma|$ and using that $\gamma>\nu_0^{-1}\sigma_c$, we obtain the desired estimate \eqref{eq:basic_lower_bound_appendix1_v2}. 

Let us finally prove \eqref{eq:basic_lower_bound_appendix1_v3}. We use again \eqref{eq:hpml_lb} and \eqref{eq:lower_bound_reocc}. The goal is to prove that for $\|\bx-\by\|>\rho_*$, $\Im s\Im \sqrt{\hpml(s+\gamma; \bx,\by)}<0$. 
For this we consider the following expression, cf. \eqref{eq:expr}:
\begin{align}
		\label{eq:hpml_sigma}
	\Im\hpml(s+\gamma; \bx,\by)
	&=-\frac{2\sigma_c\Im s}{|s+\gamma|^2}\underbrace{(\|\bp\|_{\bB}^2-\bp^{\top}\bB\ppml)}_{S}-\frac{2\sigma_c^2\Im s(\Re s+\gamma)}{|s+\gamma|^4}\|\bp-\ppml\|^2_{\bB}.
\end{align}
Next consider the term 
$$S=\|\bp\|_{\bB}^2-\bp^{\top}\bB\ppml\geq \|\bp\|_{\bB}^2-\|\bp\|_{\bB}\|\ppml\|_{\bB}.$$ 
We have  $\|\ppml\|=\radpml\|\hat{\bx}-\hat{\by}\|\leq 2\radpml$, and thus $\|\ppml\|_{\bB}<2\evmax^{1/2}\radpml$. This, combined with $\|\bp\|_{\bB}\geq \evmin^{1/2}\|\bp\|$, shows that, for  $\|\bp\|>2\evmax^{1/2}/\evmin^{1/2}\radpml$, we have that $\|\bp\|_{\bB}>\|\ppml\|_{\bB}$, and the quantity $S$ is positive.
Thus, under the condition $\|\bx-\by\|>\rho_*$, $\Im s\Im h_{\sigma}(s+\gamma;\bx,\by)\leq 0$ and since $\Im \sqrt{\hpml}\Im\hpml\geq 0$ (cf. \eqref{eq:sqrt}), we have that $\Im s\Im\sqrt{\hpml(s+\gamma;
	\bx,\by)}\leq 0$. We then have that  \eqref{eq:hpml_lb}, with the use of \eqref{eq:lower_bound_reocc} rewrites as \eqref{eq:basic_lower_bound_appendix1_v3}:
\begin{align*}
	\Re(s\sqrt{\hpml(s+\gamma; \bx,\by)})\geq c_{\bB}^{1/2}\Re s\|\bx-\by\|.
\end{align*} 
\end{proof}
\subsection{Bounds on $\hpml$ for $(\bx,\by)\in {\overline{\domext}}\times \domint$}
For $\bx\in \overline{\domext}$ and $\by \in \domint$, the expression $\hpml(s; \bx,\by)$ has a particularly simple form. Indeed, 
\begin{align}
	\label{eq:ppml_half}
	\bx_{\sigma}-\by=(\bx-\by)+\frac{\sigma_c}{s}(\bx-\radpml\hat{\bx})=(\bx-\by)+\frac{\sigma_c}{s}(\|\bx\|-\radpml)\hat{\bx},
\end{align}
so that, with $\bp=\bx-\by$ and $\xi=\|\bx\|-\radpml$, we have 
\begin{align}
	\label{eq:exprhpml_half}
	\hpml(s; \bx,\by)=\|\bp\|^2_{\bB}+2\frac{\sigma_c \xi}{s}\hat{\bx}^{\top}\bB\bp+\frac{\sigma_c^2}{s^2}\xi^2\|\hat{\bx}\|^2_{\bB}.
\end{align}

\subsubsection{Bounds for classical PMLs}
Let us list the results which are the counter-parts of Lemmas of Section \ref{section:bounds_classical_dompml_dompml}. 
Remark that the key inequality in the proof of lemmas \ref{lem:A1_part1}, \ref{lem:A2_part1} is \eqref{eq:important_bound1}. Its counterpart is 
\begin{align}
	\label{eq:important_bound2}
	\|\bp\|_{\bB}\geq  \sqrt{\evmin/\evmax}\|\xi\hat{\bx}\|_{\bB}, 
\end{align}
which follows because $\|\bp\|\geq \xi$.
We start with counterparts of Lemmas \ref{lem:A1_part1}, \ref{lem:A2_part1}. 
\begin{lem}
		\label{lem:A1_part2}
Lemma \ref{lem:A1_part1} holds for $\bx\in \overline{\domext}$, $\by\in \domint$, with possibly different constants $\gamma_{\bB}$ and $c_{\bB}$. 
\end{lem} 
\begin{lem}
	\label{lem:A2_part2}
Lemma \ref{lem:A2_part1} holds for $\bx\in \overline{\domext}$, $\by\in \domint$, with a possibly different constant $C_{\bB}$.
\end{lem}
The proof of the above result mimics the proof of Lemmas \ref{lem:A1_part1}, \ref{lem:A2_part1}, and thus we leave it to the reader. 
We will also need an analogue of Lemma \ref{lem:A3_part1}. 
\begin{lem}
	\label{lem:A3_part2}
	Let $\sigma_c>0$. Then, there exist $\tilde{\gamma}_{\bB}, \, \tilde{c}_{\bB}>0$, which depend on $\bB$,  s.t.\ for all  $s\in\setCp_{\tilde{\alpha}}$, with  $\tilde{\alpha}=\tilde{\gamma}_{\bB}\sigma_c$,
\begin{align}
	\label{eq:res4_standard_pml_part2}
	\Re(s \sqrt{\hpml(s;\bx, \by)})\geq \tilde{c}_{\bB}\Re s\|\bx-\by\|,\, \text{ whenever }\bx\in \overline{\domext} \text{ and }\by\in \domint. 
\end{align} 
\end{lem} 
\begin{proof} 
	We proceed just like in the poof of Lemma \ref{lem:A3_part1}. 
	We use the identity $$\Re (s\sqrt{\hpml})=\Re s\Re \sqrt{\hpml}-\Im s\Im \sqrt{\hpml}.$$ 
	Next, a lower bound for $\Re \sqrt{\hpml}$ can be obtained from Lemma \ref{lem:A1_part2}, because $\Re \sqrt{\hpml}>\sqrt{\Re\hpml}$, as $\Re\hpml>0$. This yields 
	\begin{align}
		\label{eq:bound0_}
		\Re (s\sqrt{\hpml(s; \bx,\by)})=\Re s\Re \sqrt{\hpml(s; \bx,\by)}-\Im s\Im \sqrt{\hpml(s; \bx,\by)}\geq c_{\bB}^{1/2}\Re s \|\bx-\by\|-|\Im s\Im \sqrt{\hpml(s; \bx,\by)}|.
	\end{align}
 We follow the path of Lemma \ref{lem:A3_part1}, which leads to the counterpart of the bound \eqref{eq:imhpml}: 
	\begin{align}
		\label{eq:imhpml2}
		|\Im \sqrt{\hpml(s; \bx,\by)}|\leq\frac{|\Im \hpml(s; \bx, \by)|}{c_{\bB}^{1/2}\|\bx-\by\|}. 
	\end{align}
	Let us examine $\Im \hpml(s; \bx,\by)$. From the expression \eqref{eq:exprhpml_half}, it follows that 
	\begin{align}
		\nonumber
		\Im \hpml(s; \bx,\by)&=-2\frac{\sigma_c\Im s\xi}{|s|^2}\hat{\bx}^{\top}\bB\bp-\frac{2\sigma_c^2\Re s\Im s}{|s|^4}\xi^2\|\hat{\bx}\|^2_{\bB}.
	\end{align}
Because $\|\bp\|_{\bB}\geq \sqrt{\evmin/\evmax}\|\xi\hat{\bx}\|_{\bB}$ as indicated by \eqref{eq:important_bound2}, we can further bound  
\begin{align*}
	|\Im \hpml(s; \bx,\by)|\leq \tilde{C}_{\bB}\frac{\sigma_c}{|s|}\left(1+\frac{\sigma_c}{|s|}\right)\|\bp\|^2.
\end{align*}
Inserting this bound into \eqref{eq:imhpml2}, and next the resulting bound into \eqref{eq:bound0_}, we obtain that 
\begin{align*}
	\Re(s\hpml(s; \bx,\by))&\geq c_{\bB}^{1/2}\Re s\|\bx-\by\|-\frac{\tilde{C}_{\bB}\sigma_c |\Im s|}{c_{\bB}^{1/2}|s|}\left(1+\frac{\sigma_c}{|s|}\right)\|\bx-\by\|\\
	&\geq c_{\bB}^{1/2}\Re s\|\bx-\by\|\left(1-\frac{\tilde{C}_{\bB}\sigma_c}{c_{\bB}\Re s}\left(1+\frac{\sigma_c}{\Re s}\right)\right).
\end{align*}
This yields the desired bound in the statement of the lemma, cf. the proof of Lemma \ref{lem:A3_part1}. 
\end{proof}
\subsubsection{Bounds for complex-shifted PMLs}
For the complex-shifted PMLs, the counterpart of \eqref{eq:exprhpml_half} reads
\begin{align}
	\label{eq:exprhpml_half2}
	\hpml(s+\gamma; \bx,\by)=	\|\bp\|^2_{\bB}+2\frac{\sigma_c \xi}{s+\gamma}\hat{\bx}^{\top}\bB\bp+\frac{\sigma_c^2}{(s+\gamma)^2}\xi^2\|\hat{\bx}\|^2_{\bB}.
\end{align}
We start with the following auxiliary result, whose proof mimics the one of Lemma \ref{lem:aux_res_part1}, and thus we leave it to the reader. 
\begin{lem}
	\label{lem:aux_res_part2}
	The statement of Lemma \ref{lem:aux_res_part1} holds for $\bx\in \overline{\domext}$, $\by\in \domint$, with possibly different constants $c_{\bB}$ and $\nu_0$. 
\end{lem}
	With the above lemma, we can prove the following counterpart of Proposition \ref{prop:cs_pml_part2}. 
	\begin{prop}
		\label{prop:cs_pml_part2}
		There exist constants $\nu_0,\, C_{\pm},\, c_{\pm} >0$ depending on $\bB$ only, s.t.\ for all $\sigma_c,\gamma>0\colon\, \sigma_c<\nu_0\gamma$, for any $s\in \mathbb{C}^+$, any $\bx\in \overline{\domext}, \, \by\in \domint$, the following bounds hold true:
		\begin{align}
			\label{eq:basic_lower_bound_appendix2}
			&|s\sqrt{\hpml(s+\gamma;\bx,\by)}|\geq C_{-}|s|\|\bx-\by\|,\\
			\label{eq:basic_upper_bound_appendix2}
			& |s\sqrt{\hpml(s+\gamma;\bx,\by)}|\leq C_+|s|\|\bx-\by\|,\\
			\label{eq:basic_lower_bound_appendix2_v2}
			&\Re(s\sqrt{\hpml(s+\gamma;\bx,\by)})\geq (c_+\Re s-c_{-}\sigma_c)\|\bx-\by\|.
		\end{align}
		Finally, if, additionally, $\|\bx\|>\rho=\radpml \evmax/\evmin$,  then a stronger bound holds true:
		\begin{align}
			\label{eq:basic_lower_bound_appendix2_v3}
			\Re \left(s\sqrt{\hpml(s+\gamma;\bx,\by)}\right)\geq c_{+}\Re s\|\bx-\by\|.
		\end{align}
	\end{prop}
The value of $\rho$ in the above is non-optimal, but we do not pursue the goal of optimizing it. 
\begin{proof}
For the proof of the bounds \eqref{eq:basic_lower_bound_appendix2}, \eqref{eq:basic_upper_bound_appendix2}, \eqref{eq:basic_lower_bound_appendix2_v2}, we refer the reader to the analogous proof of Proposition  \ref{prop:cs_pml_part1}. 
Let us now prove the bound \eqref{eq:basic_lower_bound_appendix2_v3}. We use the identity 
		\begin{align}
			\label{eq:hpml_lb2}
			\Re(s\sqrt{\hpml(s+\gamma; \bx,\by)})=\Re s\Re \sqrt{\hpml(s+\gamma; \bx,\by)}-\Im s\Im \sqrt{\hpml(s+\gamma; \bx,\by)}.
		\end{align}
We will prove that $\Im s\Im \sqrt{\hpml(s+\gamma; \bx,\by)}<0$ under the geometric condition stated before \eqref{eq:basic_lower_bound_appendix2_v3}. First let us compute
		\begin{align}
			\nonumber
\Im \hpml(s+\gamma; \bx,\by)=-2\frac{\sigma_c\Im s\xi}{|s+\gamma|^2}\hat{\bx}^{\top}\bB\bp-\frac{2\sigma_c^2(\Re s+\gamma)\Im s}{|s+\gamma|^4}\xi^2\|\hat{\bx}\|^2_{\bB}.
		\end{align}
Let us consider 
\begin{align*}
	\hat{\bx}^{\top}\bB\bp=\hat{\bx}^{\top}\bB(\|\bx\|\hat{\bx}-\|\by\|\hat{\by})\geq \evmin\|\bx\|-\evmax\|\by\|\geq \lambda_{\min}\|\bx\|-\evmax\radpml.
\end{align*}
Therefore, the above quantity is non-negative under the geometric condition of the proposition, and hence $\Im\hpml\Im s\leq 0$, which implies that $\Im\sqrt{\hpml}\Im s\leq 0$. The bound \eqref{eq:hpml_lb2} rewrites \begin{align*}	\Re(s\sqrt{\hpml(s+\gamma; \bx,\by)})&=\Re s\Re \sqrt{\hpml(s+\gamma; \bx,\by)}\geq \Re s\Re \sqrt{\hpml(s+\gamma;\bx,\by)}\\
	&\geq \Re s\sqrt{\Re\hpml(s+\gamma; \bx,\by)}\geq c_{\bB}^{1/2}\Re s\|\bx-\by\|,
	\end{align*}
as per Lemma \ref{lem:aux_res_part2}.
		
%
		\end{proof}
\subsection{Summary of the results}
\label{sec:summary_of_the_bounds}
Below we summarize all the bounds we obtained before. They are evident for $\bx,\by\in \domint$, and follow from previous results in the opposite case. 
\begin{lem}
	\label{lem:A1}
	Let $\sigma_c>0$. Then there exist $\gamma_{\bB}, \, c_{\bB}>0$, which depend on $\bB$, s.t.\ for all  $s\in\setCp_{\alpha}$, with  $\alpha=\gamma_{\bB}\sigma_c$, 
	\begin{align}
		\label{eq:first_bound}
		&\Re \hpml(s; \bx,\by)\geq c_{\bB}\|\bx-\by\|^2,\qquad \text{ for all }\bx,\by\in \setR^2. 
	\end{align}
\end{lem}
\begin{lem}
	\label{lem:A2}
	Let $\sigma_c>0$. Then there exists a constant $C_{\bB}>0$, which depends on $\bB$, s.t. for all  $s\in\setCp$, 
	\begin{align*}
		|\hpml(s; \bx,\by)|\leq C_{\bB}\left(1+\frac{\sigma_c}{|s|}\right)^2\|\bx-\by\|^2, \quad \text{for all }\bx,\by\in \setR^2.
	\end{align*}
\end{lem} 
\begin{lem}
	\label{lem:A3}
	Let $\sigma_c>0$. Then, there exist $\tilde{\gamma}_{\bB}, \, \tilde{c}_{\bB}>0$, which depend on $\bB$,  s.t.\ for all  $s\in\setCp_{\tilde{\alpha}}$, with  $\tilde{\alpha}=\tilde{\gamma}_{\bB}\sigma_c$,
	\begin{align}
		\label{eq:res4_standard_pml}
		\Re(s \sqrt{\hpml(s;\bx, \by)})\geq \tilde{c}_{\bB}\Re s\|\bx-\by\|,\, \text{ for all  }\bx,\by\in \setR^2.
	\end{align}
\end{lem} 

For the complex-scaled PMLs, the corresponding bounds are stated below.
\begin{lem}
	\label{lem:lower_bound}
There exist constants $\nu_0>0$ and $c_{\bB}>0$ depending on $\bB$, s.t.\ for all  $\sigma_c<\nu_0\gamma$, for all $s\in \mathbb{C}^+$,  
	\begin{align}
		\label{eq:bpml2}
		\Re \hpml(s+\gamma; \bx,\by)>c_{\bB}\|\bx-\by\|^2, \quad \text{ for all } \bx,\by\in \mathbb{R}^2.
	\end{align}
\end{lem}
An immediate corollary of the above lemma is the following result. 
\begin{cor}
	\label{cor:Gsigma}
	Under conditions of Lemma \ref{lem:lower_bound}, for all $\bx\neq \by$, the function $s\mapsto G_{\sigma,\gamma}(s+\gamma; \bx,\by)$ is analytic in $\mathbb{C}^+$.
\end{cor}
\begin{proof}
See the reasoning before \eqref{eq:omegainst}.
\end{proof}
\begin{prop}
		\label{prop:cs_pml_orig_appendix}
		There exist constants $\nu,\, C_{\pm},\, c_{\pm}, C_{\bB} >0$ depending on $\bB$ only, s.t.\ for all $\sigma_c,\gamma>0\colon\, \sigma_c<\nu\gamma$, and for all $s\in \mathbb{C}^+$, $h_{\sigma}(s+\gamma; \bx,\by)\notin \mathbb{R}^{-}_0$ for $\bx\neq \by$, and the following bounds hold true.\\
		1. For $(\bx,\by)\in \setR^2\times \setR^2$,
		\begin{align}
			\label{eq:bound1_app}
			C_+|s|\|\bx-\by\|\geq |s\sqrt{\hpml(s+\gamma;\bx,\by)}|\geq C_{-}|s|\|\bx-\by\|.
		\end{align}
			2. For all $(\bx,\by)\in {\overline{\domext}}\times \setR^2$, s.t. $\|\bx-\by\|>\rho=C_{\bB}\radpml$, it holds that 
			\begin{align}
				\label{eq:bound2_app}
				\Re \left(s\sqrt{\hpml(s+\gamma;\bx,\by)}\right)\geq c_{+}\Re s\|\bx-\by\|.
			\end{align} 
			3. For $(\bx,\by)\in \overline{\domext}\times\setR^2$, 
			\begin{align}
				\label{eq:bound3_app}
				\Re (s\sqrt{\hpml(s+\gamma; \bx,\by)})\geq c_{+}\Re s\|\bx-\by\|-c_{-}\sigma_c\|\bx-\by\|.
			\end{align}
\end{prop}
\begin{proof}
The analyticity of $h_{\sigma}$ follows from Lemma \ref{lem:lower_bound}. 

The bounds \eqref{eq:bound1_app} for $\bx,\by\in \domint$ are obvious. For $\bx, \by\in \domext$ they stem from Proposition \ref{prop:cs_pml_part1}, and for $\bx\in \domext$, $\by\in \domint$, they follow from Proposition \ref{prop:cs_pml_part2}. 
The remaining bounds follow for $\bx, \by\in \overline{\domext}$  from Proposition \ref{prop:cs_pml_part1}, and from $\bx\in \overline{\domext}, \by\in\domint$ from Proposition \ref{prop:cs_pml_part2}. Remark in particular that the condition $\|\bx-\by\|>\rho$ for $\bx\in\overline{\domext}$ and $\by\in \domint$ implies that that $\|\bx\|\geq\rho-\|\by\|\geq \rho-\radpml$; if $\rho=\radpml({\evmax}/{\evmin}+1)$, then \eqref{eq:basic_lower_bound_appendix2_v3} holds true.
\end{proof}
\section{Proof of Proposition \ref{prop:fs}}
\label{appendix:Gsigma}
To prove Proposition \ref{prop:fs}, it suffices to show that for all $s\in \setCp_{\alpha}$, 
    	\begin{itemize}
    		 \item the expression in the right-hand side of \eqref{eq:usigma_g} is in $H^1(\setR^2)$;
    		\item the expression in the right-hand side of \eqref{eq:usigma_g} satisfies the PDE \eqref{eq:rad_lpl_v0} in the weak sense.
    	\end{itemize}

\subsection{Proof that the RHS of \eqref{eq:usigma_g} is in $H^1(\setR^2)$}
We will need the following auxiliary bound on $G_{\sigma}$. 
\begin{lem}
	\label{lem:B1}
	Let $s\in \setCp_{\alpha}$, for some $\alpha>0$ sufficiently large. Then, for all $\bx,\,\by\in \setR^{2}$, the following bound holds true, with some constants $C_1, \, C_2, \,c$ depending on $\bB$, 
	\begin{align*}
		&|G_{\sigma}(s; \bx,\by)|\lesssim G^{(1)}_{\sigma}(\|\bx-\by\|)+G^{(2)}_{\sigma}(\|\bx-\by\|), \\ 
		&G^{(1)}_{\sigma}(r):=C_1(1+|\log|s||+|\log r|)\mathbbm{1}_{r\leq r_*},\\ 
		&G^{(2)}_{\sigma}(r):=C_2\mathrm{e}^{-c r\Re s }\mathbbm{1}_{r>r_*}, \quad r_*=(|s|+\sigma_c)^{-1}.
	\end{align*}
In a similar manner, 
	\begin{align*}
	&\|\nabla_{\bx}G_{\sigma}(s; \bx,\by)\|\leq |s|\left(1+\frac{\sigma_c}{|s|}\right)^2\left(\tilde{G}_{\sigma}^{(1)}(\|\bx-\by\|)+G_{\sigma}^{(2)}(\|\bx-\by\|)\right), \quad 
	\tilde{G}_{\sigma}^{(1)}(r):=C_1\frac{1}{|s|r}\mathbbm{1}_{r\leq r_*}.
\end{align*}
\end{lem}
\begin{proof}
\textit{Step 1. A bound on $G_{\sigma}$. }
We start by recalling that $z\mapsto K_0(z)$ is analytic in $\mathbb{C}\setminus \setR^{-}$, and is bounded according to  \eqref{eq:K0Z}, with $C(a)>0$, 
\begin{align}
	\label{eq:K0Z2_appendix}
	|K_0(z)|\leq C(a)\left\{
	\begin{array}{ll}
	\max(1,|\log|z||) & |z|\leq a, \\
		\mathrm{e}^{-\Re z}, & |z|>a.
	\end{array}
	\right.
\end{align}
Let us split $\setR^2\times\setR^2$ into two intersecting subsets, following the results of Lemmas \ref{lem:A1} and \ref{lem:A2}, which will facilitate the use of the above bound:
\begin{align}
	\label{eq:O1_O2}
	\begin{split}
	&\mathcal{O}_1:=\left\{(\bx,\by)\in \setR^2\times\setR^2: \, (|s|+\sigma_c)\|\bx-\by\|\leq 1\right\}, \\
	&\mathcal{O}_2:=\left\{(\bx,\by)\in \setR^2\times\setR^2: \, (|s|+\sigma_c)\|\bx-\by\|> 1\right\}.
	\end{split}
\end{align}
Remark that $\mathcal{O}_1\cup\mathcal{O}_2=\setR^2\times \setR^2$.
In $\mathcal{O}_1$, by Lemma \ref{lem:A2}, $|s\sqrt{\hpml(s; \bx,\by)}|\leq \sqrt{C_{\bB}}$, and from \eqref{eq:K0Z2_appendix} it follows that
\begin{align*}
	|G_{\sigma}(s; \bx,\by)|\leq C\left(|\log |s\hpml||+1\right)\leq C\left(|\log  |s||+|\log\|\bx-\by\||+1\right), 
\end{align*} 
where the last bound follows from Lemma \ref{lem:A1}. 
In $\mathcal{O}_2$, with Lemma \ref{lem:A1}, $|s\sqrt{\hpml(s; \bx,\by)}|\geq c_{\bB}\frac{|s|}{|s|+\sigma_c}$. Because we can choose $s\in \setCp_{\alpha}$ with $\alpha$ sufficiently large, we remark that for $\alpha\geq\sigma_c$, $\frac{|s|}{|s|+\sigma_c}> 1/2$. This enables us to use the bound \eqref{eq:K0Z2_appendix}, which, in turn, yields 
\begin{align*}
	|G_{\sigma}(s; \bx,\by)|\leq C\mathrm{e}^{-\Re (s\sqrt{\hpml(s; \bx,\by)})}\leq C\mathrm{e}^{-\tilde{c}_{\bB}\Re s\|\bx-\by\|}, 
\end{align*}
 cf.\ Lemma \ref{lem:A3}. This allows to obtain the bound in the statement of the lemma.
  
\textit{Step 2. A bound on $\nabla_{\bx}G_{\sigma}(s; \bx,\by)$. }
Let us compute, with \cite[10.29]{nist}, where $K_1$ is a modified Bessel function, 
\begin{align}
	\label{eq:nablaxG}
	\nabla_{\bx}G_{\sigma}(s;\bx,\by)=-\frac{s}{2\pi}(\operatorname{det}\bB)^{1/2}K_1(s\sqrt{\hpml(s;\bx,\by)})\frac{\nabla_{\bx}\hpml(s;\bx,\by)}{2\sqrt{\hpml(s;\bx,\by)}}.
\end{align}
Using the definition  \eqref{eq:js} of $\Js$ we obtain
\begin{align}
	\label{eq:nbx}
	\nabla_{\bx}\hpml(s;\bx,\by)=\left\{
	\begin{array}{ll}
		2\bB(\bx-\by_{\sigma}), & \bx\in \domint, \\
		2\Js(s,\bx)\bB(\bx_{\sigma}-\by_{\sigma}), & \bx\in \overline{\domext}. 
	\end{array}	
	\right.
\end{align}
Let us bound $\nabla_{\bx}\hpml(s; \bx,\by)$ in different regions. First of all, we have
\begin{align*}
	\|\bx_{\sigma}-\by_{\sigma}\|\leq \left\{
	\begin{array}{ll}
	\|\bx-\by\|, & \bx,\by\in \domint,\\
	2\left(1+\frac{\sigma_c}{|s|}\right)\|\bx-\by\|, & \bx\in \mathbb{R}^2, \, \by\in \overline{\domext}.
	\end{array}
	\right.
\end{align*}
Indeed, when $\bx\in \domint$ and $\by\in \overline{\domext}$, it holds that 
\begin{align*}
	\|\bx-\by_{\sigma}\|\leq \|\bx-\by\|+\frac{\sigma_c}{|s|}(\|\by\|-\radpml)\|\hat{\by}\|\leq \left(1+\frac{\sigma_c}{|s|}\right)\|\bx-\by\|,
\end{align*}
with the last bound holding true because  $\|\by\|-\radpml=\operatorname{dist}(\by,\domint)\leq \|\by-\bx\|$.
In a similar manner, when $\bx, \, \by\in \overline{\domext}$, we use \eqref{eq:ppml} and the bound before \eqref{eq:important_bound1} (which gives the factor $2$ in the final bound).

Therefore, from \eqref{eq:nbx} and the definition  \eqref{eq:js} of $\Js$ (which implies that $|J_{\sigma}(\bx)|<1+\frac{\sigma_c}{|s|}$) we conclude that 
\begin{align*}
	\|\nabla_{\bx}\hpml(s; \bx,\by)\|\leq C_{\bB}\left(1+\frac{\sigma_c}{|s|}\right)^2\|\bx-\by\|. 
\end{align*}
Plugging in the above bound into \eqref{eq:nablaxG} and using Lemma \ref{lem:A1} to bound $|\sqrt{\hpml(s; \bx,\by)}|$ from below yields 
\begin{align*}
	\|G_{\sigma}(s; \bx, \by)\|\leq c(\bB)|s|\left(1+\frac{\sigma_c}{|s|}\right)^2\left|K_1(s\sqrt{\hpml(s; \bx,\by)})\right|, \quad c(\bB)>0. 
\end{align*}
%
%
%
%
%
%
It remains to bound $K_1$.  From \cite{nist}[10.30, 10.25.3] it follows that
\begin{align}
	\label{eq:K1Z}
	|K_1(z)|\leq C\left\{
	\begin{array}{ll}
	|z|^{-1} & |z|\leq 1/2, \\
		\mathrm{e}^{-\Re z}, & |z|>1/2.
	\end{array}
	\right.
\end{align}
We split $\setR^2\times \setR^2$ as per \eqref{eq:O1_O2} and proceed like in the derivation of the bound for $|G_{\sigma}|$. 
Using Lemmas \ref{lem:A1}, \ref{lem:A2}, \ref{lem:A3}, we obtain the desired bound:
	\begin{align*}
	\|\nabla_{\bx}G_{\sigma}(s; \bx,\by)\|\leq C(\bB) |s|\left(1+\frac{\sigma_c}{|s|}\right)^2\left(\frac{1}{|s|\|\bx-\by\|}\mathbbm{1}_{\mathcal{O}_1}+\mathrm{e}^{-c\Re s\|\bx-\by\|}\mathbbm{1}_{\mathcal{O}_2}\right).
	\end{align*}
\end{proof}
With these two bounds, it is straightforward to obtain the desired regularity result. 
\begin{lem}
	Let $\hat{F}\in L^2(\setR^2)$. There exists $\alpha>0$, s.t. for all $s\in \setCp_{\alpha}$, it holds that  $\bx\mapsto\hat{u}^{\sigma}(s,\bx)=\int_{\setR^2}G_{\sigma}(s;\bx,\by)\hat{F}(\by)d\by\in H^1(\setR^2)$.
\end{lem}
\begin{proof}
By Lemma \ref{lem:B1},
\begin{align*}
	|\hat{u}^{\sigma}(s,\bx)|\lesssim \left|G_{\sigma}^{(1)}*|\hat{F}|\right|+\left|G_{\sigma}^{(2)}*|\hat{F}|\right|,
\end{align*}
and with the use of the Young convolution inequality we conclude that 
\begin{align*}
\|\hat{u}^{\sigma}(s)\|_{L^2}\lesssim \left(\|G_{\sigma}^{(1)}\|_{L^1}+\|G_{\sigma}^{(2)}\|_{L^1}\right)\|\hat{F}\|_{L^2}.
\end{align*}
The right hand side is obviously finite because $G_{\sigma}^{(1)}$ is $L^1(\mathbb{R}^2)$, compactly supported and $G_{\sigma}^{(2)}$ is exponentially decreasing. Similarly, we can show that $\nabla_{\bx}\hat{u}^{\sigma}(s)\in L^2(\setR^2)$. 
\end{proof}

\subsection{Proof that the RHS of \eqref{eq:usigma_g} satisfies \eqref{eq:rad_lpl_v0} weakly}
The next lemma shows that 
\begin{lem}	
The expression	\eqref{eq:usigma_g} satisfies the PDE \eqref{eq:lpl_u} with the right-hand side $\hat{F}$ weakly.
\end{lem}
\begin{proof}
We proceed like in 	\cite{BermudezHervellaNPrietoRodriguez:08}. Remark that $\by\mapsto G(s; \bx, \, \by)\in H^1_{loc}(\setR^2)$, cf. Lemma \ref{lem:B1}, for all $\bx\in \setR^2$. Our goal is to prove that for all $\varphi\in \mathcal{D}(\setR^2)$, 
\begin{align*}
	&\int_{\setR^2}\bA_{\sigma}(s,\bx)\left(\int_{\setR^2}\nabla_{\bx}G(s; \bx,\by)\hat{F}(\by)d\by\right)\, \nabla_{\bx}\varphi(\bx)d\bx\\
	&+s^2	\int_{\setR^2}\left(\int_{\setR^2}G(s; \bx,\by)\hat{F}(\by)d\by \right)\varphi(\bx)\operatorname{det}\Js(s,\bx) d{\bx}=\int_{\setR^2}\hat{F}(\by)\varphi(\by)d\by. \quad \bx\in \setR^2. 
\end{align*}
The above can be rewritten as follows
\begin{align*}
	\int_{\setR^2}\hat{F}(\by)\left(\int_{\setR^2}\bA_{\sigma}(s,\bx)\nabla_{\bx}G(s; \bx,\by)\nabla_{\bx}\varphi(\bx)d\bx+s^2\int_{\setR^2}\operatorname{det}\Js(s,\bx)G(s; \bx,\by)\varphi(\bx)d\bx \right)d\by=\int_{\setR^2}\hat{F}(\by)\varphi(\by)d\by.
\end{align*}
 Therefore, it is sufficient to prove that, for all $\varphi\in \mathcal{D}(\setR^2)$,   
\begin{align}
	\label{eq:phiG}
	\int_{\setR^2}\bA_{\sigma}(s,\by)\nabla_{\by}G(s; \by,\bx)\nabla_{\by}\varphi(\by)d\by+s^2	\int_{\setR^2}G(s; \by,\bx)\varphi(\by)\operatorname{det}\Js(s,\by) d{\by}=\varphi(\bx), \quad \bx\in \setR^2.
\end{align}
The first integral in the right-hand side can be rewritten as follows:
\begin{align*}
	\int_{\setR^2}\bA_{\sigma}(s,\by)\nabla_{\by}G(s; \by,\bx)\nabla_{\by}\varphi(\by)d\by&=	\lim\limits_{\varepsilon\rightarrow 0+}\int_{\setR^2\setminus B_{\varepsilon}(\bx)}\bA_{\sigma}(s,\by)\nabla_{\by}G(s; \by,\bx)\nabla_{\by}\varphi(\by)d\by\\
	&=-\lim\limits_{\varepsilon\rightarrow 0+}\int_{\partial B_{\varepsilon}(\bx)}\bA_{\sigma}(s,\by)\nabla_{\by}G(s; \by,\bx)\cdot \boldsymbol{n}_{\by}\varphi(\by)d\Gamma_{\by}\\
	&-
	\lim\limits_{\varepsilon\rightarrow 0+}\int_{\setR^2\setminus B_{\varepsilon}(\bx)}\operatorname{div}_{\by}\left(\bA_{\sigma}(s,\by)\nabla_{\by}G(s; \by,\bx)\right)\varphi(\by)d\by, 
\end{align*}
where $\boldsymbol{n}_{\by}$ is the normal to $\partial B_{\varepsilon}(\bx)$ pointing into $\mathbb{R}^2\setminus B_{\varepsilon}(\bx)$. 

Remark that $\by\mapsto G(s; \by,\bx)$ is smooth outside of $\by=\bx$, and, moreover, by construction, satisfies the equation \eqref{eq:basigma_main} with $\hat{f}=0$ strongly whenever $\by\neq \bx$. Therefore, comparing to \eqref{eq:phiG}, we have to prove that 
\begin{align*}
	\lim\limits_{\varepsilon\rightarrow 0+}\int_{\partial B_{\varepsilon}(\bx)}\bA_{\sigma}(s,\by)\nabla_{\by}G(s; \by,\bx)\cdot \boldsymbol{n}_{\by}\varphi(\by)d\Gamma_{\by}=-\varphi(\bx). 
\end{align*}
Let us introduce 
\begin{align*}
	I_{\varepsilon}(\bx)&=\int_{\partial B_{\varepsilon}(\bx)}\bA_{\sigma}(s,\by)\nabla_{\by}G(s; \by,\bx)\cdot \boldsymbol{n}_{\by}d\Gamma_{\by}, \\
	E_{\varepsilon}(\bx)&=\int_{\partial B_{\varepsilon}(\bx)}\bA_{\sigma}(s,\by)\nabla_{\by}G(s; \by,\bx)\cdot \boldsymbol{n}_{\by}(\varphi(\by)-\varphi(\bx))d\Gamma_{\by},
\end{align*}
and rewrite
\begin{align*}
&\int_{\partial B_{\varepsilon}(\bx)}\bA_{\sigma}(s,\by)\nabla_{\by}G(s; \by,\bx)\cdot \boldsymbol{n}_{\by}\varphi(\by)d\Gamma_{\by}=\varphi(\bx)I_{\varepsilon}(\bx)+E_{\varepsilon}(\bx).
\end{align*}
Our goal is to prove that, as $\varepsilon\rightarrow 0+$, $I_{\varepsilon}(\bx)\rightarrow -1$, while $E_{\varepsilon}(\bx)\rightarrow 0$. For this we will use an explicit expression
\begin{align}
	\label{eq:derG}
\nabla_{\by}G(s; \by, \bx)=-\frac{s}{2\pi}(\operatorname{det}\bB)^{1/2}K_1(s\sqrt{\hpml(s;\by,\bx)})\frac{\nabla_{\by}\hpml(s;\by,\bx)}{2\sqrt{\hpml(s;\by,\bx)}}, \quad h_{\sigma}(s; \by, \bx)=(\bx_{\sigma}-\by_{\sigma})^{\top}\bB(\bx_{\sigma}-\by_{\sigma}). 
\end{align}
In the above, cf. \eqref{eq:nbx}, 
\begin{align}
	\label{eq:derH}
	\nabla_{\by}\hpml(s;\by,\bx)=\left\{
	\begin{array}{ll}
		2\bB(\by-\bx_{\sigma}), & \by\in \domint, \\
		2\Js(s,\by)\bB(\by_{\sigma}-\bx_{\sigma}), & \by\in {\domext}. 
	\end{array}	
	\right.
\end{align}
Additionally, we will need asymptotics of modified Bessel functions \cite[10.31]{nist}: 
\begin{align}
	\label{eq:asK1}
	K_1(z)=z^{-1}+O(\log z), \quad |z|\rightarrow 0.
\end{align}
We will first prove the corresponding result for $\bx\in \domext$; the proof for the case $\bx\in \domint$ will be basically the same (hence we skip it here), and next consider the case $\bx\in \overline{\domint}\cap\overline{\domext}=\Sigma$. \\
\textit{Case 1. $\bx \in \domext$. }\\
\textit{Step 1. Proof that $I_{\varepsilon}(\bx)\rightarrow -1$ as $\varepsilon\rightarrow 0+$. }\\
Using the asymptotics of the McDonald function $K_1$ \eqref{eq:asK1} and the expression \eqref{eq:derG}, we conclude that, as $\varepsilon\rightarrow 0+$, 
\begin{align}
	\label{eq:nablaG}
	\nabla_{\by}G(s; \by,\bx)=-\frac{(\operatorname{det}\bB)^{1/2}}{4\pi \hpml(s; \by,\bx)}\nabla_{\by}\hpml(s;\by,\bx)+r_{1}(\bx,\by), \quad \|r_{1}(\bx,\by)\|\leq C_1\left\|\frac{\nabla_{\by}\hpml(s; \bx,\by)}{\sqrt{\hpml(s; \bx,\by)}}\log \hpml(s; \bx,\by)\right\|.
\end{align}
Next, we parametrize $\by\in \partial B_{\varepsilon}(\bx)$ as follows:
\begin{align*}
	\by=\bx+\varepsilon\hat{\bz}, \quad \hat{\bz}\in \mathbb{S}^2.
\end{align*}
With this parametrization, for sufficiently small $\varepsilon$, we have that (where we use the fact that $\by\mapsto \by_{\sigma}$ is regular inside $\domext$)
\begin{align*}
	\by_{\sigma}-\bx_{\sigma}&=\by\dtpml(\|\by\|)-\bx\dtpml(\|\bx\|)=\left(\bx\dtpml(\|\bx\|)+\Js(s,\bx)\varepsilon\hat{\bz}+O(\varepsilon^2)\right)-\bx\dtpml(\|\bx\|)\\
	&=
\Js(s,\bx)\varepsilon\hat{\bz}+r_{1,\varepsilon}(\hat{\bz}), \quad \|r_{1,\varepsilon}\|_{\infty}\leq C_1\varepsilon^2.
\end{align*}
With \eqref{eq:derH}, the equation \eqref{eq:nablaG} rewrites further 
\begin{align*}
	\nabla_{\by}G(s;\by,\bx)=-\frac{(\operatorname{det}\bB)^{1/2}}{2\pi \varepsilon \hat{\bz}^{\top}\Js^{\top}(s,\bx)\bB\Js(s,\bx)\hat{\bz}}\Js(s,\by)\bB \Js(s,\bx)\hat{\bz}+r_{2,\varepsilon}(\bx,\by), \quad \|r_{2,\varepsilon}\|_{\infty}\leq C_2\left|\log\varepsilon\right|, \quad C_2>0.
\end{align*}
Let us now consider $\bA_{\sigma}(s,\by)\nabla_{\by}G(s; \by,\bx)$. 
Recall that $\bA_{\sigma}(s,\by)=\Js^{-1}(s,\by)\bA\Js(s,\by)^{-\top}\det\Js(s,\by)$. Since it holds that $\Js$ is a symmetric matrix, $\bA\bB=\operatorname{Id}$,  $\Js^{-1}(\by)=\Js^{-1}(\bx)+r(\by,\bx)$, $\|r\|_{\infty}<\varepsilon$, 
we have 
\begin{align*}
\bA_{\sigma}(s,\by)\Js^{\top}(s,\by)\bB\Js(s,\bx)=\operatorname{det}\Js(s,\by)\operatorname{Id}+O(\varepsilon)=\operatorname{det}\Js(s,\bx)\operatorname{Id}+O(\varepsilon).
\end{align*}
We then obtain
\begin{align*}
	\bA_{\sigma}(s,\by)\nabla_{\by}G(s; \by,\bx)=-\frac{(\det\bB)^{1/2}\operatorname{det}\Js(s,\bx)}{2\pi \varepsilon\hat{\bz}^{\top}\Js^{\top}(s,\bx)\bB\Js(s,\bx)\hat{\bz}}\hat{\bz}+r_{3,\varepsilon}(\bx,\by), \quad \|r_{3,\varepsilon}\|_{\infty}\leq C_3|\log\varepsilon|.
\end{align*}
Applying Lemma \ref{lem:integral} to the above yields the desired identity $I_{\varepsilon}\rightarrow -1$ as $\varepsilon\rightarrow 0+$. \\
\textit{Step 2. Proof that $E_{\varepsilon}(\bx)\rightarrow 0$ as $\varepsilon\rightarrow 0+$. }
The result follows from the arguments similar to the above, by remarking that $\|\varphi(.)-\varphi(\bx)\|_{\infty}\leq C_{\varphi}\varepsilon$.\\
\textbf{Case 2. $\bx\in \domint$. } The result follows from Case 1 with $\sigma_c=0$. \\
\textbf{Case 3. $\bx\in \bdpml$. }
\textit{Step 1. Proof that $I_{\varepsilon}(\bx)\rightarrow -1$ as $\varepsilon\rightarrow 0+$. }
Remark that there is no issues in defining $G_{\sigma}(s; \by,\bx)$ in this case, since $\bx_{\sigma}(s)\equiv \bx$. 
\begin{figure}
	\centering
	\includegraphics[width=0.15\textwidth]{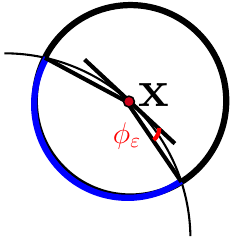}
	\caption{Splitting of $\partial B_{\varepsilon}(\bx)$ into two parts. The part of the curve $\partial B_{\varepsilon}(\bx)\cap \domint$ is marked in blue. }
	\label{fig:three_parts}
\end{figure}
We split the integral $I_{\varepsilon}(\bx)$ into two parts, see Figure \ref{fig:three_parts}, 
\begin{align*}
I_{\varepsilon}(\bx)=I^{\sigma}_{\varepsilon}(\bx)+I^{0}_{\varepsilon}(\bx), 
\end{align*}
where 
\begin{align*}
	I^{\sigma}_{\varepsilon}(\bx)=\int_{\partial B_{\varepsilon}(\bx)\cap \domext}\bA_{\sigma}(\by)\nabla_{\by}G(s;\by,\bx)\cdot \boldsymbol{n}_{\by}d\by, \quad 	I^{0}_{\varepsilon}(\bx)=\int_{\partial B_{\varepsilon}(\bx)\cap \domint}\bA\nabla_{\by}G(s;\by,\bx)\cdot \boldsymbol{n}_{\by}d\by. 
\end{align*}
We start by considering $I_{\varepsilon}^0(\bx)$: 
\begin{align*}
	I^0_{\varepsilon}(\bx)=\int_{\partial B_{\varepsilon}(\bx)\cap \domint}\bA\nabla_{\by}G(s;\by,\bx)\cdot \boldsymbol{n}_{\by}d\by.
\end{align*}
We use the expressions \eqref{eq:nablaG} and \eqref{eq:derH}, which yield
\begin{align*}
	\bA\nabla_{\by}G(s; \by, \bx)=-\frac{(\operatorname{det}\bB)^{1/2}\bB(\by-\bx)}{2\pi\|\by-\bx\|_{\bB}^2}+r_1(\bx,\by), \quad |r_1(\bx,\by)|\leq C_1|\log\|\bx-\by\|_{\bB}|, \quad C_1>0.
\end{align*}
For $\by=\bx+\varepsilon\hat{\bz}$, the above rewrites 
\begin{align}
	\bA\nabla_{\by}G(s; \by, \bx)=-\frac{(\operatorname{det}\bB)^{1/2}\bB\hat{\bz}}{2\pi\varepsilon\hat{\bz}^{\top}\bB\hat{\bz}}+r_1(\bx,\by), \quad \|r_1\|_{\infty}\leq \tilde{C}_1|\log\varepsilon|.
\end{align}
We then rewrite, with $\hat{\bz}=\hat{\bz}_{\phi}=(\cos\phi,\sin\phi)^\top$, cf. Figure \ref{fig:three_parts}, 
\begin{align*}
	I_{\varepsilon}^0(\bx)=-\int_{-\pi+\phi_{\varepsilon}}^{-\phi_{\varepsilon}}\left(\frac{(\operatorname{det}\bB)^{1/2}\bB\hat{\bz}_{\phi}}{2\pi\hat{\bz}^{\top}_{\phi}\bB\hat{\bz}_{\phi}}+\varepsilon r_1(\bx,\bx+\varepsilon\hat{\bz}_{\phi})\right)d\phi,
\end{align*}
where $\phi_{\varepsilon}\rightarrow 0$, as $\varepsilon\rightarrow 0$. Using the Lebesgue's dominated convergence theorem and Lemma \ref{lem:integral}, we conclude that, as $\varepsilon\rightarrow 0$, $I_{\varepsilon}^0\rightarrow -\frac{1}{2}$.

To study $I^{\sigma}_{\varepsilon}(\bx)$, we rewrite \eqref{eq:nablaG}:
\begin{align*}
	\nabla_{\by}G(s; \by, \bx)=-\frac{(\operatorname{det}\bB)^{1/2}\Js(s,\by)\bB(\by_{\sigma}-\bx)}{2\pi (\by_{\sigma}-\bx)^{\top}\bB(\by_{\sigma}-\bx)}+\tilde{r}_1(\bx,\by), \quad \|\tilde{r}_1\|_{\infty}\leq \tilde{C}_1|\log\varepsilon|. 
\end{align*}
With $\bA_{\sigma}(s,\by)=\Js^{-1}(s,\by)\bA \Js(s,\by)^{-1}\det\Js(s,\by)$, we have that 
\begin{align}
\bA_{\sigma}(s,\by)\nabla_{\by}G(s; \by, \bx)=-\frac{(\operatorname{det}\bB)^{1/2}\det\Js(s,\by)\Js^{-1}(s,\by)(\by_{\sigma}-\bx)}{2\pi (\by_{\sigma}-\bx)^{\top}\bB(\by_{\sigma}-\bx)}+\tilde{r}_1(\bx,\by), \quad \|\tilde{r}_1\|_{\infty}\leq \tilde{C}_1|\log\varepsilon|. 
\end{align}
Next, for $\by=\bx+\varepsilon\hat{\bz}$, $\hat{\bz}\in \mathbb{S}^2$, we have that (where we use the regularity of the mapping $\by\mapsto \by_{\sigma}$ in ${\domext}$): 
\begin{align*}
	\by_{\sigma}-\bx=\Js(s,\by)\varepsilon\hat{\bz}+r_{1,\varepsilon}(\bx,\by), \quad \|r_{1,\varepsilon}\|_{\infty}\leq \varepsilon^2.
\end{align*}
This allows to rewrite 
\begin{align}
	\bA_{\sigma}(s,\by)\nabla_{\by}G(s; \by, \bx)=-\frac{(\operatorname{det}\bB)^{1/2}\det\Js(s,\by)\hat{\bz}}{2\pi \varepsilon \hat{\bz}^{\top}\Js(s,\by)^{\top}\bB\Js(s,\by)\hat{\bz}}+r_{2,\varepsilon}(\bx,\by), \quad \|r_{2,\varepsilon}\|_{\infty}\leq \tilde{C}_1|\log\varepsilon|. 
\end{align}
We rewrite further, with $\Pi_{\by}$ being an orthogonal (Eucledian) projection on $\by$ and $\Pi_{\by}^{\perp}=\operatorname{Id}-\Pi_{\by}$, 
\begin{align*}
	\Js(s,\by)=\dpml(s,{\bx+\varepsilon\hat{\bz}}) \Pi_{\bx+\varepsilon\hat{\bz}}+\dtpml(s,\bx+\varepsilon\hat{\bz}) \Pi^{\perp}_{\bx+\varepsilon\hat{\bz}}.
\end{align*}
As $\varepsilon\rightarrow 0$, we have that 
\begin{align}
	&\Pi_{\bx+\varepsilon\hat{\bz}}=\Px+R^{\varepsilon,1}_{\bx,\hat{z}},\qquad  \Pi^{\perp}_{\bx+\varepsilon\hat{\bz}}=\Pxp+R^{\varepsilon,2}_{\bx,\hat{z}},\qquad \|R^{\varepsilon,j}_{\bx,\hat{z}}\|\lesssim \varepsilon,\, j\in \{1,2\},\\
	&\dpml(s, \bx+\varepsilon\hat{\bz})=1+\frac{\sigma_c}{s}, \quad 
	\dtpml(s,\bx+\varepsilon\hat{\bz})=1+r_{3,\varepsilon}(\bx,\hat{\bz}), \quad |r_{3,\varepsilon}(\bx,\hat{\bz})|=O(\varepsilon).
\end{align}
Therefore, with the matrix $J$ depending on $s$ and $\bx$ only, we have the decomposition
\begin{align*}
	\Js(s,\by)=J(s,\bx)+R(s,\bx,\by), \quad \|R(s,\bx,\by)\|_{\infty}<\varepsilon,
\end{align*}
and 
\begin{align*}
	\bA_{\sigma}(s,\by)\nabla_{\by}G(s; \by, \bx)=-\frac{(\operatorname{det}\bB)^{1/2}\det J(s,\bx)}{\varepsilon \hat{\bz}^{\top}J(s,\bx)^{\top}\bB J(s,\bx)\hat{\bz}}\hat{\bz}+r_{3,\varepsilon}(\bx,\by), \quad \|r_{3,\varepsilon}\|_{\infty}\leq C_3|\log\varepsilon|. 	
\end{align*}
It remains to compute the final integral, with $\hat{\bz}=\hat{\bz}_{\phi}=(\cos\phi,\sin\phi)^{\top}$, 
\begin{align*}
	I_{\varepsilon}^{\sigma}(\bx)=-\int_{-\phi_{\varepsilon}}^{\pi+\phi_{\varepsilon}}\frac{(\operatorname{det}\bB)^{1/2}\det J(s,\bx)}{\hat{\bz}^{\top}_{\phi}J(s,\bx)^{\top}\bB J(s,\bx)\hat{\bz}_{\phi}}d\phi+O(\varepsilon\log\varepsilon).
\end{align*}
As before, $\phi_{\varepsilon}\rightarrow 0$, as $\varepsilon\rightarrow 0$. 
With the dominated convergence theorem and Lemma \ref{lem:integral}, we conclude that $\lim\limits_{\varepsilon\rightarrow 0}I_{\varepsilon}^{\sigma}=-\frac{1}{2}$, and hence the desired result.

\textit{Step 2. Proof that $E_{\varepsilon}(\bx)\rightarrow 0$ as $\varepsilon\rightarrow 0+$. }
The result follows from the arguments similar to the above, by remarking that $\|\varphi(.)-\varphi(\bx)\|_{\infty}\leq C_{\varphi}\varepsilon$.

\end{proof}

\begin{lem}
	\label{lem:integral}
For all $s\in \setCp$, $\bx\in \domext$, it holds that
	\begin{align*}
		\int_{\mathbb{S}^2}\frac{1}{\hat{\bz}^{\top} \Js^{\top}(s,\bx) \bB \Js(s,\bx)\hat{\bz}}d\hat{\bz}=\frac{2\pi\sqrt{\operatorname{det\bA}}}{\operatorname{det}\Js(s,\bx)}, 
	\end{align*}
and, moreover, 
\begin{align*}
	\int_{\mathbb{S}^2: \, \phi\in (0, \pi)}\frac{1}{\hat{\bz}^{\top}_{\phi} \Js^{\top}(s,\bx) \bB \Js(s,\bx)\hat{\bz}_{\phi}}d\hat{\bz}_{\phi}=\frac{\pi\sqrt{\operatorname{det\bA}}}{\operatorname{det}\Js(s,\bx)}. 
\end{align*}
\end{lem}
\begin{proof}
The function $s\mapsto f(s; \hat{\bz})=(\hat{\bz}^{\top} \Js^{\top}(s,\bx) \bB \Js(s,\bx)\hat{\bz})^{-1}$ is analytic for $s\in \setCp$ with $\Re s$ sufficiently large, and for all $\hat{\bz}\in \mathbb{S}^2$. By the dominated convergence theorem, so is the function defined by $s\mapsto \int_{\mathbb{S}^2}f(s; \hat{\bz})d\hat{\bz}$. It remains to compute this expression for $s>0$ sufficiently large and use the analytic continuation. We assume for now that $s>0$, and thus the matrix $\Js(s,\bx)$ is positive definite. 
Denoting by $\mathbf{H}:=\bB^{1/2}\Js(s,\bx)$, we obtain
\begin{align*}
	\int_{\mathbb{S}^2}\frac{1}{\hat{\bz}^{\top} \Js^{\top}(s,\bx) \bB \Js(s,\bx)\hat{\bz}}d\hat{\bz}=\int_{\mathbb{S}^2}\frac{1}{\|\mathbf{H}\hat{\bz}\|^2}d\hat{\bz}. 
\end{align*}
Because $\mathbf{H}$ is symmetric positive definite, we have that $\mathbf{H}=\mathbf{U}\operatorname{diag}(h_1,h_2)\mathbf{U}^{\top}$, with $h_i>0$ and $\mathbf{U}$ is an orthogonal matrix (and thus is either a rotation 
 or a reflection matrix
 ). Therefore, 
\begin{align*}
	\int_{\mathbb{S}^2}\frac{1}{\hat{\bz}^{\top} \Js^{\top}(s,\bx) \bB \Js(s,\bx)\hat{\bz}}d\hat{\bz}=\int_{0}^{2\pi}\frac{1}{h_1^2\cos^2\varphi+h_2^2\sin^2\varphi}d\varphi=\frac{2\pi}{h_1h_2}=\frac{2\pi}{\operatorname{det}\mathbf{H}},
\end{align*}
where the last identity follows by rewriting the integral as an integral of a rational function on the unit circle and using the calculus of residues. 
Finally, the respective result for the integral over the half-circle follows by parity considerations. 
\end{proof}

\section{Auxiliary technical results}
\label{appendix:auxiliary_technical_results}
\begin{lem}
	\label{lem:upperbound_cxx_cxy}
Assume that $\bx,\by\in \setR^2\setminus \{0\}$. Then, with $\gamma_{xx}=\hat{\bx}^{\top}\bB\hat{\bx}$, $\gamma_{xy}=\hat{\bx}^{\top}\bB\hat{\by}$, we have that  
	\begin{align*}
		\left| \frac{\gamma_{xy}}{\gamma_{xx}} \right|\leq  \mu_*=\frac{\evmin+\evmax}{2\sqrt{\evmax\evmin}}.
	\end{align*}
\end{lem} 
\begin{proof}
	We are looking for the maximum of 
	\begin{align*}
		(\hat{\mathbf{x}}, \hat{\mathbf{y}})\mapsto F(\hat{\mathbf{x}}, \hat{\mathbf{y}})=\frac{\hat{\mathbf x}^\top\mathbf B\hat{\mathbf y}}{\hat{\mathbf x}^\top\mathbf B\hat{\mathbf x}}. 
	\end{align*}
	Taking $\hat{\bx}=(\cos\alpha,\sin\alpha)$ and $\hat{\by}=(\cos\beta, \sin\beta)$ in the basis of the eigenvectors of $\mathbf B$ yields 
	\begin{align*}
		F(\hat{\mathbf{x}}, \hat{\mathbf{y}})=\frac{\evmin\cos\alpha\cos\beta+\evmax\sin\alpha\sin\beta}{\evmin\cos^2\alpha+\evmax\sin^2\alpha}=\frac{\cos\alpha\cos\beta+\mu\sin\alpha\sin\beta}{\cos^2\alpha+\mu\sin^2\alpha}, \, \mu:=\evmax\evmin^{-1}. 
	\end{align*}
	The critical points of $F$ are given by 
	\begin{align*}
          (-\sin\alpha\cos\beta+\mu\cos\alpha\sin\beta)(\cos^2\alpha+\mu\sin^2\alpha)-(\cos\alpha\cos\beta+\mu\sin\alpha\sin\beta)(2\mu-2)\cos\alpha\sin\alpha&=0,\\
                -\cos\alpha\sin\beta+\mu\sin\alpha\cos\beta&=0. 
	\end{align*}
Replacing in the first equation $\sin\beta$ re-expressed from the second equation yields the following identity:
	\begin{align*}
          \sin\alpha\cos\beta(\mu^2-1)(\cos^2\alpha+\mu\sin^2\alpha)-(2\mu-2)\cos\beta\sin\alpha\left(\cos^2\alpha+\mu^2\sin^2\alpha\right)&=0, \text{ or }\\
                \sin\alpha\,\cos\beta(\mu-1)^2(\cos^2\alpha-\mu\sin^2\alpha)&=0.
	\end{align*}
	It is then easy to check that $F(\hat{\bx},\hat{\by})=1$ if $\sin\alpha\cos\beta=0$. On the other hand, if $\cos\alpha=\pm \mu^{1/2}\sin\alpha$, we have that $\mu^{1/2}\cos\beta=\pm \sin\beta$, and also  
	\begin{align*}
		F(\hat{\bx}, \hat{\by})\in \left\{\frac{1\pm \mu}{2\sqrt{\mu}}, \frac{-1\pm \mu}{2\sqrt{\mu}}\right\}.
	\end{align*}
\end{proof}
\begin{lem}
	\label{lemma:fx}
	Assume that $\bx: \, \|\bx\|=1$, $\bc=\bx-\by$, $\|\by\|<1$ and $\bx^\top\bB\bc<0$. Then
	\begin{align*}
		f(\bx,\bc)=\frac{(\bx^\top \bB \bc)^2}{(\bx^\top\bB\bx)\bc^\top\bB \bc}\leq \frac{(\evmax-\evmin)^2}{(\evmax+\evmin)^2}
	\end{align*} 
\end{lem}
\begin{proof}
Without loss of generality, let $\bB=\operatorname{diag}(\evmax, \evmin)$, and let $\hat{\bx}=(\cos\alpha,\sin\alpha)$, $\hat{\bc}=(\cos\beta,\sin\beta)$. Because $\bx=\bc+\by$, with $\|\by\|\leq 1$, we have that \footnote{Indeed, in the complex number notation $\mathrm{e}^{i\alpha}=\|\bc\|\mathrm{e}^{i\beta}+\by$, and we have that $|\mathrm{e}^{i\alpha}-\|\bc\|\mathrm{e}^{i\beta}|<1$ which amounts to $|\mathrm{e}^{i(\alpha-\beta)}-\|\bc\||<1$, for some $\|\bc\|>0$. Further computations yield $1-2\|\bc\|\cos(\alpha-\beta)+\|\bc\|^2<0$. This in turn requires that $\cos(\alpha-\beta)>0$, and hence the conclusion. }
\begin{align}
	\label{eq:cond1}
	\cos(\alpha-\beta)=\cos\alpha\cos\beta(1+\tan\alpha\tan\beta)>0.
\end{align}
The condition $\bx^\top\bB\bc<0$ rewrites 
\begin{align}
	\label{eq:cond2}
\evmax\cos\alpha\cos\beta+\evmin\sin\alpha\sin\beta=\evmin\cos\alpha\cos\beta\left(\frac{\evmax}{\evmin}+\tan\alpha\tan\beta\right)<0.
\end{align}
 We then have that 
\begin{align*}
	f(\bx,\bc)=\frac{(\bx^\top \bB \bc)^2}{(\bx^\top\bB\bx)\bc^\top\bB \bc}&=\frac{(\evmax\cos\alpha\cos\beta+\evmin\sin\alpha\sin\beta)^2}{(\evmax\cos^2\alpha+\evmin\sin^2\alpha)(\evmax\cos^2\beta+\evmin\sin^2\beta)}\\
	&=\frac{(\mu+\tan\alpha\tan\beta)^2}{(\mu+\tan^2\alpha)(\mu+\tan^2\beta)}, \quad \mu:=\frac{\evmax}{\evmin}.
\end{align*}
Next, for $\cos\alpha\neq 0$, $\cos\beta\neq 0$, let us set $\tan\alpha=p$, $\tan\beta=q$ and look for the critical points of 
\begin{align*}
	\mathcal{F}(p,q)=\frac{(\mu+p q)^2}{(\mu+ p^2)(\mu+ q^2)}.
\end{align*}
We are interested in a particular subset of $p,q$ which satisfy \eqref{eq:cond1} and \eqref{eq:cond2}. Remark that $\cos\alpha=0$ implies by \eqref{eq:cond1} that $\sin\alpha\sin\beta>0$ which is incompatible with \eqref{eq:cond2}, and same holds for $\cos\beta=0$. 
Therefore, we limit ourselves to the set
\begin{align}
	\label{eq:setpq}
	(\mu+ pq)(1+pq)<0 \iff pq\in (-\mu, -1).
\end{align} 
The critical points of $\mathcal{F}$ solve
\begin{align*}
	(\mu+ p q)(q(\mu+ p^2)-(\mu+ p q)p)=0, \quad 	(\mu+ p q)(p(\mu+ q^2)-(\mu+ p q)q)=0.
\end{align*}
Thus either $\mu+ p q=0$, which implies that $\mathcal{F}=0$, or $q=p$. The latter set of points does not belong to \eqref{eq:setpq}, and thus it remains to study the behaviour of $\mathcal{F}$ on the boundary of the set $pq\in (-\mu, -1)$. In particular, $\mathcal{F}=0$ for $pq=-\mu$, while for $pq=-1$, 
\begin{align*}
	\mathcal{F}\left(q,-\frac{1}{q}\right)=\frac{(\mu-1)^2}{(\mu+q^{-2})(\mu+q^2)}=\frac{(\mu-1)^2}{\mu^2+1+(q^{-2}+q^2)\mu}\leq \frac{(\mu-1)^2}{(\mu+1)^2}.
\end{align*}
\end{proof}

\section{Proof of  \eqref{eq:esigmat}}
\label{appendix:error}
In the Laplace domain, the PML problem reads 
\begin{align}
	\label{eq:lpl_shift}
	&s^2(1+\frac{\sigma}{s+\gamma})\hat{u}^{\sigma}-\partial_x\left((1+\frac{\sigma}{s+\gamma})^{-1}\partial_x\hat{u}^{\sigma}\right)=0, \quad x\in (0, \radpml+\lpml),\\
&\left. \hat{u}^{\sigma}\right|_{x=0}=\hat{g},\quad \left.\hat{u}^{\sigma}\right|_{x=\radpml+\lpml}=0.
\end{align}
The exact solution to \eqref{eq:lpl_shift} can be computed explicitly:
\begin{align*}
	\hat{u}^{\sigma}(x)=\frac{\hat{g}(s)}{1-\mathrm{e}^{-2 sx_{\sigma}^*}}(\mathrm{e}^{- sx_{\sigma}(s,x)}-\mathrm{e}^{-2 sx_{\sigma}^*+ sx_{\sigma}(s,x)}), \quad x_{\sigma}^*=x_{\sigma}(s,\radpml+L). 
\end{align*}
For $x\in (0,\radpml)$, the error between the exact solution and the above PML solution $\hat{u}(s,x)=\mathrm{e}^{-sx}\hat{g}(s)$ is then given by
\begin{align}
	\label{eq:error_lpl_1d}
	\hat{e}^{\sigma}&=\hat{u}^{\sigma}-\hat{u}=\frac{\hat{g}(s)\mathrm{e}^{-2 sx_{\sigma}^*}}{1-\mathrm{e}^{-2 s x_{\sigma}^*}}\left(\mathrm{e}^{- sx}-\mathrm{e}^{ sx}\right)=\hat{E}(s,x)\hat{g}(s)\left(\mathrm{e}^{- sx}-\mathrm{e}^{ sx}\right), \\
	\nonumber
	\hat{E}(s,x)&=\mathrm{e}^{-2 sx_{\sigma}^*}(1-\mathrm{e}^{-2 sx_{\sigma}^*})^{-1}=\sum\limits_{\ell=1}^{\infty}\mathrm{e}^{-2 s(\radpml+\lpml)\ell}\exp\left(-2\frac{\sigma_c L  s}{s+\gamma}\ell\right).
\end{align}
In the latter expansion we used $|\mathrm{e}^{-2  s x_{\sigma}^*}|<1$, which follows from $\Re \frac{s}{s+\gamma}>0$. 
Using \cite[p.197, formula (16)]{erdelyi}, we find that 
\begin{align*}
	\left(\mathcal{L}^{-1}\exp\left(-2\frac{\sigma_c L   s}{s+\gamma}\ell\right)\right)(t)=\mathrm{e}^{-2\sigma_c L\ell  }\mathrm{e}^{-\gamma t}\left(\left(\frac{2L\gamma\sigma_c \ell}{t}\right)^{1/2}I_1(2\sqrt{2\sigma_c  L\gamma\ell t})+\delta(t)\right),
\end{align*}
where $I_1$ is a modified Bessel function. 
The above allows to rewrite the error explicitly. Let us introduce   
\begin{align*}
	T_{\ell}g(t,x)=g(t- x-2 (\radpml+L)\ell)-g(t+ x-2 (\radpml+L)\ell), 	\quad \alpha_{\ell}=\sqrt{2L\gamma\sigma_c \ell}.
\end{align*}
Remark that because of the causality of $g$, $T_{\ell}g=0$ for 
\begin{align}
	\label{eq:ell_t}
	\ell\geq (t+ \radpml)/(2(\radpml+L)).
\end{align}
Then, for each $t>0$, we write (where the sums below are finite for $t<\infty$ because of the causality of $g$):
\begin{align*}
	\begin{split}
		e^{\sigma}(t,\bx)&=\mathrm{e}^{-\gamma t}\sum\limits_{\ell=1}^{\infty}\mathrm{e}^{-2\sigma_c L\ell}T_{\ell}g\\
		&+\sum\limits_{\ell=1}^{\infty}\mathrm{e}^{-2\sigma_c L\ell}\alpha^{1/2}_{\ell}\int_0^{t}\mathrm{e}^{-\gamma \tau} I_1(2\alpha_{\ell}\tau^{1/2})\tau^{-1/2}T_{\ell}g(t-\tau,x)d\tau. 
	\end{split}
\end{align*}
\section{Implementation: infinite elements}
\label{app:implementation}

For this section we assume that Assumption \ref{assump:piecewise_const} holds and additionally that $\radpml = 1$, i.e., we have
\begin{align*}
	\tilde\sigma(r) = \frac{(r-1)\sigma_c}{r}, \quad \bdpml=\mathbb{S}^1.
\end{align*}
Moreover we merely describe the discretization in the exterior domain $\domext$, since the discretization in $\domint$ is standard.
%
%
As discussed in Section \ref{sec:numerics_ie}, the implementation of the complex scaled infinite elements  necessitates rewriting the bilinear forms in the weak formulation \eqref{eq:rad_time_weak} in polar coordinates. Let us define the following coordinate transformation:
\begin{align*}
	\Psi:(1,\infty)\times\bdpml\to\domext,\quad 
	(r,\hat\vecx)\mapsto r\hat\vecx, \quad \hat\vecx=(\cos\theta,\sin\theta),\quad \theta\in [0,\, 2\pi).
\end{align*}
{Let us additionally introduce the notation $\rad{v}:=v\circ \Psi$. We denote by $\nabla_{\Sigma}\rad{v}$ the tangential gradient of $\rad{v}$ on the surface $\Sigma$ (i.e. $\partial_{\theta}\rad{v}\hat{\bx}_{\perp}$).  
	As seen from \eqref{eq:rad_time_weak}, we need to express only the following two types of the integrals in polar coordinates, for real-valued integrands:
\begin{align*}
	\langle v,v^\dagger\rangle&=\int_{(1,\infty)\times\bdpml}r \rad{v}\; \rad{v}^{\dagger}drd\hat\bx,\qquad
	\langle\nabla v,\bq^\dagger\rangle=\int_{(1,\infty)\times\bdpml} r\partial_{r}\rad{v}\hat \vecx\cdot \rad{\bq}^\dagger+\nabla_{\Sigma}\rad{v}\cdot \Pxp \rad{\bq}^\dagger drd\hat\bx.
\end{align*}
}
Plugging in basis functions (Section cf.\ \ref{sec:numerics_ie}) of the form $\rad{v} (r,\hat\vecx)=\tilde v(r)\hat v(\hat\vecx)$, $v^\dagger_\mathrm{rad} (r,\hat\vecx)=\tilde v^\dagger(r)\hat v^\dagger(\hat\vecx)$, $\bq^\dagger_\mathrm{rad} (r,\hat\vecx)=\tilde q^\dagger(r)\hat \bq^\dagger(\hat\vecx)$ leads to
\begin{align}
	\label{eq:weak_polar}
	\begin{split}
	\langle v,v^\dagger\rangle&=\int_{(1,\infty)}r\tilde v(r)\tilde v^\dagger(r)dr\int_{\bdpml}\hat v(\hat\vecx)\hat v^\dagger(\hat\vecx)d \hat\vecx,\\
	\langle\nabla v,\bq^\dagger\rangle&=\int_{(1,\infty)}r\partial_r\tilde v(r)\tilde q^\dagger(r)dr\int_{\bdpml}\hat v(\hat\vecx)\hat\vecx\cdot\hat\bq^\dagger(\hat\vecx) d\hat \vecx+\int_{(1,\infty)}\tilde v(r)\tilde q^\dagger(r)dr\int_{\bdpml}\nabla_{\Sigma}\hat v(\hat\vecx)\cdot\Pxp\hat\bq^\dagger(\hat\vecx) d\hat \vecx.
	\end{split}
\end{align}
The surface integrals can be evaluated by applying numerical integration on the surface mesh on $\bdpml$ and lead to sparse, well-conditioned matrices if the underlying volume finite element spaces are chosen properly. 
It remains to choose a basis for the radial component of the space given in \eqref{eq:radspace}. We proceed to explain how this is done in the case of standard complex-scaled infinite elements and state the resulting radial discretization matrices. Subsequently we extend this approach to the case of two-scale infinite elements. 

\subsection{Hardy space infinite elements}
Let us consider the following space, which is a special case of \eqref{eq:radspace_onepole}:
\begin{align}
	\label{eq:radspace_c1}
	\scalspextrad^M=\vecspextrad^M=\spa\left\{\exp(-r)p(r): p\in\mathcal P^M\right\}.
\end{align}
{One basis in this space is given by, cf.\ \cite[Section 6.3.2]{WessDiss} and \cite[Section 4.2]{NW22}, }
\begin{align*}
  \left\{r\mapsto \varphi_n(r)=\exp(1-r)L_{n,-1}(2r-2), n=0,\ldots,M\right\},
\end{align*}
where the functions $L_{n,-1}$ are generalized Laguerre polynomials, defined by 
\begin{align}
\label{eq:Ln1}
L_{n,-1}(x)=\left\{
\begin{array}{ll}
1, & n=0, \\
  \sum_{k=1}^n\frac{k}{n}\binom{n}{k}\frac{(-x)^k}{k!}=\frac{x\mathrm{e}^x}{n}\left(\frac{d}{dx}\right)^n\left(x^{n-1}\mathrm{e}^{-x}\right), & n\geq 1.
\end{array}
\right. 
\end{align}
The above defines a basis in the space \eqref{eq:radspace_c1}, because the  degree of the polynomial $L_{n,-1}$ is equal exactly to $n$. 
%
%
This choice of the basis ensures that the radial discretization matrices \eqref{eq:weak_polar} are sparse, see \cite[Section 6.3]{WessDiss}. Moreover, their entries can be computed in a closed form using the orthogonality of Laguerre polynomials $L_{n,-1}$ with respect to the weight $\mathrm{e}^{-x}x^{-1}$. Alternatively the matrices can be computed using numerical quadrature, i.e., weighted Gauss rules with respect to the weight $\exp(-2x)$. 
 
However, since later on we will use a two-scale version of the infinite elements, for which such computations are far from obvious, we choose an alternative approach. The approach\footnote{Historically, this is how Hardy space infinite elements were introduced and analyzed in the first place in \cite{hsm}} is based on passing to the Laplace domain \textit{with respect to the radial variable} $r-1$, and next considering the expressions for the discretization matrices fully in the Laplace domain. This technique can then be generalized to compute the entries of the discretization matrices for the two-scale method in a fairly straightforward manner. Because it is simpler to understand the main ideas for the classical Hardy space infinite elements, we choose this, more pedagogical, way of presenting the method we use in our implementation.
 
Given a function $r\mapsto\varphi(r)$, with $r\in [1, \infty)$, let us introduce its (shifted) Laplace transform (where the index $r$ indicates that we work with the radial variable, rather than with time): $$(\mathcal{L}_r\varphi)(p):=\Phi(p):=\int_0^{\infty}\mathrm{e}^{-p\xi}\varphi(\xi+1)d\xi.$$
Remark that by capital letters we will denote $\mathcal{L}_r$-transformed quantities. 
With the above definition, we can prove the following result. 
\begin{lem}
	\label{lem:lpl_transforms}
	Let $\Re p>-1$. Then the Laplace transforms of $\varphi_n$ are given by the expressions
\begin{align}
	\label{eq:lpl_phi}
	&(\mathcal{L}_r\varphi_0)(p)=\Phi_0(p)=\frac{1}{p+1},\quad
	&(\mathcal{L}_r\varphi_n)(p)=\Phi_n(p)=-\frac{2}{(p+1)^2}\left(\frac{p-1}{p+1}\right)^{n-1}, \quad n\geq 1.
\end{align}
Also, for their derivatives it holds that
\begin{align}
	\label{eq:lpl_der}
	&(\mathcal{L}_r\varphi_0')(p)=-\frac{1}{p+1},\quad
	&(\mathcal{L}_r\varphi_n')(p)=-\frac{2p}{(p+1)^2}\left(\frac{p-1}{p+1}\right)^{n-1}, \quad n\geq 1.
\end{align}
Finally, 
\begin{align}
	\label{eq:lpl_rphi}
	&(\mathcal{L}_r(r\varphi_n))(p)=D_p\Phi_n(p), \quad D_p=1-\partial_p.
\end{align}
\end{lem}
\begin{proof}
\textit{Step 1. Proof of \eqref{eq:lpl_phi}. }First of all, for $n=0$ the result follows by a direct computation. 
	For $n\geq 1$, we rewrite 
	\begin{align}
		\label{eq:phi_n_def}
	\Phi_n(p)=\int_0^{\infty}\mathrm{e}^{-px}\mathrm{e}^{-x}L_{n,-1}(2x)dx.
\end{align}
Replacing $L_{n,-1}$ by its expression \eqref{eq:Ln1} and using the identity
	\begin{align*}
		\int_0^{\infty}\mathrm{e}^{-p x}x^k\mathrm{e}^{-x}dx=k!(p+1)^{-k-1}, \quad k\in \setN,
	\end{align*}
yields the following expression: 
\begin{align*}
	\Phi_n(p)=\sum\limits_{k=1}^n\frac{k}{n}\binom{n}{k}\frac{(-2)^k}{(p+1)^{k+1}}=\frac{1}{p+1}\left.\left(\frac{x}{n}\frac{d}{dx}\sum\limits_{k=0}^n \binom{n}{k}x^k\right)\right|_{x=-\frac{2}{p+1}}=\frac{1}{p+1}\left.x(1+x)^{n-1}\right|_{x=-\frac{2}{p+1}},
\end{align*}
from which the desired expression follows immediately.\\
  \textit{Step 2. Proof of \eqref{eq:lpl_der}. } The expressions for the Laplace transforms of the derivatives are immediate, if we recall additionally that $\varphi_0(1)=1$ and $\varphi_n(1)=0$ for all $n\geq 1$ (cf. \eqref{eq:Ln1}).\\
  \textit{Step 3. Proof of \eqref{eq:lpl_rphi}. }The identity follows by a direct computation: $$\mathcal{L}_r(r\varphi_n)=\int_0^{\infty}\mathrm{e}^{-p\xi}(\xi+1)\varphi_n(\xi+1)d\xi=(1-\frac{d}{dp})\int_0^{\infty}\mathrm{e}^{-p\xi}\varphi_n(\xi+1)d\xi.$$
\end{proof}
In order to relate \eqref{eq:weak_polar} evaluated at radial basis functions $\varphi_n$ to their Laplace domain expressions, we use the following lemma, see also \cite[Lemma A.1]{hsm}. It follows from the Plancherel identity, by continuing $u, \, v$ to causal functions and by making use of the fact that they are real-valued. 
\begin{lem}
		\label{lem:ssq}
Let $u, v\in L^2((1, \infty); 
\setR)$. Then $\int_{1}^{\infty}u(r)v(r)dr=\frac{1}{2\pi i}\int_{i\setR}U(p)V(-p)dp$.
\end{lem}
The above lemma states that evaluating the entries of the mass matrix reduces to computing the integrals in the right-hand side of
\begin{align}
	\label{eq:intrhs}
	\int_{(1,\infty)}\varphi_n(r)\varphi_m(r)dr=\frac{1}{2\pi i}\int_{i\setR}\Phi_n(p)\Phi_{m}(-p)dp.
\end{align}
The main idea is then to rewrite the matrix in the right-hand side by using an alternative basis, in which its representation will be sparse. Even better, this representation will provide us a factorization of the matrix in the right-hand side with the help of tri-diagonal matrices. 
We introduce 
\begin{align*}
  \mathcal{S}_M^{\Phi}:=\spa\{\Phi_n,\, n=0,1,\ldots,M\}.
\end{align*}
Moreover,we define new basis functions $\Psi_n$ and related basis functions $\bbeta_n$ (the latter ones constitute a basis used in the Hardy space methods): 
\begin{align}
	\label{eq:aux_fs}
	\Psi_n(p):=-\frac{2}{p+1}\left(\frac{p-1}{p+1}\right)^n, n\geq 0, \quad \Psi_{-1}(p):=0, \quad
	\boldsymbol{\beta}_n:=\begin{pmatrix}\delta_{n,-1}\\ \Psi_n\end{pmatrix}, \quad n\in \setN\cup\{-1\},
\end{align}
with $\delta_{n,-1}=1$ if $n=-1$ and $0$ otherwise.
We will denote by $\mathcal{S}_M^{\Psi}:=\spa\{\Psi_n, \, n=0,\ldots,{M}\}$ and by $\mathcal{S}_M^{\bbeta}:=\spa\{\bbeta_n, \, n=-1, \ldots, {M}-1\}$.

\begin{rmk}
	The reason to introduce the basis functions $\Psi_n$ is that, as we will see later (Lemma \ref{lemma:e7}), the matrix defined by the right-hand side of \eqref{eq:intrhs} becomes diagonal in this basis.
        This allows to avoid evaluation of the integrals over infinite lines.  The basis $\bbeta_{n}$ is introduced to make a connection to the Hardy-space methods. 
\end{rmk}

The main goal is then to rewrite the 'radial' matrices arising in \eqref{eq:weak_polar} using Lemma \ref{lem:ssq} and the new basis $\mathcal{S}_M^{\Psi}$ (resp. $\mathcal{S}_M^{\bbeta}$). We start by expressing $\Phi_n$ via $\bbeta_n$.
\begin{lem}
  \label{lem:laplace_basis_new}
  Let the operators $\mathcal T_\pm:\mathcal{S}_{M}^{\bbeta}\mapsto \mathcal{S}_M^{\Phi}$ be given by, for $\bbeta=(\beta, B)^{\top}$, 
  \begin{align*}
    (\mathcal T_-\bbeta)(p)&:=\frac{\beta+B(p)}{p+1},&
    (\mathcal T_+\bbeta)(p)&:=\frac{-\beta+pB(p)}{p+1}.
  \end{align*}
  Then, for all $0\leq n\leq M$,
  \begin{align}
    \label{eq:lrphi}
    \mathcal{L}_r\varphi_n=\mathcal T _{-}\boldsymbol{\beta}_{n-1}, \quad  \mathcal{L}_r\varphi_n'=\mathcal T_{+}\boldsymbol{\beta}_{n-1}.
\end{align}
\end{lem}
\begin{proof}
  First of all, remark that indeed, $\operatorname{Range}(\mathcal T_{\pm})\subseteq\mathcal{S}_M^{\Phi}$. This is straightforward from \eqref{eq:aux_fs} and equalities of Lemma \ref{lem:lpl_transforms}. Both identities in \eqref{eq:lrphi} then follow by an explicit computation, using the results of Lemma \ref{lem:lpl_transforms}.
\end{proof}
With the above we also have the following.
\begin{lem}
  \label{lem:operators_tpm}
  $\mathcal{T}_{\pm}(\mathcal{S}_{M}^{\bbeta})\subseteq \mathcal{S}_{M}^{\Psi}$. 
\end{lem}
The above follows from $\mathcal{S}_M^{\Phi}\subseteq\mathcal{S}_{M}^{\Psi}$. We however express $\mathcal{T}_{\pm}\bbeta_n$ in a basis of  $\mathcal{S}_{M}^{\Psi}$ directly, since we will need the corresponding expressions later.
\begin{proof}
Let us prove the result for $\mathcal{T}_{-}$.
  It is straightforward to see that $\mathcal{T}_{-}\bbeta_{-1}=-\frac{1}{2}\Psi_0$. For $n\neq {-1}$, we have that $(\mathcal{T}_{-}\bbeta_n)(p)=\frac{\Psi_n(p)}{p+1}$, and we rewrite 
\begin{align}
	\label{eq:psinp}
\frac{\Psi_n(p)}{p+1}=-\frac{2}{(p+1)^2}\left(\frac{p-1}{p+1}\right)^{n}, \text{ and }\frac{1}{p+1}=-\frac{1}{2}\left(\frac{p-1}{p+1}-1\right).
\end{align}
Thus, for $n\geq 0$, the result follows from the expression:
\begin{align}
	\label{eq:psinp2}
  (\mathcal{T}_{-}\bbeta_n)(p)=-\frac{1}{2}(\Psi_{n+1}(p)-\Psi_n(p)).
\end{align}
  For $\mathcal{T}_+$ we have that $\mathcal{T}_{+}\bbeta_{-1}=\frac{1}{2}\Psi_0$. For $n\neq {-1}$, we have that $(\mathcal{T}_{+}\bbeta_n)(p)=\Psi_n(p)-\frac{\Psi_n(p)}{p+1}=\frac{1}{2}(\Psi_n(p)+\Psi_{n+1}(p))$.  
\end{proof}

These relations, together with Lemma \ref{lem:ssq}, allow us to express the matrices occurring in \eqref{eq:weak_polar} with the help of the 
bilinear form 
\begin{align}
	\label{eq:defq}
	q(U,V):=\frac{1}{2\pi i}\int_{i\setR}U(p)V(-p)dp, 
\end{align}
defined for $U,\, V$ being the Laplace transforms of one of the functions $\varphi_n$ or their derivatives. We then have the following result, whose proof is based on Lemmas \ref{lem:laplace_basis_new} and \ref{lem:ssq}, and thus is left to the reader.
\begin{lem}
	\label{lem:hsm_ie}
	The following identities hold true for all $n, m\geq 0$ (where $q$ is defined in \eqref{eq:defq} and $D_p$ in \eqref{eq:lpl_rphi}):
		\begin{align}
			\label{eq:mx_ex}
			\begin{split}
		\int_{(1,\infty)}\varphi_n\varphi_m &= q(\mathcal T_-(\bbeta_{n-1}),\mathcal T_-(\bbeta_{m-1})),\\
		\int_{(1,\infty)}r\varphi_n(r)\varphi_m(r)dr &= q(D_p\mathcal T_-(\bbeta_{n-1}),\mathcal T_-(\bbeta_{m-1})), \\
		\int_{(1,\infty)}r\varphi_n'(r)\varphi_m(r)dr &= q(D_p\mathcal T_+(\bbeta_{n-1}),\mathcal T_-(\bbeta_{m-1})).
		\end{split}
	\end{align}
\end{lem}
Therefore, instead of directly computing  the matrices in the left-hand side of \eqref{eq:mx_ex}, we can factorize them with the help of the expressions in the right-hand side of \eqref{eq:mx_ex}, by computing the matrix representations of the operators $\mathcal{T}_{\pm}: \mathcal{S}_{M-1}^{\bbeta}\rightarrow \mathcal{S}_M^{\Psi}$, $D_p: \,\mathcal{S}_M^{\Phi}\rightarrow \mathcal{S}_{M+1}^{\Phi}$\footnote{Remark that the multiplication operator $u(r)\mapsto ru(r)$ sends $\mathcal{V}_{\operatorname{rad}}^M$ into $\mathcal{V}_{\operatorname{rad}}^{M+1}$. We can thus show that the shifted derivative $D_p$ maps $ \mathcal{S}_{M}^{\Psi}\rightarrow \mathcal{S}_{M+1}^{\Psi}$.}, and of the bilinear form $q: \, \mathcal{S}_M^{\Psi}\times \mathcal{S}_{M}^{\Psi}\mapsto \mathbb{C}$. Importantly, in this new basis $q$ becomes a diagonal matrix. 
To state the result that follows, let us introduce the $q$-orthogonal projection operator $P_M^{\Psi}: \, \mathcal{S}_{M+1}^{\Psi}\rightarrow \mathcal{S}_M^{\Psi}$, defined by, for $U\in \mathcal{S}_{M+1}^{\Psi}$,
\begin{align*}
	q(P_M^{\Psi}U, V)=q(U,V), \qquad \text{for all }V\in \mathcal{S}_{M}^{\Psi}.
\end{align*}
We then have the following result.
\begin{lem}
	\label{lemma:e7}
	The matrix representations 
	$\mathbb{T}_{\pm}$ of $\mathcal T_\pm: \mathcal{S}_{M-1}^{\bbeta}\rightarrow \mathcal{S}_M^{\Psi}$, $\mathbb{Q}$ of $q:\mathcal{S}_M^{\Psi}\times \mathcal{S}_M^{\Psi}\to\setR$, $\tilde{\mathbb{D}}_p$ of $P_{M}^{\Psi}D_p: \, \mathcal{S}^{\Psi}_M\mapsto \mathcal{S}_{M}^{\Psi}$  are given by:
	\begin{align*}
		\mathbb{T}_\pm&=\frac{1}{2}\begin{pmatrix}
			\pm 1& 1 &&&\\&\pm 1& 1 &&\\&&\ddots&\ddots&\\&&&\pm 1& 1\\ &&&& \pm 1\end{pmatrix},&
			\tilde{\mathbb{D}}_p&=\operatorname{Id}-\frac{1}{2}\begin{pmatrix}-1&1&&&\\1&-3&2&&\\&2&-5&3\\&&\ddots&\ddots&\ddots\\&&&\ddots&\ddots&M\\&&&&M&-2M-1\end{pmatrix},
	& \mathbb Q=2\operatorname{Id}.
	\end{align*}
\end{lem}
\begin{proof}
\textit{Expressions of $\mathbb{T}_{\pm}$. }See Lemma \ref{lem:operators_tpm}. 

\textit{Expression of $\mathbb{Q}$. }By a direct computation it follows that
\begin{align*}
	\mathbb{Q}_{nm}&=\frac{1}{2\pi i}\int_{i\setR}\Psi_n(p)\Psi_m(-p)dp=\frac{1}{2\pi i}\int_{i\setR}\frac{4}{(p+1)(-p+1)}\frac{(p-1)^{n}(-p-1)^{m}}{(p+1)^n(-p+1)^m}dp
	=-\frac{2}{\pi i}\int_{i\setR}\frac{(p-1)^{n-m-1}}{(p+1)^{n-m+1}}dp.
\end{align*}
Let us assume without loss of generality that $n\geq m$. 
For $n\geq m+1$, we can use the Cauchy residue theorem to show that the corresponding integral vanishes, i.e. $\mathbb{Q}_{nm}=0$. For $n=m$, we have that 
\begin{align*}
	\mathbb{Q}_{mm}=4\operatorname{Res}_{p=1}\frac{1}{(p+1)(p-1)}=2.
\end{align*}
\textit{Expression of $\tilde{\mathbb{D}}_p$. } Let us first of all compute 
\begin{align*}
	\frac{d}{dp}\Psi_n(p)=\frac{2(n+1)(p-1)^n}{(p+1)^{n+2}}-\frac{2n(p-1)^{n-1}}{(p+1)^{n+1}}.%
\end{align*}
With \eqref{eq:psinp},  we express $(p+1)^{-1}=-\frac{1}{2}\left(\frac{p-1}{p+1}-1\right)$. This yields
\begin{align*}
	\frac{d}{dp}\Psi_n(p)&=-\frac{(n+1)}{(p+1)}\left(\left(\frac{p-1}{p+1}\right)^{n+1}-\left(\frac{p-1}{p+1}\right)^{n}\right)+\frac{n}{p+1}\left(\left(\frac{p-1}{p+1}\right)^{n}-\left(\frac{p-1}{p+1}\right)^{n-1}\right)\\
	&=\frac{1}{2}\left((n+1)\Psi_{n+1}-(2n+1)\Psi_n+n\Psi_{n-1}\right).
\end{align*}
Finally, as $\mathbb{Q}_{nm}=2\delta_{nm}$, we conclude with the stated form of the matrix expression of $P_{M}^{\Psi}D_p$. 
\end{proof}
	Using the result of the above lemma and \eqref{eq:mx_ex}, we conclude that the entries of the discrete matrices in \eqref{eq:weak_polar} can be written as follows:
\begin{align}
	\begin{split}
		\int_{(1,\infty)}\varphi_n\varphi_m &= \left(\mathbb{T}_{-}^{\top} \mathbb Q\mathbb {T}_-\right)_{n,m},\\
		\int_{(1,\infty)}r\varphi_n(r)\varphi_m(r)dr &= \left(\mathbb{T}_{-}^{\top}\widetilde{\mathbb{D}}_p^{\top} \mathbb Q\mathbb {T}_-\right)_{n,m},\\
		\int_{(1,\infty)}r\varphi_n'(r)\varphi_m(r)dr &= \left(\mathbb{T}_{+}^{\top}\widetilde{\mathbb{D}}_p^{\top} \mathbb Q\mathbb {T}_-\right)_{n,m}.
	\end{split}  
	\label{eq:mat_comp}
\end{align}
\subsection{Two scale Hardy space infinite elements}
Let us consider the two-scale space \eqref{eq:radspace}.
In this case, up to our knowledge, no convenient closed form for the spacial basis functions as for the classical Hardy space infinite elements is known. In the following, we assume that $\eta_0=1$, set $\eta_1=\eta$, and define the basis in the Laplace domain first (following \cite{HallaNannen:18,HallaHohageNannenSchoeberl:16}). The corresponding spacial basis functions can then be found by computing the inverse Laplace transform. This is in general not necessary:  as we are not interested in the solution outside of $\domint$, the explicit representation of these functions is not needed.
We define the two-scale Hardy basis  of \eqref{eq:aux_fs} by 
\begin{align}
	\label{eq:def_psin}
		&\Psi_n^\iepar(p):=-\frac{1+\iepar}{p+\iepar}\left(\frac{p-1}{p+1}\right)^{\lfloor(n+1)/2\rfloor}\left(\frac{p-\iepar}{p+\iepar}\right)^{\lfloor n/2\rfloor}, \quad n\geq 0, \quad \Psi_{-1}^{\iepar}(p):=0, \\ &\bbeta^{\iepar}_n(p):=\left(\begin{matrix}
			\delta_{n,-1}\\
			\Psi_n^{\iepar}(p)
		\end{matrix}\right), \quad n\in \mathbb{N}\cup\{-1\}.
\end{align}
The set of $\{\Psi_n^{\eta}, n\in\mathbb{N}_0\}$ can be shown to form a Riesz basis in an appropriate Hardy space, see \cite{HallaHohageNannenSchoeberl:16}. We will comment more on a particular form of the above basis functions below, see Remark \ref{rem:E1} after Lemma  \ref{lem:E1}. 
We then define the basis functions $\varphi_n^{\eta}$ by making use of an analogue of Lemma \ref{lem:laplace_basis_new}. In particular, first of all we define their shifted Laplace transform $\mathcal{L}_r\varphi_n^{\eta}=\Phi_n^{\eta}$ by a counterpart of \eqref{eq:lrphi}: 
\begin{align*}
  \mathcal L_r\varphi_n^{\iepar}:=\mathcal{T}^{\eta}_{-}\bbeta^{\iepar}_{n-1}, \quad n\geq 0.
\end{align*}
with the operator $\mathcal{T}^{\eta}_{-}$ defined as follows.
Provided the spaces 
\begin{align*}
	\mathcal{S}_M^{\bbeta^{\eta}}&:=\operatorname{span}
	\{\bbeta_n^{\eta}, \quad n=-1, \ldots, M-1\}, \quad \mathcal{S}_M^{\Psi^{\eta}}:=\operatorname{span}\{\Psi^{\eta}_n, \quad n=0, \ldots, M\}, \\
	\mathcal{S}_M^{\Phi^{\eta}}&:=\operatorname{span}\{\Phi^{\eta}_n, \quad n=0, \ldots, M\},
\end{align*}
the operator $\mathcal{T}^{\eta}_{-}: \, \mathcal{S}^{\bbeta^{\eta}}_{M}\rightarrow \mathcal{S}_{M}^{\Phi^{\eta}}$ is defined via
\begin{align*}
  (\mathcal{T}^{\eta}_{-}\bbeta)(p):=\frac{\beta+B(p)}{p+1}, \quad \bbeta=(\beta, B)^{\top}.
\end{align*}
With the definition of the operator $\mathcal{T}_{-}^{\eta}$, it is straightforward to verify that $\varphi_n^{\eta}$ are linearly independent. Moreover, our first result shows that we have indeed constructed the basis of $\scalspextrad(1, \eta, N)$.
\begin{lem}
  It holds that $\operatorname{span}\{\mathcal{L}_r^{-1}\Phi_n^{\eta}, \, n=0, \ldots, 2N+1\}=\scalspextrad(1, \eta, N)$. 
\end{lem}
\begin{proof}
Because $\Phi_{n}^{\eta}$ are linearly independent, it suffices to show that $\operatorname{span}\{\mathcal{L}_r^{-1}\Phi_n^{\eta}, \, n=0, \ldots, 2N+1\}\subset \scalspextrad(1, \eta, N)$. 
We will show how the result is proven for the simplest cases (for $N=0$ and $N=1$), and the result for $N\geq 2$ will follow by the same reasoning. 
Let us first study the case $N=0$. Indeed, in this case $\scalspextrad(1,\eta,0)=\operatorname{span}\{\mathrm{e}^{-p}, \mathrm{e}^{-\eta p}\}$. 
On the other hand, 
\begin{align*}
	\Phi_0^{\eta}(p)=\frac{1}{p+1}, \quad \Phi_1^{\eta}(p)=-\frac{1+\eta}{(p+\eta)(p+1)}, \,\text{ therefore, }\phi_0^{\eta}(r)=\mathrm{e}^{-(r-1)}, \quad \phi_1^{\eta}(r)=-\frac{1+\eta}{1-\eta}\left(\mathrm{e}^{-\eta(r-1)}-\mathrm{e}^{-(r-1)}\right).
\end{align*}
Thus the result of the lemma holds for $\scalspextrad(1, \eta, 0)$. 
Let us now consider the case $N=1$. We rewrite 
\begin{align}
	\label{eq:main_relation}
	&\Phi_2^{\eta}(p)=\frac{p-1}{p+1}\Phi_1^{\eta}(p)=\Phi_1^{\eta}(p)-\frac{2}{p+1}\Phi_1^{\eta}(p),\,
	&\Phi_3^{\eta}(p)=\frac{p-\eta}{p+\eta}\Phi_2^{\eta}(p)=\Phi_2^{\eta}(p)-\frac{2\eta}{p+\eta}\Phi_2^{\eta}(p).
\end{align}
Evidently, $\phi_1^{\eta}\in \scalspextrad(1, \eta, 1)$. Therefore, let us now consider $u(r)=\mathcal{L}_r^{-1}\frac{1}{p+1}\Phi_1^{\eta}(p)$. 
 Applying the shifted inverse Laplace transform, we remark that $u$ solves
\begin{align*}
	\left(\frac{d}{dr}+1\right)u(r+1)&=\phi_{1}^{\eta}(r+1),\, u(1)=0,&\implies u(r+1)&=\int_0^r\mathrm{e}^{-(r-r')}\phi_{1}^{\eta}(r'+1)dr'.
\end{align*}
Now remark that $\phi_{1}^{\eta}(r)=\mathrm{e}^{-r}p_{0}^{(1)}(r)+\mathrm{e}^{-\eta r}p_{0}^{(2)}(r)$, where $p_{0}^{(j)}\in \mathcal{P}_{0}$, and $j=1, 2$ (i.e., constants). Then, it can  be verified by plugging in this explicit form of $\phi_{1}^{\eta}$ in the above expression that $u(r)=\mathrm{e}^{-r}p_{1}^{(1)}(r)+\mathrm{e}^{-\eta r}\tilde{p}_{0}^{(2)}(r)$, where $p_{1}^{(1)}\in \mathcal{P}_{1}$, and $\tilde{p}_0^{(2)}\in \mathcal{P}_0$. Crucially, the polynomial degree in the multiplier of $\mathrm{e}^{-\eta r}$ did not increase after this integration.
The above computation shows that indeed, $\varphi_2^{\eta}\in \scalspextrad(1, \eta, 1)$, and, moreover, writes 
\begin{align*}
\varphi_2^{\eta}(r)=\mathrm{e}^{-r}p_{1}(r)+\mathrm{e}^{-\eta r}{p}_{0}(r), \text{with }p_1\in\mathcal{P}_1, \quad p_0\in \mathcal{P}_0.
\end{align*}
In a similar manner, cf. \eqref{eq:main_relation}, we can show that $\phi_3^{\eta}\in \scalspextrad(1, \eta, 1)$.  
Because \begin{align}
	\label{eq:rec_relation}
	\Phi_{2N+2}^{\eta}(p)=\frac{p-1}{p+1}\Phi_{2N+1}^{\eta}(p), \quad \Phi_{2N+3}^{\eta}(p)=\frac{p-\eta}{p+\eta}\Phi_{2N+2}^{\eta}(p),
\end{align}
 cf. the relation \eqref{eq:main_relation}, the desired result can be proven for $N>1$ by using similar ideas. 
\end{proof}

In order to formulate a counterpart of Lemma \ref{lem:operators_tpm}, it will be more convenient to introduce an auxiliary basis, very similar to $\Psi_n^{\eta}$, namely
\begin{align}
	\label{eq:def_tilde_psi_eta}
	\tilde{\Psi}_{n}^{\eta}(p):=-\frac{1+\eta}{p+1}\left(\frac{p-\eta}{p+\eta}\right)^{\lfloor (n+1)/2\rfloor}\left(\frac{p-1}{p+1}\right)^{\lfloor n/2\rfloor}, \quad n\geq 0.
\end{align}
We define $\tilde{\mathcal{S}}^{\Psi^{\eta}}_{M}:=\operatorname{span}\{\tilde{\Psi}^{\eta}_n, \, n=0, \ldots, M\}$. We then have the following result, which clarifies relations between different spaces.  
\begin{lem}
	\label{lem:E1}
        $\tilde{\mathcal{S}}^{\Psi^{\eta}}_{M}=\mathcal{S}_{M}^{\Phi^{\eta}}$, for all $M\geq 0$. Therefore $\mathcal{T}_{-}^{\eta}(\mathcal{S}^{\bbeta^{\eta}}_{M})=\tilde{\mathcal{S}}^{\Psi^{\eta}}_{M}$.
\end{lem}
\begin{proof}
Remark that $\Phi_0^{\eta}(p)=\frac{1}{p+1}=-\frac{1}{1+\eta}\tilde{\Psi}_0^{\eta}(p)$. For $n\geq 1$, with odd $n=2k+1$, we have  
\begin{align*}
  \Phi_n^{\eta}(p)&=\frac{1}{p+1}\Psi_{2k}^{\eta}(p)=-\frac{(1+\eta)}{(p+\eta)(p+1)}\left(\frac{p-1}{p+1}\right)^{k}\left(\frac{p-\eta}{p+\eta}\right)^{k}\\
	&=-\frac{(1+\eta)}{2\eta(p+1)}\left(1-\frac{p-\eta}{p+\eta}\right)\left(\frac{p-1}{p+1}\right)^{k}\left(\frac{p-\eta}{p+\eta}\right)^{k}=\frac{1}{2\eta}\left(\tilde{\Psi}_{2k}^{\eta}(p)-\tilde{\Psi}_{2k+1}^{\eta}(p)\right).
\end{align*}
In a similar manner, for even $n=2(k+1)$, we have that
	\begin{align*}
          \Phi_n^{\eta}(p)&=\frac{1}{p+1}\Psi_{2k+1}^{\eta}(p)=-\frac{(1+\eta)}{(p+\eta)(p+1)}\left(\frac{p-1}{p+1}\right)^{k+1}\left(\frac{p-\eta}{p+\eta}\right)^{k}\\
		&=\frac{(p-1)}{(p+\eta)(p+1)}\frac{-(1+\eta)}{p+1}\left(\frac{p-1}{p+1}\right)^{k}\left(\frac{p-\eta}{p+\eta}\right)^{k}.
	\end{align*}
Next, let us develop 
\begin{align*}
	\frac{p-1}{(p+\iepar)(p+1)}=
	\frac{\iepar-1}{2\iepar(\iepar+1)}+\frac{p-\iepar}{2\iepar(p+\iepar)}-\frac{(p-1)(p-\iepar)}{(1+\iepar)(p+1)(p+\iepar)}.
\end{align*}
With the above, we obtain 
\begin{align*}
		\Phi_{2k+2}^{\eta}&=\frac{1}{2\iepar}\left(\frac{\iepar-1}{\iepar+1}\tilde{\Psi}_{2k}^{\eta}+\tilde{\Psi}_{2k+1}^{\eta}-\frac{2\eta}{1+\eta}\tilde{\Psi}_{2k+2}^{\eta}\right). 
\end{align*}
\end{proof}
\begin{rmk}
	\label{rem:E1}
        The above lemma partially justifies the particular definition of functions $\Psi_n^{\eta}(p)$ \eqref{eq:def_psin} and $\tilde{\Psi}_n^{\eta}(p)$ \eqref{eq:def_tilde_psi_eta}. Remark that  $\tilde{\Psi}_n^{\eta}$ are defined like  $\Psi_n^{\eta}$ with $\eta$ and $1$ being interchanged. The functions $\tilde{\Psi}_n^{\eta}$ are introduced in order to ensure  $\tilde{\mathcal S}^{\Psi^{\eta}}_M={\mathcal S}_M^{\Phi^{\eta}}$, which does not hold true for $\tilde{\mathcal S}^{\Psi^{\eta}}_M$ replaced by ${\mathcal S}^{\Psi^{\eta}}_M$.
\end{rmk}
To formulate the next result, let us introduce the operator $\mathcal{T}_{+}^{\eta}: \, \mathcal{S}_M^{\bbeta^{\eta}}\rightarrow \tilde{\mathcal{S}}_M^{\Psi^{\eta}}$, defined by 
\begin{align*}
  (\mathcal T_+^{\eta}\bbeta)(p)&:=\frac{-\beta+pB(p)}{p+1}, \quad \bbeta=(\beta, B)^{\top}.
\end{align*}
It is straightforward to verify that it holds that 
\begin{align*}
  \mathcal{L}_r\left(\varphi_n^{\eta}\right)'=\mathcal{T}^{\eta}_{+}(\bbeta_{n-1}), \quad n\geq 0. 
\end{align*}
%
%
I.e., the same relations for the integrals as in Lemma \ref{lem:hsm_ie} hold and the discrete matrices for the corresponding operators and the bilinear form $q$ can be obtained by similar computations as above as given in the following lemmas. Since the computations are very tedious and technical we omit them at this point.
\begin{lem}
	\label{lem:Topmats_hsm_2p}
        Let $M+1=2N$. Then the operators $\mathcal T_\pm^{\iepar}: \mathcal{S}_M^{\bbeta^\iepar}\rightarrow \tilde{\mathcal{S}}_M^{\Psi^\iepar}$ can be represented by the matrices
	\begin{align*}
          \mathbb T^\iepar_-&=-\frac{1}{2\iepar}\operatorname{Id} +\frac{1}{2\iepar}\mathcal T^\iepar,& \mathbb T^\iepar_+&=\frac{1}{2}\operatorname{Id} +\frac{1}{2}\mathcal T^\iepar,
        \end{align*}
        with
        \begin{align*}
          \mathcal T^\iepar&:=
          \begin{pmatrix}
                T&T^U\\
                &\ddots&\ddots\\
                &&T&T^U\\
                &&&T
		\end{pmatrix},
                &T:=\begin{pmatrix}\frac{1-\iepar}{1+\iepar}&1\\0&0\end{pmatrix},&
                &T^U:=\begin{pmatrix}\frac{\iepar-1}{\iepar+1}&0\\1&0\end{pmatrix},
	\end{align*}
        Moreover, the matrix representation of the bilinear form $q:\tilde {\mathcal S}_M^{\Psi^\iepar}\times \tilde {\mathcal S}_M^{\Psi^\iepar}\to\setR$ is given by
	\begin{align*}
		 \mathbb{Q}^\iepar&=\frac{(1+\iepar)^2}{2}\begin{pmatrix}
                   Q\\&\ddots\\&&Q
                   \end{pmatrix},
                   &Q=\begin{pmatrix}
                        1&\frac{1-\iepar}{1+\iepar}\\
			\frac{1-\iepar}{1+\iepar}&1
                 \end{pmatrix}.
        \end{align*}
\end{lem}

Similar to the classical space we use a projection $P_{M}^{\tilde\Psi^\eta}:\mathcal S^{\tilde\Psi^\eta}_{M+2}\to\mathcal S^{\tilde\Psi^\eta}_{M}$ which is given by truncating the expansion with respect to the basis functions $\tilde\Psi^\eta_n$ to $n\leq M$. Then the corresponding square matrix of the differential operator is given by the following lemma.
\begin{lem}
        The matrix of the differential operator ${P}_N^{\tilde{\Psi}}D_p:\tilde{\mathcal{S}}_N^{\Psi^{\iepar}}\to\tilde{\mathcal{S}}_N^{\Psi^{\iepar}}$ (cf. \eqref{eq:lpl_rphi}) is given by
	\begin{align*}
		\widetilde{\mathbb D}^\iepar&=\operatorname{Id}-\frac{1}{2}\begin{pmatrix}D_0&D_1^U\\
			D^L_1&D_1&D^U_2\\
			&\ddots&\ddots&\ddots\\
                        &&D^L_{N-2}&D_{N-2}&D^U_{N-1}\\
                        &&&D^L_{N-1}&D_{N-1}
                \end{pmatrix},
        \end{align*}
        with
        \begin{align*}
                D_l&=\begin{pmatrix}-\frac{1}{1+\iepar}&\frac{1}{2\iepar}\\\frac{1}{2\iepar}&-\frac{1+4\iepar+\iepar^2}{2\iepar(1+\iepar)}\end{pmatrix}+l\begin{pmatrix}\frac{1+6\iepar+\iepar^2}{2\iepar(1+\iepar)}&\frac{1+\iepar}{2\iepar}\\\frac{1+\iepar}{2\iepar}&-\frac{1+6\iepar+\iepar^2}{2\iepar(1+\iepar)}\end{pmatrix},&
                  D^U_l&=l\begin{pmatrix}\frac{\iepar-1}{2\iepar(1+\iepar)}&0\\\frac{1+\iepar}{2\iepar}&\frac{1-\iepar}{2(1+\iepar)}\end{pmatrix},&
                    D^L_{l}&=l\begin{pmatrix}\frac{1-\iepar}{2(1+\iepar)}&\frac{1+\iepar}{2\iepar}\\0&\frac{\iepar-1}{2\iepar(1+\iepar)}\end{pmatrix}.
	\end{align*}
\end{lem}

\end{document}